\documentclass[preprint,11pt]{elsarticle}
\usepackage[toc,page]{appendix}
\usepackage{amsmath,amsthm,amsfonts,mathtools,bm,hyperref}
\usepackage{multirow}
\usepackage{tikz,pgfplots}
\usetikzlibrary{patterns}
\definecolor{BurntOrange}{rgb}{0.8, 0.33, 0.0}
\usepackage{subcaption}
\pgfplotsset{compat=newest}
\usepackage{float,indentfirst,algorithm}
\usepackage{pgfplots}
\usepackage{booktabs}
\usepackage[margin=1in]{geometry}

\theoremstyle{plain}

\newtheorem{remark}{Remark}
\usepackage[capitalise]{cleveref}
\usepackage{appendix}
\usepackage{multicol}
\usepackage{nomencl}
\usepackage{ifthen}

\usepackage{todonotes}
\usepackage{url}

\newcommand{\Mcal}{\ensuremath{\mathcal{M}}}

\newcommand{\ChoL}{\mathsf{L}}
\newcommand{\Loss}{\ensuremath{\mathcal{L}}}

\newcommand{\bt}[0]{\bm{\theta}}

\newcommand{\bu}{\mathbf{u}}

\newcommand{\bP}[0]{\mathbf{P}}
\newcommand{\bF}[0]{\mathbf{f}}

\newcommand{\bM}[0]{\mathbf{M}}

\makenomenclature
\usepackage{etoolbox}
\renewcommand\nomgroup[1]{%
  \item[\bfseries
  \ifstrequal{#1}{A}{Symbols}{%
  \ifstrequal{#1}{B}{Superscripts}{%
  \ifstrequal{#1}{C}{Subscripts}{
  \ifstrequal{#1}{D}{Operators}{}}}}%
]}

\makeatletter
\def\ps@pprintTitle{%
   \let\@oddhead\@empty
   \let\@evenhead\@empty
   \def\@oddfoot{\reset@font\hfil\thepage\hfil}
   \let\@evenfoot\@oddfoot
}
\makeatother

\begin{document}

\begin{abstract}
We present the Cholesky-factored symmetric positive definite neural network (SPD-NN) for modeling constitutive relations in dynamical equations.  Instead of directly predicting the stress, the SPD-NN trains a neural network to predict the Cholesky factor of  a tangent stiffness matrix, based on which the stress is calculated in the incremental form.  
As a result of the special structure, SPD-NN weakly imposes convexity on the strain energy function, satisfies time consistency for path-dependent materials, and therefore improves numerical stability, especially when the SPD-NN is used in finite element simulations. Depending on the types of available data, we propose two training methods, namely direct training for strain and stress pairs and indirect training for loads and displacement pairs. 
We demonstrate the effectiveness of SPD-NN on hyperelastic, elasto-plastic, and multiscale fiber-reinforced plate problems from solid mechanics. The generality and robustness of the SPD-NN make it a promising tool for a wide range of constitutive modeling applications. 
\end{abstract}

\begin{keyword}
Neural Networks, Plasticity, Hyperelasticity,  Finite Element Method, Multiscale Homogenization
\end{keyword}

\begin{frontmatter}

\title{Learning Constitutive Relations using  Symmetric Positive Definite Neural Networks}
\author[rvt1]{Kailai~Xu\corref{cor1}}
\ead{kailaix@stanford.edu}

\author[rvt1]{Daniel~Z.~Huang\corref{cor1}}
\ead{zhengyuh@stanford.edu}

\author[rvt1,rvt2]{Eric~Darve}
\ead{darve@stanford.edu}

\address[rvt1]{Institute for Computational and Mathematical Engineering,
               Stanford University, Stanford, CA, 94305}
\address[rvt2]{Mechanical Engineering, Stanford University, Stanford, CA, 94305}
\cortext[cor1]{Both authors contributed equally to this work.}

\end{frontmatter}


\nomenclature[A]{$\mathsf{H}$}{tangent stiffness matrix}
\nomenclature[A]{$\mathsf{C}$}{linear tangent stiffness matrix}
\nomenclature[A]{$\mathcal{M}$}{constitutive model}
\nomenclature[A]{$\mathsf{M}_{\bt}$}{neural network based constitutive model}
\nomenclature[A]{$\sigma$}{stress}
\nomenclature[A]{$\bm{\sigma}$}{stress tensor or its Voigt notation}
\nomenclature[A]{$\bm{S}$}{second Piola-Kirchhoff  stress  tensor}
\nomenclature[A]{$w$}{strain-energy density function}
\nomenclature[A]{$\bm{F}$}{deformation gradient tensor}
\nomenclature[A]{$J$}{determinant of $\bm{F}$}
\nomenclature[A]{$\bm{C}$}{right Cauchy-Green stretch tensor}
\nomenclature[A]{$I_1, I_2, I_3$}{three  scalar  invariants of right Cauchy-Green stretch tensor}
\nomenclature[A]{$\mathbb{I}$}{identity matrix}
\nomenclature[A]{$\epsilon$}{strain}
\nomenclature[A]{$\bm{\epsilon}$}{strain tensor or its Voigt notation}
\nomenclature[A]{$\bt$}{weights and biases of the neural network}
\nomenclature[A]{$\mathsf{L}$}{lower triangular Cholesky decomposition of the tangent stiffness matrix}
\nomenclature[A]{$\mathcal{I}(t)$}{all other quantities related to material states}
\nomenclature[A]{$\rho$}{density}
\nomenclature[A]{$\bu$}{(displacement) displacement vector}
\nomenclature[A]{$\mathbf{b}$}{body force vector}
\nomenclature[A]{$\bM$}{mass matrix}
\nomenclature[A]{$\bP$}{discrete internal force vector}

\nomenclature[A]{$p$}{load parameters}
\nomenclature[A]{$d$}{transition function parameter}
\nomenclature[A]{$D$}{transition function}
\nomenclature[A]{$\bF$}{discrete external force vector}
\nomenclature[A]{$f$}{yield function}
\nomenclature[A]{$\bm\epsilon^{e}$}{elastic strain tensor}
\nomenclature[A]{$\bm\epsilon^{p}$}{plastic strain tensor}
\nomenclature[A]{$\bm x$}{coordinate}
\nomenclature[A]{$n$}{number of time step}
\nomenclature[A]{$N$}{number of data}
\nomenclature[A]{$\sigma_Y$}{yield stress}
\nomenclature[A]{$\tilde{\sigma}_Y$}{estimated yield stress}
\nomenclature[A]{$n_x, n_y, n_z$}{number of elements in each direction}
\nomenclature[A]{$T$}{total time}
\nomenclature[A]{$\alpha$}{hardening variable}
\nomenclature[A]{$E$}{Young's modulus}
\nomenclature[A]{$K$}{plastic modulus}
\nomenclature[A]{$L_x, L_y,L_z$}{length}
\nomenclature[A]{$A$}{cross section area}
\nomenclature[A]{$t$}{current time}
\nomenclature[A]{$c_1, c_2$}{Rivlin-Saunders material parameters}
\nomenclature[A]{$\bar{\bm{t}}$}{traction}
\nomenclature[A]{$\dot{\lambda}$}{plastic  rate parameter}
\nomenclature[A]{$\nu$}{Poisson's ratio}
\nomenclature[A]{$r$}{fiber radius}
\nomenclature[A]{$\lambda_J$}{Lagrange multiplier}
\nomenclature[B]{$\Box_{\bm{\theta}}$}{neural network related quantity}
\nomenclature[B]{$\Box^{i}$}{$i$-th time step}
\nomenclature[C]{$\Box_{j}$}{$j$-th data}
\nomenclature[D]{$\Delta\Box$}{incremental for one time step}
\nomenclature[D]{$\nabla\Box$}{gradient operator}
\nomenclature[D]{$\texttt{div}\Box$}{divergence}
\nomenclature[D]{$\Box^T$}{transpose}
\nomenclature[D]{$\dot\Box,\ddot\Box$}{time derivatives}
\begin{multicols}{2}
\printnomenclature
\end{multicols}

\section{Introduction}

Material modeling aims to construct constitutive models to describe the relationship between strain and stress, in which the relationship may be hysteresis. 
The constitutive relations can be derived from microscopic interactions between multiscale structures or between atoms. However, the first-principles simulations, which resolve all these interactions, remain prohibitively expensive.
Because of the computational difficulty, constitutive models are traditionally constructed  with simplified assumptions and empirical relations. Parameters in the models are then calibrated on limited and coarse-scale tensile test data. These constitutive models lead to affordable simulations,  and thus the models are very important for large-scale engineering and scientific applications. However, constructing these models is challenging because the models are usually mappings between high dimensional spaces and are different on a case by case basis.

Since the development of deep learning techniques, the deep neural network (NN) emerges as a promising alternative for constitutive modeling. 
For example, pioneering work has demonstrated the feasibility of constitutive modeling using neural networks in a wide variety of applications, such as Ghaboussi et al.~\cite{ghaboussi1991knowledge} for modeling concrete, Ellis et al.~\cite{ellis1995stress} for modeling  sands,  Shen et al.~\cite{shen2005finite} and Liang et al.~\cite{liang2008neural} for modeling hyper-elastic materials, and Furukawa et. al.~\cite{furukawa1998implicit} for modeling viscoelastic materials. 
Recently, recurrent neural networks~(RNNs), effective for history-dependent phenomena, have been applied to model multiscale multi-permeability poroplasticity~\cite{wang2018multiscale}, multiscale-plasticity~\cite{mozaffar2019deep}, and multiscale one-dimenensional bars~\cite{ghavamian2019accelerating}. 
Additionally, neural networks are used to address more complex material behaviors, such as microcracking, brittle fracture, and crack propagation~\cite{asteris2017anisotropic, goswami2019transfer, gajewski2014sensitivity,liu2020neural}. 
Moreover, attempts for new material designs are made by leveraging neural networks with microscopic structure parameters incorporated in the inputs~\cite{le2015computational, bessa2017framework, chen2020generative,zhang2019machine}. 
The cited literature shows the potential capability of neural networks to \textit{represent} complex constitutive relations. In this paper, we are more interested in using the NN-based constitutive model to \textit{predict} material behavior, which involves embedding the neural network into finite element solvers.

When solving the conservation equations embedded with NN-based constitutive models that have no constraints on the neural networks, numerical instabilities are observed~(see \cref{SEC:TRUSS} for an example). We think one reason is the violation of the basic constraints related  to these constitutive models. The idea of adding physics-related constraints to the neural network is not new. For example, 
Ling et al.~\cite{ling2016machine} enforced isotropicity and cubic symmetry on the predicted strain-energy of the crystal elastic materials through building a basis of invariant inputs.
Liu et al.~\cite{liu2019deep,liu2019exploring} designed deep material networks, as composition of simple building blocks inspired from the two-phase linearly elastic model. Heider et al.~\cite{heider2020invariance} introduced coordinate-free invariant metrics for the objective function (the so-called loss function) to model anisotropic elastoplastic materials. 
However, despite the improved accuracy on the validation set by incorporating these constraints, the performances of the resulting conservation equations equipped with these physics constrained neural networks are not studied.

In the present work, we focus on the numerical stability aspect of the resulting hybrid model---the conservation equations embedded with NN-based constitutive models. 
We propose a novel neural network architecture with customized output layers to fulfill two objectives. 
\begin{itemize}
    \item The first objective is related the convexity of the strain-energy, which plays significant role in the numerical stability. To this end, we propose to predict the tangent stiffness matrix and enforce  it to be symmetric positive definite. The symmetric positive definiteness guarantees the weak convexity of the strain-energy. 
    \item Another objective is the time consistency, related to path-dependent materials, such as elasto-plastic material. 
\end{itemize}
To achieve these objectives, instead of training a neural network to identify a nonlinear map directly between strain and stress, we train a neural network (with weights and biases $\bt$) that maps the strain and other related quantities to a \textbf{lower triangular matrix} $\ChoL_{\bt}$ (Cholesky factor). We construct the constitutive model in the following incremental form:
\begin{equation}\label{equ:delta_sigma}
    \Delta \bm{\sigma} =\ChoL_{\bt}\ChoL_{\bt}^T \Delta \bm{\epsilon}
\end{equation}
where $\Delta \bm{\epsilon}$ and $\Delta \bm{\sigma}$ are the incremental strain and stress in Voigt notation\footnote{For brevity, $\bm{\epsilon}$ and $\bm{\sigma}$ represent both tensor form and vector form in Voigt notation depending on the context.}.
In our formulation, the tangent stiffness matrix $\ChoL_{\bt}\ChoL_{\bt}^T$ is automatically symmetric positive semidefinite and the incremental form leads to time consistency.  

In the present work, results based on neural networks trained on both direct input-output data and indirect full-field data~\cite{grediac2006virtual, huang2019predictive,yang2019structural} are presented.
The robustness of our approach is demonstrated in different numerical applications, including a hyperelastic material, an elasto-plastic material, and a multiscale material.
We have developed a software library that seamlessly combines traditional finite element methods and neural networks.  The code is accessible online:
\begin{center}
\url{https://github.com/kailaix/NNFEM.jl}
\end{center}

The remainder of this paper is organized as follows. 
We first introduce the background in \Cref{SEC:BACKGROUD}, including the governing equations, different classical constitutive relations and their associated constraints. 
In \Cref{SEC:METHODS}, we present our constraint-embedded neural network architecture---SPD-NN, and the training procedures. 
Finally, we apply the learned NN-based constitutive relations for several solid mechanics problems, including a one-dimensional truss coupon, and several two-dimensional thin plate problems with hyperelastic material, elasto-plastic material and  fiber-reinforced multiscale material. Several issues related to neural networks is  discussed in \Cref{SEC:DISCUSS_NN}. We discuss a possible generalization of the approach in \Cref{SEC:CONCLUSION}.

\section{Background}
\label{SEC:BACKGROUD}

\subsection{Governing Equations}
The governing equation of a solid undergoing infinitesimal deformations can be written as
\begin{equation}
\begin{aligned}
   \rho \ddot{\bu}  &= \texttt{div}~\bm{\sigma} + \rho  \mathbf{b} &&\textrm{ in } \Omega\\
    \bu &= \bar{\bu} &&\textrm{ on } \Gamma_{u}\\
    \bm{\sigma} \bm{n} &= \bar{\bm{t}} &&\textrm{ on } \Gamma_{t}\\
\end{aligned}
\label{EQ:LINEAR_MOMENTUM}
\end{equation}
where $\rho$ is the mass density, $u$ is the displacement vector, $\bm{\sigma}$ is the stress tensor, and $\rho \mathbf{b}$ is the body force vector; $\Omega$ denotes the  computational domain.
The prescribed displacement $\bar{u}$ and the surface traction $\bar{\bm{t}}$ are imposed on the domain boundaries $\Gamma_{u}$ and $\Gamma_{t}$ with the outward unit normal $\bm{n}$, where $\Gamma_u \cap \Gamma_t = \emptyset$ and $\Gamma_u \cup \Gamma_t = \partial \Omega$.

To solve for the displacement $u$ from \cref{EQ:LINEAR_MOMENTUM}, we also need the constitutive relations, which maps the deformation history of the structure to the stress:
\begin{equation}\label{EQ:INF_CONSTITUTIVE_LAW}
    \bm{\sigma}(t) = \Mcal(\bm{\epsilon}(t), \mathcal{I}(t))
\end{equation}
Here  $\bm{\epsilon}(t)$ is the strain tensor at time $t$ related to the displacement vector and $\mathcal{I}(t)$ denotes all other quantities related to material states during time $\tau = [0,t)$, such as $\bm{\epsilon}(\tau)$, $\bm{\sigma}(\tau)$, etc. The infinitesimal strain tensor is
\begin{equation*}
    \bm{\epsilon} = \bm{\epsilon}(\bu) = \frac{1}{2}[\nabla \bu + (\nabla \bu)^T]
\end{equation*}

When we apply the finite element method to solve \cref{EQ:LINEAR_MOMENTUM} numerically, we have the following semi-discrete equation at time $t$
 \begin{equation}
 \label{EQ:LINEAR_MOMENTUM_SEMI}
 \bM \ddot{\bu}  + \bP(\bu, \Mcal(\bm{\epsilon}(\bu), \mathcal{I})) = \bF(\bu, x, p)
 \end{equation}
where we use the same notation $\bu$ to denote the spatial discretization of the displacement vector $\bu$ in \cref{EQ:LINEAR_MOMENTUM}, $\bM$ is the discrete mass matrix, $\bP$ and  $\bF$ are the discrete internal and external force vectors, $x$ is the coordinate vector, and $p$ is the parameter vector of external loads. We adopt the generalized $\alpha$-method~\cite{chung1993time} with $\alpha_m=-1$ and $\alpha_f = 0$ for temporal discretization of \cref{EQ:LINEAR_MOMENTUM_SEMI}. This generalized $\alpha$-method allows for dissipating high frequency energy  to damp high frequency modes, which is crucial for the robustness of the numerical solver when an approximate constitutive relation is used. 

\begin{remark}
\label{SEC:FINITE_DEFORMATION}
For structure undergoing finite or large deformations, the finite strain tensor reads
\begin{equation}
    \bm{\epsilon} = \bm{\epsilon}(\bu) = \frac{1}{2}[\nabla \bu + (\nabla \bu)^T + (\nabla \bu)^T  \nabla \bu]
\end{equation}
$\bm\sigma$ in \cref{EQ:LINEAR_MOMENTUM} represents the first Piola-Kirchhoff stress tensor, the constitutive relation generally relates the finite strain tensor with the symmetric second Piola-Kirchhoff stress tensor 
\begin{equation}
\label{EQ:FIN_CONSTITUTIVE_LAW}
\bm S(t) = \Mcal(\bm{\epsilon}(t), \mathcal{I}(t))
\end{equation}
The first Piola-Kirchhoff tensor and the second Piola-Kirchhoff tensor are related by 
$$\bm\sigma = \bm F  \bm S$$
here $\bm F = \nabla \bu + \mathbb{I}$ is the deformation gradient tensor, where $\mathbb{I}$ is the identity matrix. 

Our method works for both infinitesimal and finite deformations. 
Numerical examples of structures undergoing finite deformations are reported in \cref{SEC:TRUSS} and \cref{SEC:PLATE_HYPERELASTICITY}, and  examples using infinitesimal deformations are reported in other numerical examples.  

\end{remark}

\subsection{Constitutive Relations and Associated Constraints}
\label{SEC:CONSTRAINTS}
\Cref{EQ:INF_CONSTITUTIVE_LAW} represents one possible form of constitutive relations. Even when restricted to this specific form, the constitutive relations can describe a great variety of material properties. Nevertheless, most of the constitutive relations share extra constraints (i.e., objectivity, convex strain energy function, and second law of thermodynamics) that play significant roles in material stability and numerical stability. 
In what follows, several stability- and consistency-related constraints of constitutive relations are discussed. These constraints should be considered in any data-driven constitutive models.

\subsubsection{Hyperelastic Materials}
The constitutive relations of  hyperelastic materials are path-independent and are related to the strain-energy density function $\omega(\bm{\epsilon})$, 
\begin{equation}
\label{EQ:CONSTITUTIVE_LAW_ELASTIC}
    \bm{\bm{\sigma}} = \frac{\partial \omega(\bm{\epsilon})}{\partial \bm{\epsilon}}
\end{equation}
In general, the strain-energy density function is assumed to be convex, namely, the tangent stiffness matrix $\frac{\partial^2 \omega(\bm{\epsilon})}{\partial^2 \bm{\epsilon}}
$ is symmetric positive definite~(SPD). The assumption in one-dimensional is equivalent to that the strain-stress is monotonically increasing\footnote{Exceptions exist, i.e., strain-softening~\cite{bazant1984continuum}, which is beyond the scope of the present work.}. The convexity is crucial to both the stability of the material and also the numerical scheme.

\subsubsection{Elasto-plastic Materials}

The constitutive relations of elasto-plastic materials are rate-independent but path-dependent, and feature transition between elastic and plastic behaviors. Namely, under loading, the material behaves elastically until the initial yield stress~$\sigma_Y$ is attained, and then undergoes permanent irreversible plastic deformations with further loading. The onset and continuance of plastic deformation is governed by a yield function
\begin{equation}
    f(\bm{\sigma}) \leq 0
\end{equation}
For plastic deformation, the stress state must remain on the yield surface $f = 0$; and for elastic deformation, the yield function satisfies $f < 0$. 

The strain is assumed to be additively decomposed into elastic $\bm{\epsilon}^e$ and plastic $\bm{\epsilon}^p$ parts, as follows, 
\begin{equation}
    \bm{\epsilon} = \bm{\epsilon}^e + \bm{\epsilon}^p
\end{equation}
The constitutive relation relates the stress and the elastic part of the strain:
\begin{equation}
\label{EQ:CONSTITUTIVE_LAW_PLASTICITY}
    \bm{\sigma} = \mathsf{C} \bm{\epsilon}^e
\end{equation}
here $\mathsf{C}$ is the tangent stiffness tensor.
The plastic strain rate is given by a flow rule, i.e., the associative flow rule as follows,
\begin{equation*}
    \dot{\bm{\epsilon}}^p = \dot{\lambda} \frac{\partial f}{\partial \bm{\sigma}}
\end{equation*}
where $\dot{\lambda}$ is called the plastic rate parameter or the consistency parameter, which is non-zero only if $f = 0$.

Naturally, the constitutive relation is written in the rate form
\begin{equation}
\label{EQ:INF_LAW_PLASTICITY_INCRE}
    \dot{\bm{\sigma}} = \mathsf{H}\dot{\bm{\epsilon}}=
    \begin{dcases*}
    \mathsf{C} \dot{\bm{\epsilon}} & \text{ if } $f < 0$\\
    \Big[\mathsf{C} - \frac
    {(\mathsf{C} \, \frac{\partial f}{\partial \bm{\sigma}})  (\mathsf{C} \, \frac{\partial f}
    {\partial \bm{\sigma}})^T}{\left(\frac{\partial f}{\partial \bm{\sigma}}\right)^T  \mathsf{C} \, \frac{\partial f}{\partial \bm{\sigma}}} \Big]      
    \dot{\bm{\epsilon}} & \text{ if } $f = 0$
    \end{dcases*}
\end{equation}
It is worth mentioning the tangent stiffness matrix $\mathsf{H}$ is symmetric positive semidefinite, and when strain hardening is considered, it becomes SPD. 
Moreover, the form in \cref{EQ:INF_LAW_PLASTICITY_INCRE} implies time consistency, i.e., as $\Delta\bm\epsilon\rightarrow 0$, $\Delta{\bm{\sigma}}\rightarrow 0$, and also rate-independent, i.e., the tangent stiffness matrix is independent of the strain or stress rate.

\subsubsection{Associated Constraints}
Based on the discussion above, the following properties are crucial to be incorporated in data-driven constitutive models
\begin{enumerate}
    \item[(1)] Symmetry positive definiteness of the tangent stiffness matrix~(i.e., strain energy convexity), which leads to non-singular stiff matrix. We also found that the SPD property is crucial for numerical stability.
    \item[(2)] Time consistency, which is formulated as follows,
    \begin{equation}
    \displaystyle \lim_{\Delta \bm{\epsilon} \rightarrow 0} \Delta\bm{\sigma} = 0
    \end{equation}
    It is crucial for the convergence of the numerical approximation when $\Delta t \xrightarrow{} 0$.
\end{enumerate}

\section{Methodology}
\label{SEC:METHODS}

In this section, we describe our method for learning the constitutive relations~\cref{EQ:INF_CONSTITUTIVE_LAW,EQ:FIN_CONSTITUTIVE_LAW}. Our discussion will be divided into two parts: \begin{enumerate}
    \item the neural network architecture for approximating the mapping between the strain and the stress;
    \item the direct and indirect training methods based on the types of available data.
\end{enumerate}   

\subsection{Neural Network Architectures}
\label{SEC:CONSTRAINTS-EMBEDDED-NN}

We note that the neural network has already been used to approximate the constitutive relations. The neural network architectures in many literatures output stress directly~($\bm{\sigma}$-NN) with the following form (recall that $\mathcal{I}$ stands for other relevant information up to the current time)
\begin{equation}
\label{EQ:NN1}
\mbox{$\bm{\sigma}$-NN:}\qquad \bm{\sigma} =  \mathsf{NN}_{\bt}(\bm{\epsilon},\mathcal{I})
\end{equation}
or output stress increment directly~($\Delta\bm{\sigma}$-NN):
\begin{equation}
\label{EQ:NN2}
 \mbox{$\Delta\bm{\sigma}$-NN:}\qquad  \Delta\bm{\sigma} =  \mathsf{NN}_{\bt}(\bm{\epsilon},\mathcal{I}) 
\end{equation}

Those architectures are suitable for learning the constitutive relations when the strain and the stress data are available. However, these  strain-stress relations expressed by \Cref{EQ:NN1,EQ:NN2} do not satisfy certain physical constraints and thus may break numerical solvers~(see \cref{SEC:TRUSS}), when we plug them into a numerical solver. 

We propose an alternative architecture~(SPD-NN) based on the incremental form, 
\begin{equation}\label{EQ:CHOL}
  \begin{aligned}
     \Delta\bm{\sigma} &= \mathsf{H}_{\bt}  \Delta\bm{\epsilon} = \ChoL_{\bt}\ChoL_{\bt}^T \Delta\bm{\epsilon} \\
    \ChoL_{\bt} &=  \mathsf{NN}_{\bt}(\bm{\epsilon},\mathcal{I}) 
\end{aligned}
\end{equation}
Here instead of outputting the stress or stress increment directly, the neural network outputs a lower triangular matrix $\ChoL_{\bt}$, which is the Cholesky factor of the  tangent stiffness matrix $\mathsf{H}_{\bt}$. 
The numerical approximation to \cref{EQ:CHOL} in the dynamic simulations has the following form 
\begin{equation}
\label{EQ:CHOL_NN_0}
    \bm{\sigma}^{n+1}  =  \ChoL_{\bt} \ChoL_{\bt}^T (\bm{\epsilon}^{n+1} - \bm{\epsilon}^{n})  + \bm{\sigma}^{n}
\end{equation}
here, the superscript $n$ indicates the time step.

An obvious advantage of \cref{EQ:CHOL} is the guarantee that the tangent stiffness matrix is a symmetric semidefinite matrix, which is true for commonly used constitutive relations. Additionally, when $\mathsf{NN}_{\bt}$ is bounded, we have
\begin{equation}
\label{EQ:TIME_CONSISTENCTY}
    \displaystyle \lim_{\Delta \bm{\epsilon} \rightarrow 0} \Delta\bm{\sigma} = 0
\end{equation}
which indicts the time-consistency for both path-dependent/independent constitutive relations. In contrast, both the $\bm{\sigma}$-NN with $\bm{\sigma}^{n+1} = \mathsf{NN}_{\bt}(\bm{\epsilon}^{n+1}, \bm{\epsilon}^{n}, \bm{\sigma}^{n})$ and the $\Delta\bm{\sigma}$-NN with $\bm{\sigma}^{n+1} = \mathsf{NN}_{\bt}(\bm{\epsilon}^{n+1}, \bm{\epsilon}^{n}, \bm{\sigma}^{n}) + \bm{\sigma}^{n}$ fail to satisfy the time-consistency condition~\cref{EQ:TIME_CONSISTENCTY} for path-dependent constitutive relations.
Finally, when $\mathsf{NN}_{\bt}$ does not depend on strain rate, \cref{EQ:CHOL} implies that the constitutive law is rate-independent, which is a feature of the hyperelastic or elasto-plastic materials.

In the following, we discuss how to adapt SPD-NN for different materials.

\begin{itemize}

\item \textbf{Linear Elasticity}

For linear elastic material, the tangent stiffness matrix is independent of the strain or stress. Therefore, no neural network is needed, and the constitutive relation reads
\begin{equation}
\label{EQ:LINEAR}
\bm{\sigma}^{n+1} = \mathsf{C}_{\theta}\bm{\epsilon}^{n+1} 
\end{equation}
where $\mathsf{C}_{\theta}$ is the parametric tangent stiffness tensor and the unknowns are simply entries in the tensor (no neural network). 

\item \textbf{Nonlinear Elasticity}

For nonlinear elastic material, the tangent stiffness matrix~\cref{EQ:CONSTITUTIVE_LAW_ELASTIC} depends only on the strain at the current time step. The constitutive relation~\cref{EQ:CHOL_NN_0} can be formulated as 
\begin{equation}
\label{EQ:CHOL_NN_1}
\bm{\sigma}^{n+1} = \ChoL_{\theta}(\bm{\epsilon}^{n+1} )\ChoL_{\theta}(\bm{\epsilon}^{n+1})^T(\bm{\epsilon}^{n+1} -  \bm{\epsilon}^{n})  + \bm{\sigma}^{n}
\end{equation}
Note the input $\bm\epsilon^{n+1}$ must be evaluated at the current time step $n+1$ because it is impossible to determine the stress at time step $n+1$ given only the information at time step $n$.

\item \textbf{Elasto-Plasticity}

For elasto-plastic material, the material behavior features transition from elastic behavior to plastic behavior. To model this effect, we consider two types of constitutive relations: 
\begin{enumerate}
	\item In the linear elasticity region, the constitutive relation is approximated by 
\begin{equation}
\label{EQ:PLASTICIT_LINEAR}
\bm\sigma_{\mathrm{elasticity}}^{n+1} =  \mathsf{C}_{\theta}(\bm\epsilon^{n+1} -  \bm\epsilon^{n})  + \bm\sigma^{n}
\end{equation}

\item In the plasticity region, the constitutive relation is approximated by 
\begin{equation}
\label{EQ:PLASTICIT_NONLINEAR}
\bm\sigma_{\mathrm{plasticity}}^{n+1} = \ChoL_{\theta}(\bm{\epsilon}^{n+1}, \bm{\epsilon}^{n},\bm{\sigma}^{n} )\ChoL_{\theta}(\bm\epsilon^{n+1}, \bm{\epsilon}^{n},\bm{\sigma}^{n})^T(\bm\epsilon^{n+1} -  \bm\epsilon^{n})  + \bm\sigma^{n}
\end{equation}
\end{enumerate}

However, since we do not know when the transition occurs (the yield strength $\sigma_Y$ is not available and strain hardening could strength the material), we can relax the constitutive relation using a differentiable\footnote{The differentiability of the transition function is necessary since we need to evaluate the gradient of $\bm{\sigma}^{n+1}$ with respect to $\bm{\sigma}^n$ during the training with indirect data.} transition function $D(\bm\sigma^{n}, \tilde{\sigma}_Y)$~(see \cref{fig:H}), whose value is between 0 and 1, as follows,
\begin{equation}
\label{EQ:CHOL_NN_2}
\bm\sigma^{n+1} = (1 - D(\bm\sigma^{n}, \tilde{\sigma}_Y)) \bm\sigma_{\mathrm{elasticity}}^{n+1} + D(\bm\sigma^{n}, \tilde{\sigma}_Y) \bm\sigma_{\mathrm{plasticity}}^{n+1}
\end{equation}
here $\tilde{\sigma}_Y$ is the estimated yield strength, which does not need to be accurate.  When the equivalent stress estimated from $\bm\sigma^{n}$ (such as von Mises stress)\footnote{In general, $D$ should depend on $\bm\sigma^{n+1}$ instead of $\bm\sigma^{n}$. Because it is difficult to express $\bm\sigma^{n+1}$ explicitly using this form, we assume $\bm\sigma^{n+1}\approx \bm\sigma^{n}$ and thus obtain \cref{EQ:CHOL_NN_2}.} is smaller than $\tilde{\sigma}_Y$, the material behavior is assumed to be linear and described by \cref{EQ:PLASTICIT_LINEAR}; otherwise, it is assumed to be described by the plastic form \cref{EQ:PLASTICIT_NONLINEAR}. It is worth noting
the plastic form~\cref{EQ:PLASTICIT_NONLINEAR} can degenerate to the linear elastic form~\cref{EQ:PLASTICIT_LINEAR}, but not vice versa. Therefore, the estimated yield strength $\tilde{\sigma}_Y$ should be smaller than the yield strength $\sigma_Y$. 
In the present study, historical data for plastic form~\cref{EQ:PLASTICIT_NONLINEAR} span only one time step, more historical data might be needed for materials with strong hysteresis.

\begin{figure}[hbtp]
\centering
	\scalebox{0.8}{\input{figures/H}}
  \caption{An exemplary transition function $D(\sigma, \tilde{\sigma}_Y) = \texttt{sigmoid}\left(100(\sigma^2-\tilde{\sigma}^2_Y)\right)$ and $\tilde\sigma_Y=1$. Here \texttt{sigmoid} is the sigmoid function $\texttt{sigmoid}(x) = (1+e^{-x})^{-1}$. }
  \label{fig:H}
\end{figure}

\end{itemize}

\begin{remark}
Generally, we can write the strain-stress relations for a linear elastic material in Voigt notation, as follows: 
\begin{equation*}
\begin{bmatrix}
\sigma_{11}\\
\sigma_{22}\\
\sigma_{33}\\
\sigma_{23}\\
\sigma_{13}\\
\sigma_{12}\\
\end{bmatrix}
=
\begin{bmatrix}
C_{1111} & C_{1122} & C_{1133}& C_{1123}& C_{1113} &C_{1112}\\
               & C_{2222} & C_{2233}& C_{2223}& C_{2213} &C_{2212}\\
               &                 & C_{3333}& C_{3323}& C_{3313} &C_{3312}\\
               &                 &                 & C_{2323}& C_{2313} &C_{2312}\\
              &     symm            &                  &                & C_{1313} &C_{1312}\\
              &                 &                  &                &                 &C_{1212}\\
\end{bmatrix}
\begin{bmatrix}
\epsilon_{11}\\
\epsilon_{22}\\
\epsilon_{33}\\
2\epsilon_{23}\\
2\epsilon_{13}\\
2\epsilon_{12}\\
\end{bmatrix}
\end{equation*}
The corresponding Cholesky factor $\mathsf{L}_{\bt}$ of  $\mathsf{C}_{\bt}$ is a full lower triangular matrix. For the materials considered in the present work, they are orthotropic, i.e., they have three mutually orthogonal planes of reflection symmetry. Thus the tangent stiffness matrix can be expressed by a block diagonal matrix (Orth-NN)
\begin{equation*}
\mathsf{C}_{\bt} = \begin{bmatrix}
C_{1111} & C_{1122} & C_{1133}&  &       &\\
               & C_{2222} & C_{2233}&   & &\\
               &                 & C_{3333}&  & &\\
               &                 &                 & C_{2323}&  &\\
              &     symm            &                  &                & C_{1313} &\\
              &                 &                  &                &                 &C_{1212}\\
\end{bmatrix}
\end{equation*}
Therefore, the associated Cholesky factor $\mathsf{L}_{\bt}$ has the following form
\begin{equation}
\label{EQ:CHOL_ORTH}
\mathsf{L}_{\bt} = \begin{bmatrix}
L_{1111}  &  & &  &       &\\
L_{2211}  & L_{2222} & &   & &\\
 L_{3311}  &  L_{3322}               & L_{3333} &  & &\\
               &                 &                 & L_{2323}&  &\\
              &               &                  &                & L_{1313} &\\
              &                 &                  &                &                 &L_{1212}\\
\end{bmatrix}
\end{equation}
Consequently, we can further simplify the neural network outputs to only the nonzero entries in \cref{EQ:CHOL_ORTH}.
\Cref{EQ:CHOL_ORTH} is the form of the Cholesky factor used in the present work.
\end{remark}

\begin{remark}

Another restriction on plastically stable materials in small strains is given  by the Drucker's postulate, which is related to the second law of thermodynamics.
The Drucker's postulate can be stated in different ways, one of which is that the second-order plastic work is non-negative~\cite[p.~328]{belytschko2013nonlinear}
\begin{equation}
\label{EQ:DRUCKER}
0 \leq \dot{\bm{\sigma}}^T  \dot{\bm{\epsilon}}^p
\end{equation}
Bringing \cref{EQ:CONSTITUTIVE_LAW_PLASTICITY} and the rate form of \cref{EQ:CHOL} into \cref{EQ:DRUCKER} leads to 
\begin{equation*}
0 \leq \dot{\bm{\sigma}}^T  \dot{\bm{\epsilon}}^p = \dot{\bm{\sigma}}^T  (\dot{\bm{\epsilon}} - \dot{\bm{\epsilon}}^e) 
= \dot{\bm{\sigma}}^T \mathsf{H}_{\bt}^{-1} \dot{\bm{\sigma}} -  \dot{\bm{\sigma}}^T \mathsf{C}_{\bt}^{-1} \dot{\bm{\sigma}} 
\end{equation*}
The Drucker's postulate requires 
\begin{equation*}
 \mathsf{H}_{\bt} \preceq \mathsf{C}_{\bt} 
\end{equation*}
which brings an additional constraint to the tangent stiffness matrix besides SPD.
We have attempted to design another neural network architecture imposing this constraint,
\begin{equation}
\label{EQ:BOUND_CHOL_NN}
    \Delta \bm{\sigma} =  \mathsf{H}_{\bt} \Delta \bm{\epsilon}  = [\mathsf{C}_{\bt}^{-1} + \ChoL_{\bt}\ChoL_{\bt}^T]^{-1} \Delta \bm{\epsilon} 
\end{equation}
However, this architecture does not improve the prediction and brings additional challenges in training the neural networks based on our observation. 
Therefore, this architecture~\cref{EQ:BOUND_CHOL_NN} is not used in the present work.

\end{remark}

\subsection{Training Methods}

Based on the types of available data, i.e., direct and indirect, we can employ different training methods. For the following discussion, we summarize the neural network based constitutive relations as follows
\begin{equation*}
	\bm{\sigma}^{n+1} = \mathsf{M}_{\bt}(\bm{\epsilon}^{n+1}, \bm{\epsilon}^{n}, \bm{\sigma}^{n}) := \begin{cases}
		 \mathsf{C}_{\bt}\bm{\epsilon}^{n+1}  & \mbox{Linear Elasticity}\\
		\ChoL_{\bt}(\bm{\epsilon}^{n+1} )\ChoL_{\bt}(\bm{\epsilon}^{n+1})^T(\bm{\epsilon}^{n+1} -  \bm{\epsilon}^{n})  + \bm{\sigma}^{n} & \mbox{Nonlinear Elasticity}\\
		 (1 - D(\bm\sigma^{n}, \tilde{\sigma}_Y)) \bm\sigma_{\mathrm{elasticity}}^{n+1} + D(\bm\sigma^{n}, \tilde{\sigma}_Y) \bm\sigma_{\mathrm{plasticity}}^{n+1} & \mbox{Elasto-Plasticity}
	\end{cases}
\end{equation*}

\subsubsection{Direct Data}
\label{SEC:DIRECT_DATA}

\textit{Direct} data consist of the input and the output of the NN-based constitutive relations such as strain-stress pairs or strain-stress increments pairs. These data points come from experimental measurements and numerical simulation results. The comprehensive strain-stress data measurement relying on simple mechanical tests, such as tensile or bending tests, might be challenging. However, comprehensive strain-stress data generated from sub-scale simulations, such as representative volume element (RVE) simulations~\cite{hashin1983analysis, sun1996prediction, feyel2000fe2, kanit2003determination, yuan2008toward} or post-processed from direct numerical simulations, are widely used to train neural networks~\cite{le2015computational, bessa2017framework, ling2016machine, wang2018multiscale}.

Mathematically, the direct data are given in terms of $N$ sequences of strain-stress pairs at $n$ time snapshots. 
\begin{equation*}
    (\bm{\epsilon}_j^1, \bm{\sigma}_j^1),\,(\bm{\epsilon}_j^2, \bm{\sigma}_j^2),\,\cdots,\,(\bm{\epsilon}_j^n, \bm{\sigma}_j^n)\quad j = 1,2,3,\ldots,N
\end{equation*}
Here the superscripts indicate time and the subscripts indicate the sequential number.
We train the neural network by solving a minimization problem
\begin{equation}
\label{EQ:LOSS_DIRECT}
  \arg\min_{\bt}  \Loss(\bt) := \sum_{j=1}^{N} \sum_{i=2}^{n}\left(\bm{\sigma}_{j}^{i} - \mathsf{M}_{\bt}(\bm{\epsilon}_{j}^{i}, \bm{\epsilon}_{j}^{i-1}, \bm{\sigma}_{j}^{i-1}) \right)^2
\end{equation}

\subsubsection{Indirect Data}
\label{SEC:INDIRECT_DATA}
\textit{Indirect} data consist of deformation data from structure coupons under different load conditions.
    Deformations are measured by techniques such as digital image correlation or grid method~\cite{surrel1994moire}.
    These techniques can record complete heterogeneous fields, which are rich in the constitutive relations.
    The virtual fields method~\cite{grediac2006virtual, avril2008overview, geymonat2002identification, feng2001genetic} has been designed to apply the finite element method (FEM) to bridge the full-field data with parametric constitutive relations. Recently, the method is generalized to an end-to-end training procedure to learn neural-network~(or its counterparts) based constitutive relations from the full-field data~\cite{huang2019predictive}.
  
 The indirect data are given by $N$ full-field deformation-load sequential data at  $n$ time snapshots
 \begin{equation*}
    (\bu_j^1, \bF_j^1),\,(\bu_j^2, \bF_j^2),\,\cdots,\,(\bu_j^n, \bF_j^n)\quad j = 1,2,3,\ldots,N
\end{equation*}
Here the superscripts indicate time and the subscripts indicate the sequential number.

\begin{figure}[htpb]
\centering
  \includegraphics[width=0.9\textwidth]{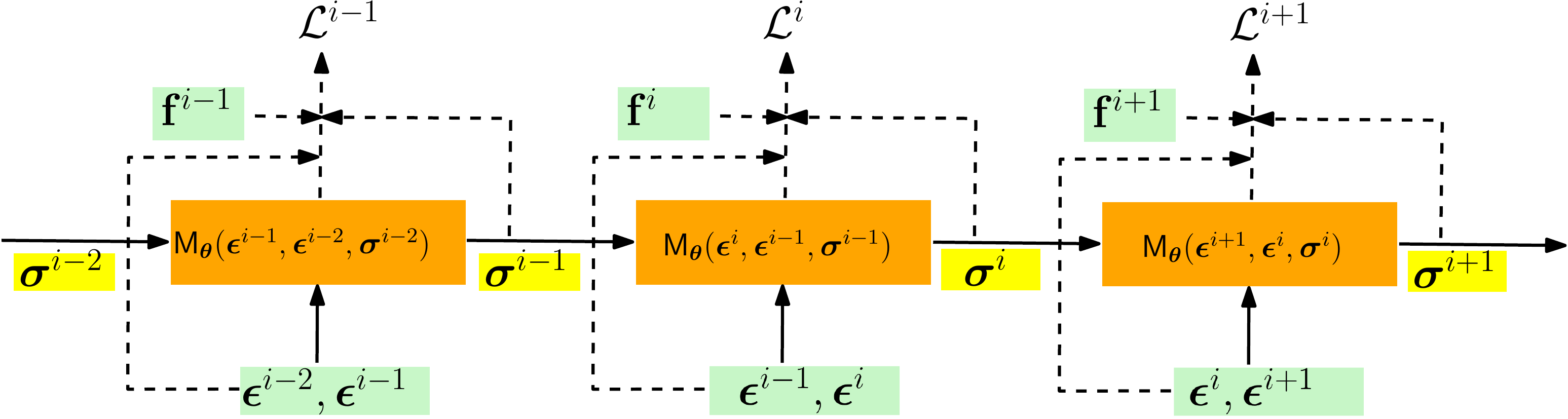}
  \caption{The procedure of evaluating the residual $\Loss^{i}$ at each time step for training with indirect data  approach in \cref{EQ:LOSS_INDIRECT}.}
  \label{fig:rnn}
\end{figure}

We can compute the acceleration and the stress using the formulas
\begin{align}
    \ddot{\bu}_j^i &:= \frac{\bu_j^{i+1} - 2\bu_j^i + \bu_j^{i-1}}{\Delta t^2} \label{EQ:ACCELERATION}\\
    \bm{\sigma}_j^{i}(\bt) &:= \mathsf{M}_{\bt}(\bm{\epsilon}(\bu_j^i),\bm{\epsilon}(\bu_j^{i-1}), \bm{\sigma}_j^{i-1}(\bt)), \quad i = 2,\ldots, n-1\\
    \ddot{\bu}_j^1 &= 0, \; {\bu}_j^1 = 0,\; \bm{\sigma}_j^{1}(\theta) = 0
\end{align}
Here $\Delta t$ is the time step, which is assumed to be constant.
We train the neural network by solving a minimization problem
\begin{equation}
\label{EQ:LOSS_INDIRECT}
\arg\min_{\bt}    \Loss(\bt) := \sum_{j=1}^{N} \sum_{i=2}^{n-1}\left(\bM\ddot{\bu}_j^i + \bP(\bu_j^i,     \bm{\sigma}_j^{i}(\theta)) - \bF_j^i \right)^2
\end{equation}
Here $i$ is the index for time and $j$ is the index for spacial degrees of freedom. 
It is worth mentioning, in \cref{EQ:LOSS_INDIRECT}, the predicted stress $\bm{\sigma}_j^{i}(\bt)= \mathsf{M}_{\bt}(\bm{\epsilon}(\bu_j^i), \bm{\epsilon}(\bu_j^{i-1}),\bm{\sigma}_j^{i-1}(\bt))$ depends on the {\it predicted} stress $\bm{\sigma}_j^{i-1}(\bt)$ from the last time step. The procedure of evaluating the residual $\Loss^{i}=\left(\bM\ddot{\bu}_j^i + \bP(\bu_j^i,     \bm{\sigma}_j^{i}(\theta)) - \bF_j^i \right)^2$ at each time step in \cref{EQ:LOSS_INDIRECT} is depicted in \cref{fig:rnn}. The procedure resembles the recurrent neural network in deep learning, where the state variables are the stresses. Like the recurrent neural network, the stress predictions at different time steps must be sequentially computed and so does the back-propagated gradients, and thus the computation is hard to be parallelized. Additionally, the training of recurrent neural network suffers from exploding and vanishing  gradients problem~\cite{pascanu2013difficulty}, which poses a challenge for indirect data training of SPD-NNs as well.

\subsection{Initialization for Indirect Data Training}
Another notable challenge in training the recurrent neural network based constitutive relation is the existence of many local minima, since the training data set is small in the present work, and different constitutive relations under the same external loads may produce similar displacements. 
A common strategy to find a good local minimum is to start from multiple random weights and biases. However, this strategy is expensive and does not guarantee a good initial guess. We thereby propose an initialization technique that yields a set of reasonably good initial weights and biases.

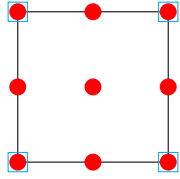
\begin{figure}[htbp]
\centering
\begin{tikzpicture}
  \draw (0,0) -- (2,0) -- (2,2) -- (0,2)-- (0,0);
\foreach \i in {0,...,2}{
    \foreach \j in {0,...,2}
    {   
        \filldraw[red] (\i,\j) circle (3pt);
    }
}
\foreach \i in {0,2}{
    \foreach \j in {0,2}
    {
        \node[draw, rectangle, cyan] at (\i,\j) {};
    }
}
\end{tikzpicture}
    \caption{Schematic of the 2D quadratic element, with quadratic nodes~(red circle) for displacements and linear nodes~(blue empty square) for stress components.}
\label{FIG: ELEMENT}
\end{figure}

The key idea is that according to \cref{EQ:LINEAR_MOMENTUM_SEMI}, we have the relation between the internal force and the stress field
\begin{equation}
\label{EQ:PREFIT}
     \bP(\bu_j^i, \bm{\sigma}_j^{i}) = \bF_j^i  - \bM\ddot{\bu}_j^i
\end{equation}
However, solving the stress field $\bm{\sigma}_j^{i}$ from \cref{EQ:PREFIT} is generally an underdetermined problem, since the number of equations is fewer than the number of stress unknowns. To alleviate this, quadratic elements are used to represent the displacement field $\bu$, and in each quadratic element, a linear stress field is assumed and approximated~(see \cref{FIG: ELEMENT}).

The number of stress unknowns is roughly $3(n_x + 1)(n_y + 1)$ for a 2D plate, and $6(n_x + 1)(n_y + 1)(n_z + 1)$ for a 3D cube, the number of equations is roughly  $2(2n_x + 1)(2n_y + 1)$ for a 2D plate, and $3(2n_x + 1)(2n_y + 1)(2n_y + 1)$ for a 3D cube. Here $n_x$, $n_y$, and $n_z$ represent the number of quadratic elements in each direction. Therefore, the number of equations outnumbers the number of unknowns, and \cref{EQ:PREFIT} becomes an overconstrained problem, the least-square fitting is applied for solving it.
And then the stress field is approximated linearly at each Gaussian point in each quadratic element.

Once we solve for $\bm{\sigma}_j^i$ from \Cref{EQ:PREFIT}, we can use the technique in \Cref{SEC:DIRECT_DATA} to pre-train the neural network. Although the least square approximation of the stress field is poor~(the error can be larger than 100\%), but the approximation is qualitatively correct and sufficient for obtaining a good initial guess.

\section{Applications}
In this section, we present numerical results from solid mechanics for the proposed NN based constitutive relations: 
\begin{itemize}
    \item Problem with a 1D truss coupon made of elasto-plastic materials under dynamic loading, which compares the proposed SPD-NN~(\ref{EQ:CHOL_NN_0}) and other neural network architectures, including $\bm{\sigma}$-NN~(\ref{EQ:NN1}) and $\Delta\bm{\sigma}$-NN~(\ref{EQ:NN2}).
    \item Problems with 2D thin plates made of hyperelastic materials, elasto-plastic materials, and multiscale fiber-reinforced materials under dynamic loading, which demonstrate the effectiveness of the proposed SPD-NN~(\ref{EQ:CHOL_NN_0}) for learning path-independent, path dependent~(hysteresis), and  multiscale constitutive relations. 
\end{itemize}
In all problems, the training data and test data of the strain/stress fields and the displacement fields  are generated numerically.

\subsection{1D Trusses with Elasto-plasticity}
\label{SEC:TRUSS}
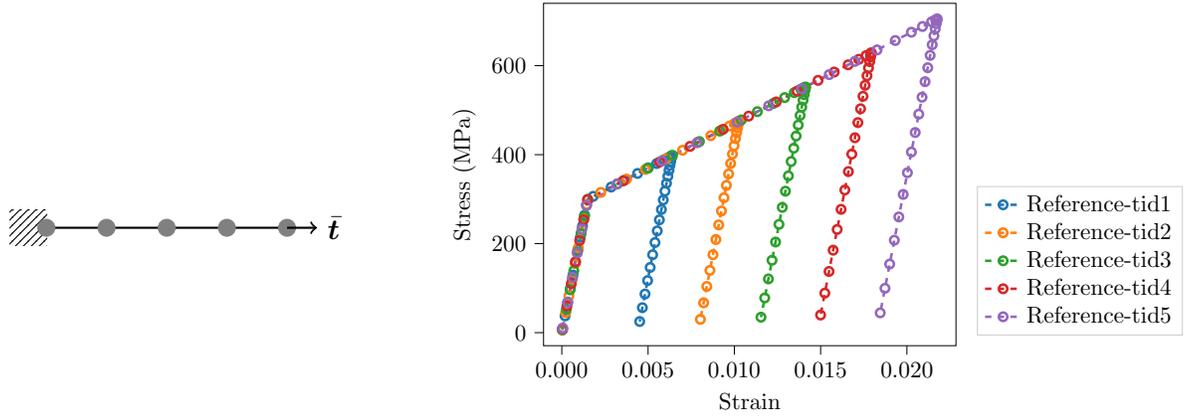
\begin{figure}[htbp]
\centering
\begin{subfigure}[b]{0.35\textwidth}
\begin{tikzpicture}[scale=0.8]
\filldraw [white] (0,-2) circle (4pt);
\draw[thick] (0,1) -- (4,1);
\fill [pattern = north east lines] (-0.6,0.7) rectangle (0.0, 1.3);
\filldraw [gray] (0,1) circle (4pt);
\filldraw [gray] (1,1) circle (4pt);
\filldraw [gray] (2,1) circle (4pt);
\filldraw [gray] (3,1) circle (4pt);
\filldraw [gray] (4,1) circle (4pt);
\draw[thick,->] (4, 1.0) -- (4.5, 1.0) node[anchor=west]{$\bar{\bm{t}}$};
\end{tikzpicture}
\caption*{}
\end{subfigure}%
\begin{subfigure}[b]{0.45\textwidth}
\scalebox{0.8}{
\begin{tikzpicture}

\definecolor{color0}{rgb}{0.12156862745098,0.466666666666667,0.705882352941177}
\definecolor{color1}{rgb}{1,0.498039215686275,0.0549019607843137}
\definecolor{color2}{rgb}{0.172549019607843,0.627450980392157,0.172549019607843}
\definecolor{color3}{rgb}{0.83921568627451,0.152941176470588,0.156862745098039}
\definecolor{color4}{rgb}{0.580392156862745,0.403921568627451,0.741176470588235}

\begin{axis}[
scaled ticks=false,
legend cell align={left},
legend style={fill opacity=0.8, draw opacity=1, text opacity=1, at={(1.55,0.03)}, anchor=south east, draw=white!80!black},
tick align=outside,
tick pos=left,
x grid style={white!69.0196078431373!black},
xlabel={Strain},
xmin=-0.00105827502561468, xmax=0.0228439737988373,
xtick style={color=black},
xtick={-0.005,0,0.005,0.01,0.015,0.02,0.025},
xticklabels={-0.005,0.000,0.005,0.010,0.015,0.020,0.025},
y grid style={white!69.0196078431373!black},
ylabel={Stress (MPa)},
ymin=-29.3374335935994, ymax=740.125757651418,
ytick style={color=black}
]
\addplot [very thick, color0, dashed, mark=*, mark size=2, mark options={solid,fill opacity=0}]
table {%
2.81908300422345e-05 5.6381660084469
0.000187683746315012 37.536749263003
0.000343668309733048 68.7336618761499
0.000497131170133144 99.4262339455293
0.000647426534739178 129.485306866374
0.000793680569222919 158.736113763161
0.000935002884876083 187.000576893856
0.00107052766967965 214.10553385463
0.0011994244900524 239.884897929236
0.00132090405991174 264.180811901153
0.0014342230610353 286.84461212591
0.00181948686220043 306.389737235902
0.00286058017109304 327.211603413754
0.00363133812017686 342.626762395431
0.00439303580216942 357.8607160268
0.00501126813769093 370.225362664648
0.00550126289692884 380.025257819741
0.005925572033542 388.511440327341
0.00616450523375492 393.29010432007
0.00635666286742942 397.133256954377
0.00637562174283685 397.512434474214
0.00636671706840759 395.731499733427
0.00634593325314733 391.574736668351
0.00631314528491512 385.017143021575
0.00626851298481862 376.090683002266
0.0062123037778227 364.848841603093
0.00614486008014034 351.36010206663
0.00606659334340647 335.706754719875
0.00597798123913278 317.98433386516
0.00587956482493916 298.301051026464
0.00577194534666045 276.777155370758
0.00565578066963469 253.544219965645
0.00553178136139091 228.744358316937
0.00540070644851217 202.529375741241
0.00526335887106663 175.05986025219
0.00512058065950207 146.504217939342
0.00497324786065024 117.037658169041
0.00482226524123333 86.841134285731
0.00466856079892522 56.1002458241831
0.0045130801125639 25.0041085519967
};
\addlegendentry{Reference-tid1}
\addplot [very thick, color1, dashed, mark=*, mark size=2, mark options={solid,fill opacity=0}]
table {%
3.3828824373981e-05 6.76576487479619
0.00022521207950716 45.0424159014331
0.000412373697774497 82.4747395109665
0.000596498231574946 119.299646255749
0.000776811507525514 155.362301445341
0.00095226595258949 190.453190458192
0.0011217943589993 224.358871740243
0.00128435920445176 256.871840830822
0.00143896557647292 287.793115235133
0.00226031264904289 315.206252974943
0.00375279743501823 345.05594869445
0.00484677614746998 366.935522943486
0.00606140677110469 391.228135409481
0.00699684643378408 409.936928473761
0.00789404011247931 427.88080201168
0.00863459415370959 442.691882790921
0.00922006330697028 454.401265812116
0.00972028436584022 464.405686766344
0.0100104725168215 470.209449777253
0.0102340035915807 474.680071150117
0.0102592188325928 475.184375875962
0.0102486223960965 473.065088739445
0.0102237979372495 468.100196955083
0.0101846176759812 460.264144701051
0.010131281114053 449.596832315412
0.0100641084586834 436.162301241509
0.00998350822929268 420.042255363374
0.009889971176608 401.334844826462
0.00978406708055914 380.154025616717
0.00966644140863844 356.628891232611
0.00953781152907055 330.902915319083
0.00939896248220849 303.133105946728
0.00925074233641499 273.489076788093
0.00909405715450375 242.152040405919
0.00892986559765286 209.313729035828
0.00875917319557656 175.175248620653
0.00858302631388977 139.945872283393
0.0084025058517499 103.841779855523
0.00821872070492297 67.0847504902448
0.00803280103136731 29.9008157792213
};
\addlegendentry{Reference-tid2}
\addplot [very thick, color2, dashed, mark=*, mark size=2, mark options={solid,fill opacity=0}]
table {%
3.94667614855144e-05 7.89335229710289
0.000262737609167358 52.5475218334735
0.000481069672157902 96.2139343718435
0.000695845603162529 139.169120551953
0.000906163111303343 181.232622179406
0.00111080123165105 222.160246249025
0.00130851635936997 261.703271792929
0.00149809974672009 299.619949263072
0.00323193202914041 334.638640574785
0.00498143865949836 369.628773181944
0.00642003916262919 398.400783244561
0.00797281755091166 429.456351001532
0.00914725186338931 452.945037251085
0.0104058499942009 478.116999823974
0.0113336195708573 496.672391260666
0.0122478758536451 514.957516869412
0.0129161596962976 528.323193692812
0.0134717550260704 539.435100075987
0.0138536915245349 547.073829857545
0.014056055634465 551.121111969057
0.0141307531614943 552.16874558956
0.0141190622100165 549.830555294001
0.014090225766778 544.063266593623
0.0140446899430925 534.956101855122
0.0139827141887345 522.560950983515
0.0139046646260146 506.951038439551
0.0138110123627919 488.220585795028
0.0137023257792617 466.48326908901
0.0135792655343288 441.871220102479
0.0134425802921606 414.534171668899
0.0132931022556121 384.638564359259
0.013131742271295 352.366567495903
0.0129594844712875 317.915007494492
0.0127773804659743 281.494206431959
0.0125865431151509 243.326736267379
0.0123881399091834 203.646095073998
0.0121833859950344 162.695312244315
0.0119735368846105 120.725490159676
0.0117598808853918 77.9942903160686
0.0115437312956711 34.7643723720768
};
\addlegendentry{Reference-tid3}
\addplot [very thick, color3, dashed, mark=*, mark size=2, mark options={solid,fill opacity=0}]
table {%
4.51046413791294e-05 9.02092827582588
0.000300260336077017 60.0520672154069
0.000549756238057296 109.951247585084
0.000795173299198921 159.034659786247
0.00103548137786558 207.096275518654
0.00126928646524576 253.85729299479
0.00149516898232001 299.033796409798
0.00355133995062165 341.026799007126
0.00567484597956317 383.496919585956
0.00743344711630352 418.668942320763
0.00934967135292372 456.993427001234
0.0108288793207543 486.577586357845
0.0124225503924084 518.451007700236
0.0136387541654704 542.775083044153
0.0148513119390466 567.026238444658
0.0157851224274399 585.702448167147
0.0165985074851352 601.970149281223
0.0172150121441176 614.300242215527
0.0176320224085769 622.640447486094
0.0178981142441878 627.962283958963
0.0179436199161506 628.872397857755
0.0179295463403247 626.057682689803
0.0178967226008992 619.492934787697
0.017844905676382 609.129549883874
0.0177743649731829 595.021409244039
0.0176855224412491 577.252902857315
0.0175789160752314 555.931629653792
0.0174551924122147 531.186897050494
0.0173151023185645 503.168878320501
0.0171594966172123 472.047738050135
0.0169893211571493 438.012646037617
0.0168056113112884 401.270676865524
0.0166094859289794 362.045600403845
0.0164021407736287 320.576569333835
0.0161848414782203 277.116710252292
0.0159589160543423 231.931625476836
0.0157257469933068 185.297813269899
0.0154867630009424 137.501014797199
0.0152434304105612 88.8344967211393
0.0149972443214169 39.5972788924832
};
\addlegendentry{Reference-tid4}
\addplot [very thick, color4, dashed, mark=*, mark size=2, mark options={solid,fill opacity=0}]
table {%
5.074246405712e-05 10.148492811424
0.000337780261017185 67.5560522034423
0.000618433400076465 123.686679981952
0.00089448133459806 178.896266851926
0.00116476633982754 232.953267896654
0.00142772171299512 285.544342530296
0.00318226176568308 333.645235306845
0.00581299437240194 386.259887441222
0.00784073221631451 426.814644319474
0.0101609446221 473.218892435184
0.011991391746272 509.827834905256
0.0139021381713766 548.042763407348
0.0154971291470777 579.942582921305
0.0169777635961664 609.555271903079
0.0182722871233024 635.44574244581
0.0193288946916413 656.577893643587
0.0202607545712155 675.215091055453
0.0209032935644708 688.065870919835
0.0214282588513031 698.565176486331
0.0216625937272154 703.251873808312
0.0217575079431804 705.150158049372
0.0217407708518332 701.802739731521
0.021703717768082 694.392122934276
0.0216456127637646 682.771122064506
0.0215665632833797 666.961225987346
0.021467009804082 647.050530127811
0.0213475492860024 623.158426511923
0.0212089042082939 595.429410970258
0.0210519144057581 564.031450463177
0.0208775317070457 529.154910720769
0.0206868143899765 491.011447307036
0.0204809210586122 449.832781034294
0.0202611039176498 405.869352841962
0.0200287014704343 359.388863399008
0.0197851306753955 310.674704391431
0.019531878599492 260.024289210925
0.0192704936106606 207.747291444842
0.0190025761546831 154.163800249565
0.0187297691652512 99.6024023634144
0.0184537481592861 44.3982011706363
};
\addlegendentry{Reference-tid5}
\end{axis}

\end{tikzpicture}}
\caption*{}
\end{subfigure}%

\caption{1D truss problem setup~(left): the truss coupon consists of 4 elements with left end fixed and force load on the right end; Extracted strain-stress curves for different loading conditions (right).}
\label{FIG:TRUSS_BC}
\end{figure}

In this section, we consider 1D truss elements with a length $L_x = 1\textrm{m}$, a cross section area $A = 0.005\textrm{m}^2$ and a density $8000 \textrm{kg/m}^3$.  The truss coupon, clamped on the left end, is tested under 5 loading conditions~(see \cref{FIG:TRUSS_BC}-left). The setting is corresponds to a uniaxial tensile test. The prescribed time-dependent load force $\bar{\bm{t}}$ consists of both loading and unloading parts and takes the form
$$\bar{\bm{t}} = p \sin\left(\frac{t\pi}{T}\right), \, p = (0.4\textrm{tid} + 1.6)\times 10^{6} \; \textrm{N}$$
here $\textrm{tid} = 1,2,3,4$ and $5$ are the test indices.  The total simulation time is $T = 0.2s$. The case $\textrm{tid} = 3$ is used as test set and all the other tests are used as training sets.

The trusses are  made of elasto-plastic materials. The young's modulus is 
$$E = 200~\textrm{GPa}$$
The yield function with isotropic hardening  has the form
$$f = |\sigma| - \sigma_Y - K\alpha$$
The yield strength is $\sigma_Y = 0.3~\textrm{GPa}$, the plastic modulus is $K = \frac{200}{9}~\textrm{GPa}$, the internal hardening variable follows the simplest evolutionary equation
$$\dot{\alpha} = \dot{\lambda}$$

The truss coupon consists of 4 truss elements, which are modeled as geometric nonlinear truss elements~\cite[p.~63]{de2012nonlinear}. The time step size is $\Delta t = 0.001~\mathrm{s}$.

As for the SPD-NN~\cref{EQ:CHOL_NN_2}, the estimated Young's modulus is assumed known, which can be separately calibrated using small external loads (e.g., using the method proposed in \cite{huang2019predictive}), and the estimated yield strength is $\tilde{\sigma}_Y = 0.1~\textrm{GPa}$, the transition function is 
\begin{equation}
\label{EQ:TRANSIT_TRUSS}
D(\sigma^{n}, \tilde{\sigma}_Y) = \texttt{sigmoid}\left(\frac{(\sigma^{n})^2 - \tilde{\sigma}_Y^2} {d\tilde{\sigma}_Y^2}\right)
\end{equation}
here $d=0.1$ is the nondimensional parameter.

Because the strain and stress pairs can be extracted from the uniaxial tensile tests, we can apply the direct input-output data training method \cref{SEC:DIRECT_DATA} to train the neural network. 
\cref{FIG:TRUSS_BC}-right shows the strain-stress curves exhibiting plastic loading and elastic unloading phenomena, which the neural networks are trained to learn.
These curves start from the origin and rise with the slope $E$ in the linear elastic region, until the stress reaches the yield stress.  Then the curves enter the strain hardening region, where plastic deformations happen, until the elastic unloading.
The slope of the unloading curve is typically equal to the slope in the elastic (initial) region of the stress-strain curve. 
The plasticity deformations result in permanent strains, which cannot be recovered by the elastic unloading.  
Since training neural networks involve highly non-convex optimization problems, we start from 10 different initial weights for all NN training.

To evaluate the quality of approximation, we propose two kinds of tests
\begin{enumerate}
\item NN test: extract the sequential strain-stress data $(\epsilon^{i} , \sigma^{i})$
at each Gaussian quadrature point from the test data, and compare the predicted stress and the reference stress for each tuple~$(\epsilon^{i+1} , \epsilon^{i},  \sigma^{i}; \sigma^{i+1})$.
\item NN-FEM test: embed the learned constitutive relation $\mathsf{M}_{\bt}$ into the finite element framework, solve the governing equation of the truss coupon under the corresponding loading condition, and finally compare the predicted strain-stress paths and the truss deformations. 
\end{enumerate} 
Because the optimization results depend on the initial guess for the neural network weights and biases, we pick the neural network with the minimal loss on the training set for the NN test and NN-FEM test. Other choices for selecting the candidate neural networks such as cross-validation \cite{friedman2001elements} can also be used. 
It is worth noting that NN-FEM test is more challenging than the NN test, but is more relevant for predictive modeling. Indeed, the NN test does not take the numerical stability in the predictive modeling into account and only evaluates the NN's ability to fit the strain-stress curve. As we will see in the following examples, although SPD-NNs and other NN methods fit the strain-stress curves equally well in the NN test, SPD-NNs are significantly preferable due to their predicting power in new scenarios.

\subsubsection{Comparison of SPD-NN, \texorpdfstring{$\bm{\sigma}$-NN}{NN1}, and \texorpdfstring{$\Delta\bm{\sigma}$-NN}{NN2}}\label{sect:compare_cholnn}
\begin{figure}[htpb]
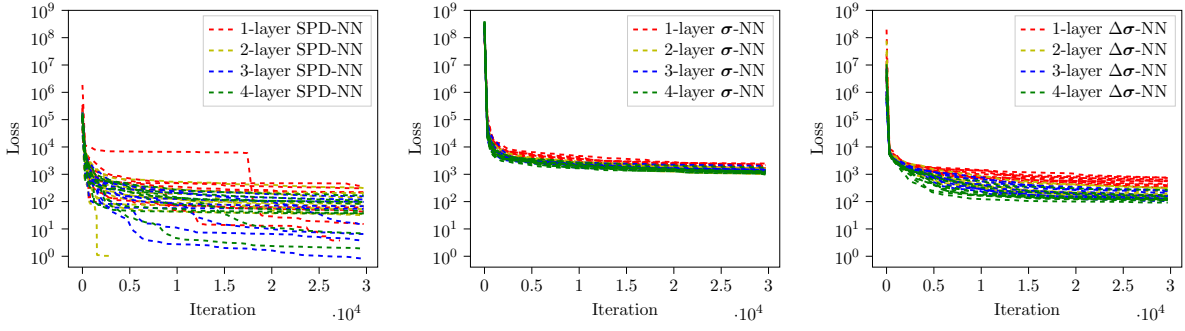

\scalebox{0.6}{\input{figures/nnlearn_piecewise_loss.tex}}~
\scalebox{0.6}{\input{figures/nnlearn_ae_scaled_loss.tex}}~
\scalebox{0.6}{\input{figures/nnlearn_ae_scaled2_loss.tex}}~
  \caption{The losses evaluated on the training set at each training iteration with SPD-NN~(left), $\bm{\sigma}$-NN~(middle), and $\Delta\bm{\sigma}$-NN~(right) with different number of layers for the 1D truss problem. Different curves correspond to different initial guesses. }
  \label{FIG:TRUSS-LOSS}
\end{figure}

The performance of the aforementioned neural network architectures, including the proposed SPD-NN, $\bm{\sigma}$-NN and $\Delta\bm{\sigma}$-NN,  with different hyper-parameters are compared. The neural networks considered contain 1, 2, 3 and 4 hidden layers with 20 neurons in each layer, and \texttt{tanh} as the activation function is used. Both input and output are 1-dimensional. 
The losses  at each training step  are reported in \cref{FIG:TRUSS-LOSS}. Different initial weights lead to different local minima, as with any neural network based data-driven approaches.
And deeper neural networks perform better in terms of the training error.
SPD-NN achieves significantly smaller training errors, compared with $\bm{\sigma}$-NN and $\Delta\bm{\sigma}$-NN.

\begin{figure}[htpb]
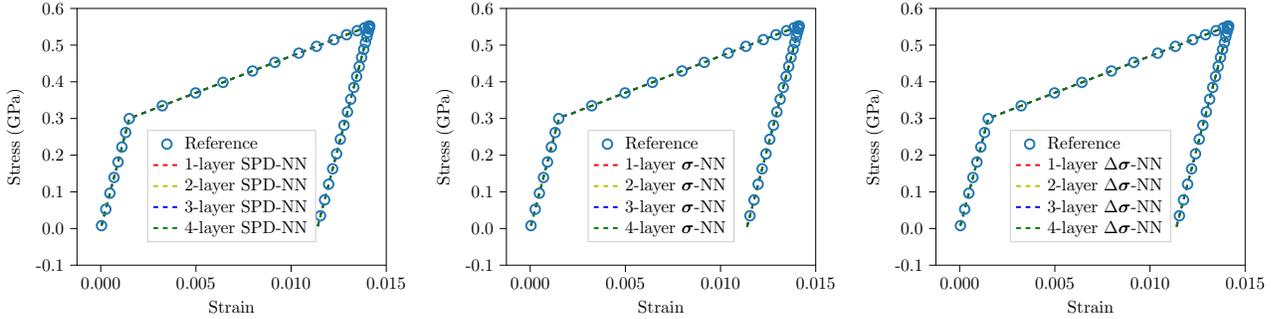

\scalebox{0.6}{\input{figures/nnlearn_piecewise_p2p_strain_stress_tid3.tex}}~
\scalebox{0.6}{\input{figures/nnlearn_ae_scaled_p2p_strain_stress_tid3.tex}}~
\scalebox{0.6}{\input{figures/nnlearn_ae_scaled2_p2p_strain_stress_tid3.tex}}~
  \caption{The strain stress curve results for SPD-NN~(left), $\bm{\sigma}$-NN~(middle), and $\Delta\bm{\sigma}$-NN~(right) in the NN test. }
  \label{FIG:TRUSS-NN-TEST}
\end{figure}

\begin{figure}[htpb]
\scalebox{0.6}{
\begin{tikzpicture}

\definecolor{color0}{rgb}{0.12156862745098,0.466666666666667,0.705882352941177}
\definecolor{color1}{rgb}{0.75,0.75,0}

\begin{axis}[
scaled ticks=false,
legend cell align={left},
legend style={fill opacity=0.8, draw opacity=1, text opacity=1, at={(0.5,0.09)}, anchor=south, draw=white!80!black},
tick align=outside,
tick pos=left,
x grid style={white!69.0196078431373!black},
xlabel={Strain},
xmin=-0.00125, xmax=0.015,
xtick style={color=black},
xtick={-0.005,0,0.005,0.01,0.015},
xticklabels={-0.005,0.000,0.005,0.010,0.015},
y grid style={white!69.0196078431373!black},
ylabel={Stress (GPa)},
ymin=-0.1, ymax=0.6,
ytick style={color=black},
ytick={-0.1,0,0.1,0.2,0.3,0.4,0.5,0.6},
yticklabels={-0.1,0.0,0.1,0.2,0.3,0.4,0.5,0.6}
]
\addplot [very thick, color0, mark=*, mark size=3, mark options={solid,fill opacity=0}, only marks]
table {%
3.94667614855144e-05 0.00789335229710289
0.000262737609167358 0.0525475218334735
0.000481069672157902 0.0962139343718435
0.000695845603162529 0.139169120551953
0.000906163111303343 0.181232622179406
0.00111080123165105 0.222160246249025
0.00130851635936997 0.261703271792929
0.00149809974672009 0.299619949263072
0.00323193202914041 0.334638640574785
0.00498143865949836 0.369628773181944
0.00642003916262919 0.398400783244561
0.00797281755091166 0.429456351001532
0.00914725186338931 0.452945037251085
0.0104058499942009 0.478116999823974
0.0113336195708573 0.496672391260666
0.0122478758536451 0.514957516869412
0.0129161596962976 0.528323193692812
0.0134717550260704 0.539435100075987
0.0138536915245349 0.547073829857545
0.014056055634465 0.551121111969057
0.0141307531614943 0.55216874558956
0.0141190622100165 0.549830555294001
0.014090225766778 0.544063266593623
0.0140446899430925 0.534956101855122
0.0139827141887345 0.522560950983515
0.0139046646260146 0.506951038439551
0.0138110123627919 0.488220585795028
0.0137023257792617 0.46648326908901
0.0135792655343288 0.441871220102479
0.0134425802921606 0.414534171668899
0.0132931022556121 0.384638564359259
0.013131742271295 0.352366567495903
0.0129594844712875 0.317915007494492
0.0127773804659743 0.281494206431959
0.0125865431151509 0.243326736267379
0.0123881399091834 0.203646095073998
0.0121833859950344 0.162695312244315
0.0119735368846105 0.120725490159676
0.0117598808853918 0.0779942903160686
0.0115437312956711 0.0347643723720768
};
\addlegendentry{Reference}
\addplot [very thick, red, dashed]
table {%
3.94679440676673e-05 0.00789334923639001
0.000262738826612118 0.0525475043331484
0.00048107271369646 0.0962139252371662
0.00069583230108652 0.139169877729095
0.000906118261649265 0.181230855412958
0.00111072725748992 0.222164398367546
0.00130846692939043 0.261704027903225
0.00149780430975328 0.299625386743327
0.00322073698058162 0.334634784283324
0.00497141759711679 0.369646603979641
0.00640905266685183 0.398397858717167
0.00796242334858364 0.429463702299746
0.0091368495401388 0.452950664115717
0.0103953086359504 0.478118581921088
0.0113235419285654 0.496682153894242
0.0122373342290999 0.514957083471732
0.0129062105418037 0.528333827514861
0.0134613353805085 0.539435472853571
0.01384408499422 0.547090103741843
0.0140446810128468 0.551103598277519
0.0141144646147547 0.552226215718619
0.0141016736526106 0.549852714889
0.0140728062971872 0.544079202739179
0.0140272203102717 0.534968027417544
0.0139652048758403 0.522566330593017
0.0138871529315994 0.506950979451367
0.0137936049163881 0.488226210538445
0.0136849584051303 0.466489061221919
0.0135619423465475 0.441876458763426
0.0134252413569952 0.414539655897098
0.0132757235344991 0.384644142695147
0.0131143517214963 0.352372159365777
0.0129421029574677 0.317917527405103
0.0127600133603791 0.281494403571133
0.0125692048976074 0.24332803962944
0.0123708353312709 0.20364748259259
0.012166113985704 0.162695975564466
0.0119562478313953 0.120721653211746
0.0117426344590417 0.0779829208922336
0.0115265955979171 0.0347627005831748
};
\addlegendentry{1-layer SPD-NN}
\addplot [very thick, color1, dashed]
table {%
3.94693857812017e-05 0.00789387684187272
0.00026273350435116 0.0525475208436555
0.000481065572526977 0.0962139345872402
0.000695826593458427 0.139169883405026
0.000906120895794703 0.181230788299116
0.00111073796867424 0.222164345053638
0.00130848851002855 0.261704207958475
0.0014980717609794 0.299612633349505
0.00322104008315237 0.334635996196469
0.00497151477192064 0.36964530481725
0.00640922212129162 0.398399262154408
0.0079623827938683 0.429462331873008
0.00913684943308818 0.452951547838999
0.0103951707076609 0.478117864717981
0.0113234156057873 0.496682677150932
0.0122371472680012 0.514957233000857
0.0129059743610152 0.528333715647913
0.0134610962252047 0.539436029802089
0.0138434471356177 0.547083039120348
0.0140455498451027 0.551124823005619
0.0141196206474413 0.552171256159152
0.014107938040025 0.549835590284423
0.0140790976116537 0.544067499660047
0.0140335770320306 0.534961799971374
0.013971652781546 0.522567891562455
0.0138935817295253 0.506955119409613
0.013799943951841 0.488222472030264
0.0136912614570749 0.466486856555672
0.0135682022022105 0.441875339599882
0.0134315345861744 0.414541176848854
0.0132820302685791 0.384644059809259
0.013120709754062 0.352372482961902
0.0129485158491792 0.317920905392218
0.0127664961780931 0.281497870170169
0.01257574852076 0.243330750322947
0.012377362951292 0.203649578802226
0.0121726930412075 0.162698521506463
0.0119628909412503 0.12072808160533
0.0117493270544831 0.077996301673673
0.0115332409868918 0.0347659127424834
};
\addlegendentry{2-layer SPD-NN}
\addplot [very thick, blue, dashed]
table {%
3.94689503843385e-05 0.00789378834391847
0.000262733955566399 0.0525475193659454
0.000481065986240417 0.0962139348306551
0.000695825683113496 0.139169880283174
0.000906118329342109 0.181230787958491
0.0011107337942216 0.222164325612008
0.00130848212280662 0.261704128266444
0.00149806442489164 0.299611929718142
0.00323118478202684 0.334629522768132
0.00498187339422014 0.369634489902273
0.00641989926962261 0.398396502555753
0.00797307823638301 0.429457636058107
0.00914739630832405 0.452943832791969
0.0104060085632686 0.478116374517017
0.0113337223081665 0.496673514608576
0.0122478967061566 0.514957002689553
0.0129161592117776 0.528323599376741
0.0134716701144671 0.539434044947705
0.013854007364762 0.547078308976729
0.0140548643385603 0.551102387234913
0.0141270319258187 0.552161967951817
0.0141152176633177 0.549834632535187
0.0140863444788218 0.544065850529339
0.014040810165613 0.534959161892728
0.0139788120361947 0.522563023406231
0.0139007432396196 0.506947150770901
0.0138071329724081 0.488219732284366
0.0136984236885573 0.466483405614551
0.0135753597727287 0.441871187419702
0.0134386247047169 0.414525937907431
0.0132893567623637 0.384637947547257
0.0131280296343113 0.352365147216158
0.0129558344897333 0.317914291754309
0.0127737647826751 0.281492761481502
0.0125829273049948 0.243324416732479
0.0123845650519812 0.203645065480671
0.0121798379684994 0.162681800772256
0.0119700903043993 0.120722723966394
0.0117564722982326 0.0779927051182828
0.0115403355155068 0.0347620662804726
};
\addlegendentry{3-layer SPD-NN}
\addplot [very thick, green!50!black, dashed]
table {%
3.94689687357629e-05 0.00789379101006046
0.000262733982946356 0.0525475193288907
0.00048106784478643 0.0962139158666054
0.000695863238952516 0.139169985052512
0.000906185736134471 0.181230840614616
0.00111081841484354 0.222164987829134
0.00130859547901673 0.261706627796227
0.0014981409656174 0.299615818548169
0.00324242757383932 0.334629032875664
0.00499374752493604 0.36961819028247
0.00643452691697689 0.398398176973487
0.00798792675485959 0.429446118199273
0.00916402584775969 0.452939804412613
0.01042341827316 0.478108699744784
0.0113524850315964 0.496663409639055
0.0122679512205035 0.514949324223006
0.0129365423042325 0.528312771812314
0.0134914906727303 0.539423822514049
0.0138728318656146 0.547063914436263
0.0139473334583165 0.549048220799923
0.0141199643798847 0.552496070582387
0.0141207049488251 0.551413969550072
0.014100425395682 0.5473512309099
0.0140632836862833 0.539920600632693
0.0140095375819849 0.529168852406695
0.0139394807497031 0.515155052596288
0.0138535549000971 0.497961402693289
0.0137524470282083 0.477722902731981
0.0136365375902008 0.454530885165868
0.0135065836390329 0.428534862807534
0.0133634071907397 0.399897188820989
0.0132079053468865 0.36879900028096
0.0130408793356927 0.335399450498303
0.0128635841271355 0.299929673989193
0.0126769862123208 0.262617749293221
0.0124822253197673 0.223660421066899
0.0122805025629049 0.183313480948412
0.0120731786097025 0.141826004798615
0.0118612893971518 0.099438421024087
0.0116462319567873 0.0564251934753326
0.0114293385203493 0.0130409982637403
};
\addlegendentry{4-layer SPD-NN}
\end{axis}

\end{tikzpicture}}~
\scalebox{0.6}{
\begin{tikzpicture}

\definecolor{color0}{rgb}{0.12156862745098,0.466666666666667,0.705882352941177}
\definecolor{color1}{rgb}{0.75,0.75,0}

\begin{axis}[
scaled ticks=false,
legend cell align={left},
legend style={fill opacity=0.8, draw opacity=1, text opacity=1, at={(0.5,0.09)}, anchor=south, draw=white!80!black},
tick align=outside,
tick pos=left,
x grid style={white!69.0196078431373!black},
xlabel={Strain},
xmin=-0.00125, xmax=0.015,
xtick style={color=black},
xtick={-0.005,0,0.005,0.01,0.015},
xticklabels={-0.005,0.000,0.005,0.010,0.015},
y grid style={white!69.0196078431373!black},
ylabel={Stress (GPa)},
ymin=-0.1, ymax=0.6,
ytick style={color=black},
ytick={-0.1,0,0.1,0.2,0.3,0.4,0.5,0.6},
yticklabels={-0.1,0.0,0.1,0.2,0.3,0.4,0.5,0.6}
]
\addplot [very thick, color0, mark=*, mark size=3, mark options={solid,fill opacity=0}, only marks]
table {%
3.94667614855144e-05 0.00789335229710289
0.000262737609167358 0.0525475218334735
0.000481069672157902 0.0962139343718435
0.000695845603162529 0.139169120551953
0.000906163111303343 0.181232622179406
0.00111080123165105 0.222160246249025
0.00130851635936997 0.261703271792929
0.00149809974672009 0.299619949263072
0.00323193202914041 0.334638640574785
0.00498143865949836 0.369628773181944
0.00642003916262919 0.398400783244561
0.00797281755091166 0.429456351001532
0.00914725186338931 0.452945037251085
0.0104058499942009 0.478116999823974
0.0113336195708573 0.496672391260666
0.0122478758536451 0.514957516869412
0.0129161596962976 0.528323193692812
0.0134717550260704 0.539435100075987
0.0138536915245349 0.547073829857545
0.014056055634465 0.551121111969057
0.0141307531614943 0.55216874558956
0.0141190622100165 0.549830555294001
0.014090225766778 0.544063266593623
0.0140446899430925 0.534956101855122
0.0139827141887345 0.522560950983515
0.0139046646260146 0.506951038439551
0.0138110123627919 0.488220585795028
0.0137023257792617 0.46648326908901
0.0135792655343288 0.441871220102479
0.0134425802921606 0.414534171668899
0.0132931022556121 0.384638564359259
0.013131742271295 0.352366567495903
0.0129594844712875 0.317915007494492
0.0127773804659743 0.281494206431959
0.0125865431151509 0.243326736267379
0.0123881399091834 0.203646095073998
0.0121833859950344 0.162695312244315
0.0119735368846105 0.120725490159676
0.0117598808853918 0.0779942903160686
0.0115437312956711 0.0347643723720768
};
\addlegendentry{Reference}
\addplot [very thick, red, dashed]
table {%
4.08889624748808e-05 0.00787229595007668
0.000260053883461613 0.0525172466562051
0.000481855936481852 0.0962132249979262
0.000698348791393998 0.139172691572173
0.000904480061668871 0.181230401267596
0.00110631425046904 0.222153005546579
0.00132280977389184 0.261645598237538
0.00171627953074338 0.299018813932196
0.00323110510077593 0.334501482807902
0.00496371992527749 0.369503263277342
0.00642389044734166 0.398748375288848
0.00793909184567857 0.429023685407099
0.00916037662295887 0.453408697869779
0.0103764500120811 0.477686735749407
0.0113473656373139 0.497019035096363
0.0122339749348333 0.51470980843387
0.0129307050160776 0.528434021672707
0.013494387635186 0.539421754677139
0.013902361149432 0.5469136635176
0.0141765986438545 0.551238059268273
0.0143254085591518 0.552086442928967
0.0143833946223522 0.549771432006098
0.014368175465508 0.543843619967535
0.0143245034952905 0.534859837699867
0.0142648226021896 0.522402040941574
0.0142055315527012 0.506774190113394
0.0141463831984271 0.488066580498613
0.0140792456141835 0.466317664816347
0.0139949569479611 0.441690125257196
0.0138859500824686 0.414364620924282
0.0137473224771247 0.384481646641895
0.0135835790599724 0.352209472230411
0.0134125515872069 0.317752179831315
0.01328508256817 0.28127255606018
0.0132944624002614 0.24303719767732
0.0134163235309526 0.203502307117742
0.0134419808354318 0.162640848403296
0.0133451789330687 0.120474026626505
0.0132250146622628 0.0779795549439368
0.0130719246005625 0.0346731462987797
};
\addlegendentry{1-layer $\bm{\sigma}$-NN}
\addplot [very thick, color1, dashed]
table {%
3.82976251778017e-05 0.00793757286745595
0.000264166425435856 0.0525475458987418
0.000480112749048315 0.0962122278783488
0.000694644796803702 0.139167802974339
0.00090751348051118 0.181249630202735
0.00111099069102603 0.222162220201086
0.00130723594001461 0.26168047967459
0.00165070787650355 0.298842245619863
0.00324281273350691 0.334735461507474
0.0049766952912733 0.369293337887303
0.00644284081209161 0.398836458278485
0.00793625993742898 0.428988535936614
0.00913460090864077 0.453416683920911
0.0103181632235266 0.477772692645043
0.0112442423049156 0.496993736341124
0.0120956015814637 0.514875826731667
0.0127275792368023 0.528473519760025
0.0132225993817813 0.539620557572589
0.0135164388978424 0.547153405576474
0.0136131418940947 0.551536914537229
0.0134633474062292 0.55277772006749
0.0128859555927634 0.550722421817506
0.0100087440594551 0.54951684243011
-0.0036432601727352 0.54408231440488
0.00000305433       0.547379246496647
0.0718556994751697 0.192721910624124
0.330226597662544 0.35625132190084
0.563061340345363 0.127173789712946
};
\addlegendentry{2-layer $\bm{\sigma}$-NN}
\addplot [very thick, blue, dashed]
table {%
3.87845306912995e-05 0.00794055296715892
0.000265720438368063 0.0525568002802273
0.000482566164474645 0.0962157609840971
0.000695021915821252 0.139163328163794
0.000908981125402213 0.181234118745827
0.00111747794614385 0.222168230368504
0.00132025718992501 0.261696948267393
0.00179399158360446 0.29944446591706
0.00334231660723389 0.336757602463651
0.00494604102652229 0.369077919681195
0.00643334973281932 0.399519084914209
0.00787184135364818 0.428302405620054
0.00915732900013565 0.45385288796198
0.0103176224083028 0.47702833056697
0.0113132485232906 0.497023863305966
0.0121626662484179 0.51413248817189
0.0128563530288213 0.528051760135191
0.0134026246420867 0.53883373912893
0.0138075621088394 0.54643645596914
0.0140830274522752 0.550868359845977
0.0142446297950874 0.552095415965308
0.0143122118121912 0.550025672167397
0.0143103024549297 0.544519484289096
0.0142670848067356 0.535457696264125
0.0142133290432334 0.522944851940063
0.0141688584356404 0.507255337651486
0.0141263803890984 0.48848108725956
0.01406806702559 0.466710770101275
0.0139799789996134 0.442089228933051
0.013857643041262 0.414761979982276
0.0137071999469 0.384900798017167
0.0135431287707213 0.35269829001891
0.0133825644946387 0.318369301770498
0.013239765175469 0.282174389927295
0.013120381865428 0.244440247241689
0.0130218341672678 0.205587577421154
0.0129360205126732 0.166017405706795
0.0128380423406705 0.125646177467344
0.0126998686613024 0.0839244494115147
0.0125056529296128 0.0402460370008097
};
\addlegendentry{3-layer $\bm{\sigma}$-NN}
\addplot [very thick, green!50!black, dashed]
table {%
4.005949765301e-05 0.00787422558920577
0.000348002307006362 0.0518468239712947
0.000847560391861523 0.0961972555628921
0.00103539489807573 0.139124878089104
0.00117315371769126 0.181174786755068
0.00131661301250579 0.222098435839601
0.00148697218089841 0.261596349136245
0.00181144347138174 0.299115336172889
0.00322354465256289 0.334154802185547
0.00501124788801422 0.369792892651595
0.00646636239670306 0.398577066156971
0.0079712580890547 0.429078280458412
0.00914475808586558 0.453446155476879
0.0103119964639373 0.47759997514697
0.011271391265729 0.497187890282223
0.0121404422988902 0.514602618460952
0.0128617272250241 0.528604482227535
0.0134415805059251 0.539313326788688
0.0139021353650075 0.547014162270236
0.014231463041999 0.551102699263272
0.0144704271423328 0.552067012302529
0.014630788346531 0.549586941139085
0.0147365887925324 0.543658986524844
0.0148122311431795 0.534643140248579
0.0148528686650691 0.522081947720016
0.0148680734447409 0.506480835323774
0.0148517942641581 0.48775399460841
0.0148027931663497 0.465993341824508
0.014723405578574 0.441382922015296
0.0146175177327654 0.414063637629792
0.0144907737603115 0.384194957601069
0.0143482440672404 0.351957733539112
0.0141933204127697 0.317534851177743
0.0140277700251545 0.281155691614053
0.0138508259598543 0.243032415294352
0.0136614596334062 0.203401597179502
0.0134594235404786 0.16249490921423
0.0132468486575854 0.120579301002968
0.013027794743737 0.0778870083753143
0.012808134548827 0.0347189414370815
};
\addlegendentry{4-layer $\bm{\sigma}$-NN}
\end{axis}

\end{tikzpicture}}~
\scalebox{0.6}{
\begin{tikzpicture}

\definecolor{color0}{rgb}{0.12156862745098,0.466666666666667,0.705882352941177}
\definecolor{color1}{rgb}{0.75,0.75,0}

\begin{axis}[
scaled ticks=false,
legend cell align={left},
legend style={fill opacity=0.8, draw opacity=1, text opacity=1, at={(0.5,0.09)}, anchor=south, draw=white!80!black},
tick align=outside,
tick pos=left,
x grid style={white!69.0196078431373!black},
xlabel={Strain},
xmin=-0.00125, xmax=0.015,
xtick style={color=black},
xtick={-0.005,0,0.005,0.01,0.015},
xticklabels={-0.005,0.000,0.005,0.010,0.015},
y grid style={white!69.0196078431373!black},
ylabel={Stress (GPa)},
ymin=-0.1, ymax=0.6,
ytick style={color=black},
ytick={-0.1,0,0.1,0.2,0.3,0.4,0.5,0.6},
yticklabels={-0.1,0.0,0.1,0.2,0.3,0.4,0.5,0.6}
]
\addplot [very thick, color0, mark=*, mark size=3, mark options={solid,fill opacity=0}, only marks]
table {%
3.94667614855144e-05 0.00789335229710289
0.000262737609167358 0.0525475218334735
0.000481069672157902 0.0962139343718435
0.000695845603162529 0.139169120551953
0.000906163111303343 0.181232622179406
0.00111080123165105 0.222160246249025
0.00130851635936997 0.261703271792929
0.00149809974672009 0.299619949263072
0.00323193202914041 0.334638640574785
0.00498143865949836 0.369628773181944
0.00642003916262919 0.398400783244561
0.00797281755091166 0.429456351001532
0.00914725186338931 0.452945037251085
0.0104058499942009 0.478116999823974
0.0113336195708573 0.496672391260666
0.0122478758536451 0.514957516869412
0.0129161596962976 0.528323193692812
0.0134717550260704 0.539435100075987
0.0138536915245349 0.547073829857545
0.014056055634465 0.551121111969057
0.0141307531614943 0.55216874558956
0.0141190622100165 0.549830555294001
0.014090225766778 0.544063266593623
0.0140446899430925 0.534956101855122
0.0139827141887345 0.522560950983515
0.0139046646260146 0.506951038439551
0.0138110123627919 0.488220585795028
0.0137023257792617 0.46648326908901
0.0135792655343288 0.441871220102479
0.0134425802921606 0.414534171668899
0.0132931022556121 0.384638564359259
0.013131742271295 0.352366567495903
0.0129594844712875 0.317915007494492
0.0127773804659743 0.281494206431959
0.0125865431151509 0.243326736267379
0.0123881399091834 0.203646095073998
0.0121833859950344 0.162695312244315
0.0119735368846105 0.120725490159676
0.0117598808853918 0.0779942903160686
0.0115437312956711 0.0347643723720768
};
\addlegendentry{Reference}
\addplot [very thick, red, dashed]
table {%
3.95130021758227e-05 0.00788875367795907
0.000262712626113764 0.052550516572439
0.000480322762818513 0.0962128010491301
0.000695614631565621 0.139168772632778
0.000905741091379289 0.181235064561868
0.0011101155360173 0.222154912387473
0.00130825522341318 0.261713068126932
0.00152135405472568 0.299320783348775
0.00327168048939588 0.334629920795616
0.00505065637152489 0.369602414595446
0.00654157870696624 0.398335723902335
0.00818438448563922 0.429349330199244
0.00949680536012493 0.452760510034721
0.0109735734887406 0.477774767751528
0.0122451108746564 0.49616544725575
0.0136837206923538 0.514071719952102
0.0151501715366699 0.527057309697663
0.0169207973938427 0.537051455102477
0.0194732830315932 0.541812563726587
0.0284421478117495 0.536400029442608
};
\addlegendentry{1-layer $\Delta\bm{\sigma}$-NN}
\addplot [very thick, color1, dashed]
table {%
3.93187256084496e-05 0.00789477601580551
0.000262573794056325 0.0525487074924115
0.00048088205316576 0.0962143028903345
0.000695864586298526 0.13916949674117
0.000906085540472365 0.181233055060841
0.00111113066049851 0.222167237206687
0.00130813796685753 0.261690000914166
0.00150639077467955 0.299395292882374
0.00324329116412875 0.334612302851805
0.00498334135776818 0.369628587914087
0.00640995816988962 0.398418564387147
0.00794481001285317 0.429470181758404
0.00909630319336443 0.452977130206087
0.0103159331351053 0.47819749052994
0.0111702271122084 0.496825212214943
0.0119435914847012 0.515508692826163
0.0125659097624697 0.535604613825834
0.014118871124149 0.575592242861556
-0.0011789299242605 0.489460962071618
-0.000581565354028127 0.520507863704395
-0.00117877203995308 0.548537643368938
-0.00117597969419723 0.539332580904296
-0.00112671517214561 0.564414647235404
-0.00135778265470918 0.569624733814876
-0.00116529588679732 0.544324553830447
-0.00124805149184246 0.51715872374629
-0.00108901341441195 0.522848447312822
-0.00108381806200877 0.513890924526039
-0.00113446720568434 0.473265665880031
-0.000928653507988639 0.453291949264466
-0.000906904600714401 0.438287228919384
-0.000722607523678767 0.386441843747322
-0.000626088155948692 0.365035264839322
-0.000455473270717177 0.312786417351218
-0.000352544054590061 0.279926969569764
-0.000265523440803925 0.229687604370844
-0.000351294040941064 0.174934441122823
-0.000455583202181006 0.120070804341369
-0.000391590075317839 0.0605401393725636
-0.00004805093773767 -0.0038425107578465
-0.0405255178123943 -0.433175369620508
-0.407056439403237 -0.882050299871246
};
\addlegendentry{2-layer $\Delta\bm{\sigma}$-NN}
\addplot [very thick, blue, dashed]
table {%
3.95692919945585e-05 0.00789389837268052
0.000262792269393327 0.0525464058293176
0.000481110358957208 0.0962135403716172
0.000695932903888241 0.139168988186463
0.000906252230509066 0.181231295300467
0.00111082601154209 0.22216102630998
0.00130851547532516 0.261720226915202
0.00150049929495619 0.299539691743276
0.00323152014757114 0.33463376600032
0.0049630991516277 0.369641088141056
0.00636400928943928 0.398446351518384
0.00783104889549769 0.429537750742566
0.00880610922424098 0.453161004038061
0.00946822739433056 0.478827968423915
0.0077009359961647 0.499878513198204
0.0029634925357734 0.523849473566502
0.0025426133955674 0.572540041021841
0.00355174280385212 0.414465749527953
0.00259315261359579 0.607118918989053
0.00228552064945204 0.497292199235546
0.000545313007851565 0.275963904262393
0.050496679928726 0.689464938077412
0.0186449159344137 0.585745202555795
0.0128784793022557 0.43024612061831
0.0127901827226839 0.40999022207589
0.0175186967374931 0.532905717116113
0.0188708373982143 0.553190255575055
0.0215402773799861 0.620295329074573
0.0103921529923891 0.37823113866625
0.015940245921105 0.526634358155159
};
\addlegendentry{3-layer $\Delta\bm{\sigma}$-NN}
\addplot [very thick, green!50!black, dashed]
table {%
3.93493215735238e-05 0.00788746663985415
0.000262695550227629 0.052547654664549
0.00048108666360434 0.0962144299211311
0.000695824842237532 0.139168653531705
0.000906154213240006 0.181234620316752
0.00111076533593697 0.222160383074112
0.00130850903752255 0.261704963245932
0.00149894784303355 0.299592414128628
0.00200281022875173 0.308073389735921
0.00371199791598534 0.336290879709034
0.00792768239238414 0.3681493903945
0.109929014788497 0.332596790291565
};
\addlegendentry{4-layer $\Delta\bm{\sigma}$-NN}
\end{axis}

\end{tikzpicture}}
\caption{The strain stress curve of the NN-FEM test results obtained by using SPD-NN~(left), $\bm{\sigma}$-NN~(middle), and $\Delta\bm{\sigma}$-NN~(right).}
\label{FIG:TRUSS-NN-FEM-TEST-STRESS}
\end{figure}
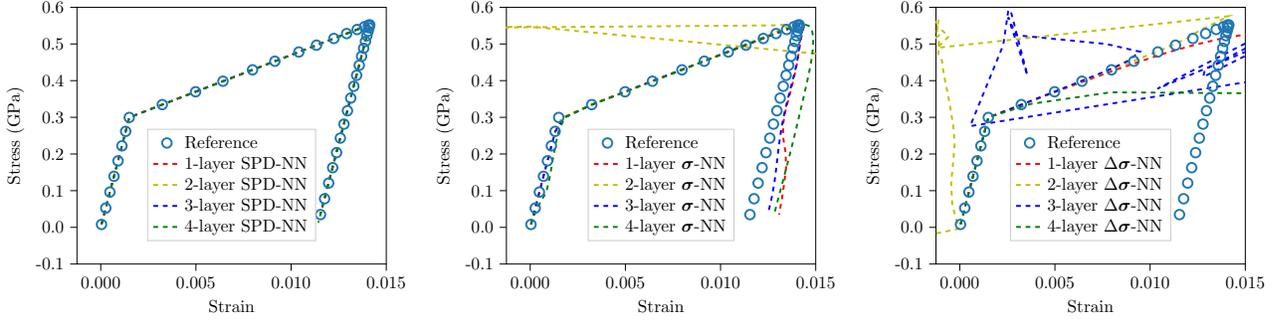

For the NN test, the predicted strain-stress relations are reported in \cref{FIG:TRUSS-NN-TEST}. All kinds of NN architectures give good results.
For the NN-FEM test, the predicted displacement trajectories of the right end point and the predicted strain-stress curves for one Gaussian quadrature point on the right end element are reported in \cref{FIG:TRUSS-NN-FEM-TEST-DISP} and  \cref{FIG:TRUSS-NN-FEM-TEST-STRESS}.
Most of $\bm{\sigma}$-NN and $\Delta\bm{\sigma}$-NN architectures are found to be numerically unstable. 
We attribute the instability to the deviation of the strain and stress pairs in the test loading condition from the training set, where the numerical errors are accumulated during the prediction process. More specifically, the \textbf{violation of SPD constraints} on the tangent stiffness matrix for $\bm{\sigma}$-NN or $\Delta\bm{\sigma}$-NN makes them vulnerable to these errors and less robust to even slight extrapolations.
In contrast, the proposed SPD-NN delivers numerically stable results and the predicted displacements and strain-stress curves overlap the exact ones.

\begin{figure}[htpb]
\scalebox{0.6}{
\begin{tikzpicture}

\definecolor{color0}{rgb}{0.12156862745098,0.466666666666667,0.705882352941177}
\definecolor{color1}{rgb}{0.75,0.75,0}

\begin{axis}[
scaled ticks=false,
legend cell align={left},
legend style={fill opacity=0.8, draw opacity=1, text opacity=1, at={(0.97,0.03)}, anchor=south east, draw=white!80!black},
tick align=outside,
tick pos=left,
x grid style={white!69.0196078431373!black},
xlabel={Time (s)},
xmin=-0.01, xmax=0.21,
xtick style={color=black},
xtick={-0.05,0,0.05,0.1,0.15,0.2,0.25},
xticklabels={-0.05,0.00,0.05,0.10,0.15,0.20,0.25},
y grid style={white!69.0196078431373!black},
ylabel={Displacement (m)},
ymin=-0.005, ymax=0.06,
ytick style={color=black},
ytick={-0.01,0,0.01,0.02,0.03,0.04,0.05,0.06},
yticklabels={-0.01,0.00,0.01,0.02,0.03,0.04,0.05,0.06}
]
\addplot [very thick, color0, mark=*, mark size=3, mark options={solid,fill opacity=0}, only marks]
table {%
0 0
0.005 0.000879071318966359
0.01 0.00175178603128696
0.015 0.00261242457318852
0.02 0.00345675778867308
0.025 0.00427942738576081
0.03 0.0050753033163228
0.035 0.00583950109613968
0.04 0.0106107839172756
0.045 0.019541230444874
0.05 0.0236237531592033
0.055 0.0313289373424314
0.06 0.0350380478185725
0.065 0.0408174614436271
0.07 0.0442730136401229
0.075 0.0480651190399302
0.08 0.0509821771299189
0.085 0.053019447971144
0.09 0.0550077364727381
0.095 0.0555728014239678
0.1 0.0563137678270884
0.105 0.056260502065091
0.11 0.0561555202057213
0.115 0.0559882518807527
0.12 0.0557565179273646
0.125 0.055461132096511
0.13 0.0551037742624527
0.135 0.0546865831978064
0.14 0.0542120745133923
0.145 0.053683114752886
0.15 0.0531029032911068
0.155 0.0524749538344366
0.16 0.0518030742630135
0.165 0.0510913446992788
0.17 0.0503440938727484
0.175 0.0495658738865951
0.18 0.0487614335058666
0.185 0.0479356900987087
0.19 0.0470937003732258
0.195 0.0462406300637518
0.2 0.0453817227311881
};
\addlegendentry{Reference}
\addplot [very thick, red, dashed]
table {%
0 0
0.005 0.000879078401218516
0.01 0.00175179752124534
0.015 0.00261243856778289
0.02 0.00345641319513846
0.025 0.00427923454464217
0.03 0.00507493901990294
0.035 0.00583931345891475
0.04 0.0105766759051313
0.045 0.0195182597228965
0.05 0.0235977998938799
0.055 0.0313008153608988
0.06 0.0350168344823498
0.065 0.0407875480351868
0.07 0.0442530673108374
0.075 0.0480360713418296
0.08 0.0509613797625199
0.085 0.0529923455093924
0.09 0.0549852668506537
0.095 0.0555199346446133
0.1 0.0562671805578271
0.105 0.0562044105853059
0.11 0.0560985608316154
0.115 0.0559307715925234
0.12 0.055699010877408
0.125 0.0554039765462809
0.13 0.0550464657359777
0.135 0.0546295412859708
0.14 0.0541550463172466
0.145 0.0536261501183365
0.15 0.0530456542871657
0.155 0.0524173986415327
0.16 0.0517453147140339
0.165 0.0510332859144629
0.17 0.0502858557412443
0.175 0.0495075750710385
0.18 0.048703127487631
0.185 0.047877412619905
0.19 0.0470357487990184
0.195 0.0461831494118014
0.2 0.0453243294285531
};
\addlegendentry{1-layer SPD-NN}
\addplot [very thick, color1, dashed]
table {%
0 0
0.005 0.000879061538242488
0.01 0.00175177616526263
0.015 0.0026124143154725
0.02 0.00345642168976587
0.025 0.00427927515005798
0.03 0.00507501437052192
0.035 0.00583938965858001
0.04 0.0105780964061268
0.045 0.0195169620643884
0.05 0.0235985305004029
0.055 0.031299369713911
0.06 0.0350164963588067
0.065 0.0407863868922361
0.07 0.0442517989041527
0.075 0.04803502030439
0.08 0.0509595545943314
0.085 0.0529912233361132
0.09 0.0549833138110029
0.095 0.0555478417331359
0.1 0.0562918120275845
0.105 0.0562377014954141
0.11 0.0561326089763396
0.115 0.0559654035976417
0.12 0.05573428833152
0.125 0.0554388415333137
0.13 0.0550812886853278
0.135 0.0546638263329486
0.14 0.0541890040738581
0.145 0.0536602906561461
0.15 0.053079954527676
0.155 0.0524521389438237
0.16 0.0517805476623829
0.165 0.0510688762078209
0.17 0.050321815854599
0.175 0.049543769757971
0.18 0.0487396848138011
0.185 0.0479143302871915
0.19 0.0470728554406128
0.195 0.0462203000983228
0.2 0.0453614454582803
};
\addlegendentry{2-layer SPD-NN}
\addplot [very thick, blue, dashed]
table {%
0 0
0.005 0.000879062690600044
0.01 0.00175177728020032
0.015 0.00261241104928193
0.02 0.00345641193587466
0.025 0.0042792589108844
0.03 0.00507499078107346
0.035 0.0058393620825063
0.04 0.0106062937618303
0.045 0.0195451878826163
0.05 0.023621300806923
0.055 0.0313292058289549
0.06 0.0350386686014015
0.065 0.0408171766504609
0.07 0.0442739746918498
0.075 0.0480646897602867
0.08 0.0509830608059845
0.085 0.0530183293492297
0.09 0.0550086196171651
0.095 0.055558606536534
0.1 0.0563180206486562
0.105 0.0562632645283923
0.11 0.0561580024990627
0.115 0.0559906942191262
0.12 0.0557590078515599
0.125 0.0554634260291996
0.13 0.0551058917269182
0.135 0.0546883263531129
0.14 0.0542134636304129
0.145 0.0536842321881547
0.15 0.053104655977548
0.155 0.0524769174638035
0.16 0.0518051102336995
0.165 0.051093207068687
0.17 0.0503458993713291
0.175 0.0495676832276342
0.18 0.0487632893190593
0.185 0.0479379248625543
0.19 0.0470960500454008
0.195 0.0462429135494835
0.2 0.0453839466901184
};
\addlegendentry{3-layer SPD-NN}
\addplot [very thick, green!50!black, dashed]
table {%
0 0
0.005 0.000879062758800149
0.01 0.00175177959624474
0.015 0.00261253993115254
0.02 0.00345666582619255
0.025 0.00427959210686884
0.03 0.00507539643499307
0.035 0.00583963003945318
0.04 0.0106237066809708
0.045 0.019554400081667
0.05 0.0236461407628962
0.055 0.0313500200472626
0.06 0.0350631516286393
0.065 0.0408470644666885
0.07 0.0443015633669648
0.075 0.0480995777896689
0.08 0.0510147369600692
0.085 0.053050535550267
0.09 0.0550381499609952
0.093 0.0552217780230334
0.0975 0.056103153925734
0.1025 0.0562497489165849
0.1075 0.0561967414468371
0.1125 0.0560657528396898
0.1175 0.0558671379853565
0.1225 0.0556032859554117
0.1275 0.0552768445450757
0.1325 0.054889837723282
0.1375 0.0544439367436295
0.1425 0.0539417145636868
0.1475 0.0533864207677814
0.1525 0.0527816350762976
0.1575 0.0521312815175881
0.1625 0.051438966497999
0.1675 0.0507087114772472
0.1725 0.0499456299154162
0.1775 0.0491538934134204
0.1825 0.048338360029485
0.1875 0.0475042307025414
0.1925 0.0466564157693146
0.1975 0.0457997432187567
};
\addlegendentry{4-layer SPD-NN}
\end{axis}

\end{tikzpicture}}~
\scalebox{0.6}{
\begin{tikzpicture}

\definecolor{color0}{rgb}{0.12156862745098,0.466666666666667,0.705882352941177}
\definecolor{color1}{rgb}{0.75,0.75,0}

\begin{axis}[
scaled ticks=false,
legend cell align={left},
legend style={fill opacity=0.8, draw opacity=1, text opacity=1, at={(0.97,0.03)}, anchor=south east, draw=white!80!black},
tick align=outside,
tick pos=left,
x grid style={white!69.0196078431373!black},
xlabel={Time (s)},
xmin=-0.01, xmax=0.21,
xtick style={color=black},
xtick={-0.05,0,0.05,0.1,0.15,0.2,0.25},
xticklabels={-0.05,0.00,0.05,0.10,0.15,0.20,0.25},
y grid style={white!69.0196078431373!black},
ylabel={Displacement (m)},
ymin=-0.005, ymax=0.06,
ytick style={color=black},
ytick={-0.01,0,0.01,0.02,0.03,0.04,0.05,0.06},
yticklabels={-0.01,0.00,0.01,0.02,0.03,0.04,0.05,0.06}
]
\addplot [very thick, color0, mark=*, mark size=3, mark options={solid,fill opacity=0}, only marks]
table {%
0 0
0.005 0.000879071318966359
0.01 0.00175178603128696
0.015 0.00261242457318852
0.02 0.00345675778867308
0.025 0.00427942738576081
0.03 0.0050753033163228
0.035 0.00583950109613968
0.04 0.0106107839172756
0.045 0.019541230444874
0.05 0.0236237531592033
0.055 0.0313289373424314
0.06 0.0350380478185725
0.065 0.0408174614436271
0.07 0.0442730136401229
0.075 0.0480651190399302
0.08 0.0509821771299189
0.085 0.053019447971144
0.09 0.0550077364727381
0.095 0.0555728014239678
0.1 0.0563137678270884
0.105 0.056260502065091
0.11 0.0561555202057213
0.115 0.0559882518807527
0.12 0.0557565179273646
0.125 0.055461132096511
0.13 0.0551037742624527
0.135 0.0546865831978064
0.14 0.0542120745133923
0.145 0.053683114752886
0.15 0.0531029032911068
0.155 0.0524749538344366
0.16 0.0518030742630135
0.165 0.0510913446992788
0.17 0.0503440938727484
0.175 0.0495658738865951
0.18 0.0487614335058666
0.185 0.0479356900987087
0.19 0.0470937003732258
0.195 0.0462406300637518
0.2 0.0453817227311881
};
\addlegendentry{Reference}
\addplot [very thick, red, dashed]
table {%
0 0
0.005 0.000870923083515745
0.01 0.00175110987876382
0.015 0.00262388826833368
0.02 0.00345395720964754
0.025 0.00426050381856151
0.03 0.0050974132309487
0.035 0.00630411179888409
0.04 0.0108824168476538
0.045 0.0191334648089029
0.05 0.0239898301811749
0.055 0.0308903855334429
0.06 0.0353425469428658
0.065 0.0405281079457713
0.07 0.0444363724898835
0.075 0.0479847359458777
0.08 0.0510004518020758
0.085 0.0531885356218133
0.09 0.0550493963619537
0.095 0.0561265731533079
0.1 0.0568716476992647
0.105 0.0571278153215153
0.11 0.0571134739426435
0.115 0.0569768222210431
0.12 0.0567291579710444
0.125 0.0565046632894853
0.13 0.0562771547384909
0.135 0.0560216094751349
0.14 0.0557077124102496
0.145 0.0553001659151164
0.15 0.0547785321214855
0.155 0.0541480323024984
0.16 0.0534682843201617
0.165 0.0528862191368136
0.17 0.0527537885018381
0.175 0.0532924726408608
0.18 0.0536750864693818
0.185 0.0532397991570891
0.19 0.0528392830100497
0.195 0.0522277478048262
0.2 0.051586438267376
};
\addlegendentry{1-layer $\bm{\sigma}$-NN}
\addplot [very thick, color1, dashed]
table {%
0 0
0.005 0.000887542864332501
0.01 0.00175040667211076
0.015 0.00260614840895145
0.02 0.00346031497993567
0.025 0.0042829037103671
0.03 0.00506580321448695
0.035 0.00605544090317005
0.04 0.0110156174558964
0.045 0.0191535866042352
0.05 0.0240524610875334
0.055 0.0309242026270479
0.06 0.0352060891713938
0.065 0.0403960603371704
0.07 0.0440018073398862
0.075 0.0475500835226629
0.08 0.0502492293489506
0.085 0.0522499100365265
0.09 0.0537456442559479
0.095 0.0542040295736954
0.1 0.0541464779819701
0.105 0.0525445534780952
0.11 0.0484914817841192
0.115 0.0243370225731414
0.11775 0.00403294706336368
0.118375 0.0290419575878821
0.1195 0.167744318627946
0.12025 0.303217146035237
0.121375 0.578123776605413
0.123 0.957882045226973
0.12375 1.01567498951727
0.1255 1.06289593175781
0.1285 1.33394423812031
0.1335 1.80658018825285
0.1385 2.20578172306239
0.1435 2.65322125774372
0.1465 3.01345819052708
0.14875 3.42611276095451
0.15175 4.01644820076663
0.15475 4.05882762050406
0.15875 3.60887226752663
0.16375 2.85720701254155
0.16875 2.12626987368083
0.171 1.65667906857169
0.173 0.960158706416766
0.17525 0.294493944134351
0.17625 0.0366410121699136
0.177375 -0.226601700017042
0.179625 -0.623070680256852
0.181125 -0.796557246519802
0.183625 -1.02740985702912
0.185125 -1.18402676056051
0.187375 -1.36879988253216
0.188625 -1.43977089997284
0.190875 -1.58942840899353
0.195875 -2.13495704962868
0.197875 -2.45023205065117
};
\addlegendentry{2-layer $\bm{\sigma}$-NN}
\addplot [very thick, blue, dashed]
table {%
0 0
0.005 0.000889454628766658
0.01 0.00175740549521075
0.015 0.00260301888089238
0.02 0.0034524282700439
0.025 0.00428133969297919
0.03 0.00503955797824655
0.035 0.00587398891104507
0.04 0.00560182886040498
0.045 -0.00267757478564986
0.05 -0.0145738497916282
0.055 -0.0272169456391506
0.06 -0.04043206528308
0.065 -0.0532034909112403
0.07 -0.0652679814043808
0.075 -0.0759057065164164
0.08 -0.0845974306293142
0.085 -0.0907083987202652
0.09 -0.093802225559081
0.095 -0.0939249040543532
0.1 -0.0919028514052585
0.105 -0.0893923554776519
0.11 -0.0884888888509743
0.115 -0.0909290763681326
0.12 -0.0972655093079304
0.125 -0.107001220828857
0.13 -0.119636754223258
0.135 -0.134965868058347
0.14 -0.152874159180811
0.145 -0.173315436900081
0.15 -0.196308762901637
0.155 -0.221974558333073
0.16 -0.250619751130913
0.165 -0.282888545478805
0.17 -0.320065304855019
0.175 -0.364667749216237
0.18 -0.421473406352163
0.185 -0.498394169389592
0.19 -0.604976404219031
0.195 -0.747949982607296
0.2 -0.925995237723045
};
\addlegendentry{3-layer $\bm{\sigma}$-NN}
\addplot [very thick, green!50!black, dashed]
table {%
0 0
0.005 0.00104581948001214
0.01 0.00319977231311563
0.015 0.00403207098394598
0.02 0.00458331285239867
0.025 0.00514411166638543
0.03 0.00578798815282637
0.035 0.00682511278627297
0.04 0.0106954712633054
0.045 0.0192646644107229
0.05 0.0242355308347113
0.055 0.0309410615378751
0.06 0.0354569131721403
0.065 0.0401883911434918
0.07 0.0442767139619851
0.075 0.0475489035172016
0.08 0.0508174918958074
0.085 0.052936557056539
0.09 0.0551048482055555
0.095 0.0563078109707412
0.1 0.0574755418421484
0.105 0.0580777334955774
0.11 0.0585901432215657
0.115 0.0589186879136113
0.12 0.0590933365039828
0.125 0.059198109068947
0.13 0.0591570879971247
0.135 0.0589876903184528
0.14 0.0586984563454094
0.145 0.0583006561094183
0.15 0.0578152300891613
0.155 0.0572640961561306
0.16 0.0566623044892784
0.165 0.0560174698517329
0.17 0.0553283643607624
0.175 0.0545905691081023
0.18 0.0538025207623086
0.185 0.0529696133013274
0.19 0.0521075864213888
0.195 0.0512381490918972
0.2 0.0503802006559669
};
\addlegendentry{4-layer $\bm{\sigma}$-NN}
\end{axis}

\end{tikzpicture}}~
\scalebox{0.6}{
\begin{tikzpicture}

\definecolor{color0}{rgb}{0.12156862745098,0.466666666666667,0.705882352941177}
\definecolor{color1}{rgb}{0.75,0.75,0}

\begin{axis}[
scaled ticks=false,
legend cell align={left},
legend style={fill opacity=0.8, draw opacity=1, text opacity=1, at={(0.97,0.03)}, anchor=south east, draw=white!80!black},
tick align=outside,
tick pos=left,
x grid style={white!69.0196078431373!black},
xlabel={Time (s)},
xmin=-0.01, xmax=0.21,
xtick style={color=black},
xtick={-0.05,0,0.05,0.1,0.15,0.2,0.25},
xticklabels={-0.05,0.00,0.05,0.10,0.15,0.20,0.25},
y grid style={white!69.0196078431373!black},
ylabel={Displacement (m)},
ymin=-0.005, ymax=0.06,
ytick style={color=black},
ytick={-0.01,0,0.01,0.02,0.03,0.04,0.05,0.06},
yticklabels={-0.01,0.00,0.01,0.02,0.03,0.04,0.05,0.06}
]
\addplot [very thick, color0, mark=*, mark size=3, mark options={solid,fill opacity=0}, only marks]
table {%
0 0
0.005 0.000879071318966359
0.01 0.00175178603128696
0.015 0.00261242457318852
0.02 0.00345675778867308
0.025 0.00427942738576081
0.03 0.0050753033163228
0.035 0.00583950109613968
0.04 0.0106107839172756
0.045 0.019541230444874
0.05 0.0236237531592033
0.055 0.0313289373424314
0.06 0.0350380478185725
0.065 0.0408174614436271
0.07 0.0442730136401229
0.075 0.0480651190399302
0.08 0.0509821771299189
0.085 0.053019447971144
0.09 0.0550077364727381
0.095 0.0555728014239678
0.1 0.0563137678270884
0.105 0.056260502065091
0.11 0.0561555202057213
0.115 0.0559882518807527
0.12 0.0557565179273646
0.125 0.055461132096511
0.13 0.0551037742624527
0.135 0.0546865831978064
0.14 0.0542120745133923
0.145 0.053683114752886
0.15 0.0531029032911068
0.155 0.0524749538344366
0.16 0.0518030742630135
0.165 0.0510913446992788
0.17 0.0503440938727484
0.175 0.0495658738865951
0.18 0.0487614335058666
0.185 0.0479356900987087
0.19 0.0470937003732258
0.195 0.0462406300637518
0.2 0.0453817227311881
};
\addlegendentry{Reference}
\addplot [very thick, red, dashed]
table {%
0 0
0.005 0.000878468894492047
0.01 0.00174933810622466
0.015 0.00261046040356858
0.02 0.0034563796915877
0.025 0.00427531746788229
0.03 0.00507689024278174
0.035 0.00584410310870029
0.04 0.0106983052285374
0.045 0.0197077863029154
0.05 0.0239483996070539
0.055 0.0319126303755402
0.06 0.0360661401437262
0.065 0.0424591085158404
0.07 0.0469815590641886
0.075 0.0523181594430365
0.08 0.0577139896230987
0.085 0.0633421378747367
0.09 0.0723828507293698
0.095 0.0957418892565093
0.1 0.164969494329073
0.105 0.292781023077036
0.11 0.504276231962213
0.115 0.835102253988104
0.12 1.3358057399425
};
\addlegendentry{1-layer $\Delta\bm{\sigma}$-NN}
\addplot [very thick, color1, dashed]
table {%
0 0
0.005 0.000878428664667465
0.01 0.00175099607731404
0.015 0.00261225554361125
0.02 0.00345681351252729
0.025 0.00427941009295231
0.03 0.00507400965326273
0.035 0.0058342823173674
0.04 0.0106334680068155
0.045 0.0195292899203659
0.05 0.0235694487803592
0.055 0.0312086788768744
0.06 0.0348113086077762
0.065 0.0404376046457931
0.07 0.0435038117443872
0.075 0.0464094687201033
0.08 0.0452522456998963
0.084 0.0291848375007716
0.08675 0.00588031562005591
0.08825 0.017464976127418
0.09075 0.0777884985650093
0.09375 0.17320666327903
};
\addlegendentry{2-layer $\Delta\bm{\sigma}$-NN}
\addplot [very thick, blue, dashed]
table {%
0 0
0.005 0.000879151812440363
0.01 0.00175203481036431
0.015 0.00261265214397808
0.02 0.00345726897695145
0.025 0.00427935548115949
0.03 0.00507606423467361
0.035 0.00583841295512305
0.04 0.010627497062903
0.045 0.0195180375989903
0.05 0.0235415404911276
0.055 0.0310915543736259
0.06 0.0344773685330522
0.065 0.0394465308385032
0.07 0.0399822644742155
0.0745 0.0330710185390567
0.07625 0.0254528123705376
0.077125 0.0163478817850465
0.078125 0.0208319513775599
0.078875 0.016968291521269
0.0795 0.0177727674509247
0.081875 0.115429992417118
0.084375 0.0783605642887939
0.086875 -0.00329027061807474
0.089 -0.00767543103660163
0.089625 0.0207831500332035
0.090375 0.0580230972399927
0.094875 0.103463048156662
0.097875 -0.102882118393108
0.099625 -0.129440049832539
0.1005 -0.116533771436854
0.10175 -0.015292120721858
0.10325 0.0780748005805463
0.1045 0.263349603680436
};
\addlegendentry{3-layer $\Delta\bm{\sigma}$-NN}
\addplot [very thick, green!50!black, dashed]
table {%
0 0
0.005 0.000878977001350152
0.01 0.00175184893348432
0.015 0.00261235332234743
0.02 0.00345686337078439
0.025 0.00427927974991005
0.03 0.00507526515591396
0.035 0.00583903787234587
0.03725 0.00659962573928253
0.04025 0.0115048247554527
0.04525 0.0224380469206622
0.05025 0.0801375478258914
0.05525 0.378117928782923
};
\addlegendentry{4-layer $\Delta\bm{\sigma}$-NN}
\end{axis}

\end{tikzpicture}}
  \caption{Displacement trajectories of the right end point for SPD-NN~(left), $\bm{\sigma}$-NN~(middle), and $\Delta\bm{\sigma}$-NN~(right) in the NN-FEM test. }
  \label{FIG:TRUSS-NN-FEM-TEST-DISP}
\end{figure}
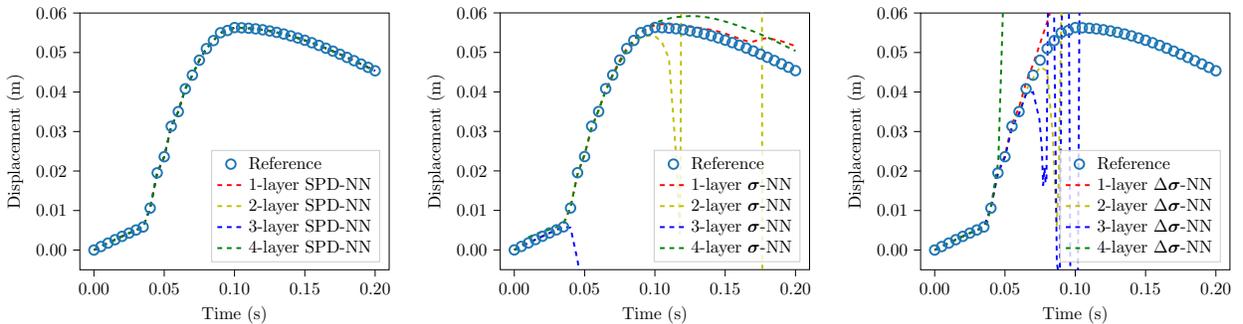

\subsubsection{Comparison of Different Neural Network Architectures for SPD-NN}\label{sect:cnn}

In this experiment, we consider different neural network architectures. We vary widths, depths and activation functions of the neural network used in SPD-NN, while keeping other settings the same as \cref{sect:compare_cholnn}. The errors in the NN-FEM test in terms of mean squared errors of the predicted displacements at the final time $T$
are shown in \cref{tab:nncompare}. We see that the \texttt{tanh} activation function is in general more accurate than others if appropriate widths and depths of the neural network are chosen~(For this case, the training set consists of 800 data, the number of neural network parameters is not supposed to outnumbers it too much). However, ReLU and leaky ReLU are more robust for deep and wide neural networks, despite being less accurate. 

\begin{table}[]
\centering
\begin{tabular}{@{}l|l|llll@{}}
\toprule
 & Depth\textbackslash{}Width & 2 & 10 & 20 & 40 \\ \midrule
\multirow{4}{*}{\texttt{tanh}} & 1 & 1.1$\times10^{-5}$ & 4.2$\times10^{-5}$ & 2.0$\times10^{-5}$ & 3.0$\times10^{-5}$ \\
 & 3 & 4.7$\times10^{-5}$ & 1.8$\times10^{-6}$ & 1.3$\times10^{-6}$ & NaN \\
 & 8 & 5.5$\times10^{-4}$ & 1.3$\times10^{-5}$ & NaN & NaN \\
 & 20 & 8.8$\times10^{-6}$ & 5.9$\times10^{-5}$ & NaN & NaN \\ \midrule
\multirow{4}{*}{ReLU} & 1 & 9.3$\times10^{-3}$ & 3.7$\times10^{-5}$ & 6.1$\times10^{-5}$ & 4.2$\times10^{-3}$ \\
 & 3 & 9.3$\times10^{-3}$ & 3.1$\times10^{-3}$ & 4.7$\times10^{-3}$ & 5.3$\times10^{-3}$ \\
 & 8 & 9.1$\times10^{-3}$ & 1.8$\times10^{-4}$ & 6.1$\times10^{-3}$ & 1.0$\times10^{-2}$ \\
 & 20 & 3.5 & 5.9$\times10^{-6}$ & 3.3$\times10^{-3}$ & NaN \\ \midrule
\multirow{4}{*}{leaky ReLU} & 1 & 4.6$\times10^{-3}$ & 4.8$\times10^{-3}$ & 4.3$\times10^{-3}$ & 1.2$\times10^{-2}$ \\
 & 3 & 8.9$\times10^{-3}$ & 1.2$\times10^{-2}$ & 1.0$\times10^{-2}$ & 2.1$\times10^{-3}$ \\
 & 8 & 5.3$\times10^{-3}$ & 5.0$\times10^{-3}$ & 6.6$\times10^{-3}$ & 5.4$\times10^{-3}$ \\
 & 20 & 9.3$\times10^{-3}$ & 4.4$\times10^{-3}$ & 8.6$\times10^{-3}$ & 1.8$\times10^{-2}$ \\ \midrule
\multirow{4}{*}{SELU} & 1 & 1.0$\times10^{-2}$ & 6.6$\times10^{-4}$ & 2.2$\times10^{-4}$ & 2.2$\times10^{-4}$ \\
 & 3 & 9.6$\times10^{-3}$ & 9.8$\times10^{-5}$ & 4.3$\times10^{-4}$ & 3.8$\times10^{-4}$ \\
 & 8 & 9.2$\times10^{-5}$ & 7.6$\times10^{-3}$ & 4.8$\times10^{-5}$ & NaN \\
 & 20 & 8.3$\times10^{-6}$ & 5.8$\times10^{-5}$ & NaN & NaN \\ \midrule
\multirow{4}{*}{ELU} & 1 & 4.2$\times10^{-3}$ & 4.6$\times10^{-3}$ & 5.0$\times10^{-3}$ & 6.6$\times10^{-3}$ \\
 & 3 & 6.8$\times10^{-5}$ & 5.7$\times10^{-3}$ & 4.7$\times10^{-4}$ & 4.1$\times10^{-3}$ \\
 & 8 & 8.0$\times10^{-6}$ & 8.4$\times10^{-6}$ & NaN & NaN \\
 & 20 & 6.6$\times10^{-5}$ & NaN & NaN & NaN \\ \bottomrule
\end{tabular}
\caption{Mean squared errors of the predicted displacements at the final time $T$ for different neural network architectures, including \texttt{tanh}, ReLU~(rectified linear unit), ELU~(exponential linear unit), SELU~(scaled exponential linear unit), and leaky ReLU~(leaky rectified linear unit) with leakage 0.1. The width is the number of activation nodes in each hidden layer. The depth is the number of hidden layers. ``NaN'' denotes the numerical simulation fails (e.g., due to numerical instability) using the trained NN-based constitutive relations.}
\label{tab:nncompare}
\end{table}

\subsubsection{Choice of $\tilde\sigma_Y$ and $d$}

We consider the impact of $\tilde\sigma_Y$ and $d$ in the transition function $D$~(see~\cref{EQ:TRANSIT_TRUSS})
\begin{equation*}
D(\sigma^{n}, \tilde{\sigma}_Y) = \texttt{sigmoid}\left(\frac{(\sigma^{n})^2 - \tilde{\sigma}_Y^2} {d\tilde{\sigma}_Y^2}\right)
\end{equation*}
on the accuracy of the SPD-NN. As mentioned before, $\tilde\sigma_Y$ controls where the transition from the elastic form~(\cref{EQ:PLASTICIT_LINEAR}) to the plastic form~(\cref{EQ:PLASTICIT_NONLINEAR}) happens  and $d$ controls the sharpness of the transition~(see \cref{fig:H}). In this numerical experiment, we use a fully connected neural network with 3 hidden layers, 20 neurons in each layer, and the \texttt{tanh} activation function. We vary $\tilde\sigma_Y$ and $d$ in \cref{EQ:TRANSIT_TRUSS}. 

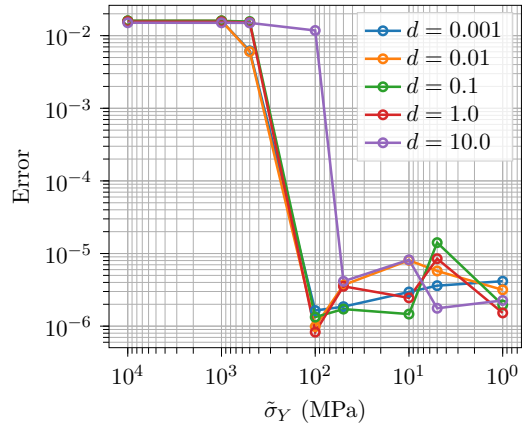
\begin{figure}[hbtp]
\centering
	\scalebox{0.8}{
\begin{tikzpicture}

\definecolor{color0}{rgb}{0.12156862745098,0.466666666666667,0.705882352941177}
\definecolor{color1}{rgb}{1,0.498039215686275,0.0549019607843137}
\definecolor{color2}{rgb}{0.172549019607843,0.627450980392157,0.172549019607843}
\definecolor{color3}{rgb}{0.83921568627451,0.152941176470588,0.156862745098039}
\definecolor{color4}{rgb}{0.580392156862745,0.403921568627451,0.741176470588235}

\begin{axis}[
legend cell align={left},
legend style={fill opacity=0.8, draw opacity=1, text opacity=1, draw=white!80.0!black},
log basis x={10},
log basis y={10},
tick align=outside,
tick pos=left,
x dir=reverse,
x grid style={white!69.01960784313725!black},
xlabel={$\tilde\sigma_Y$ (MPa)},
minor tick num = 4,
xmajorgrids,
xminorgrids,
yminorgrids,
xmin=0.630957344480193, xmax=15848.9319246111,
xmode=log,
xtick style={color=black},
y grid style={white!69.01960784313725!black},
ylabel={Error},
ymajorgrids,
ymin=5.05339185513665e-07, ymax=0.0261649678164883,
ymode=log,
ytick style={color=black},
]
\addplot [very thick, color0, mark=o, mark size=2, mark options={solid}]
table {%
1 4.15699948486241e-06
5 3.6083790255892e-06
10 2.95371683045814e-06
50 1.85466123641828e-06
100 1.64164448607221e-06
500 0.00618071448107817
1000 0.0159750137242565
10000 0.0159750137242565
};
\addlegendentry{$d=0.001$}
\addplot [very thick, color1, mark=o, mark size=2, mark options={solid}]
table {%
1 3.16095040698029e-06
5 5.74895139696578e-06
10 8.0750416555111e-06
50 3.67602155194961e-06
100 9.86340766302514e-07
500 0.00603245778170436
1000 0.0159750137242565
10000 0.0159750137242565
};
\addlegendentry{$d=0.01$}
\addplot [very thick, color2, mark=o, mark size=2, mark options={solid}]
table {%
1 1.95494956142115e-06
5 1.41241394095072e-05
10 1.4655402256274e-06
50 1.71758906325267e-06
100 1.33229021223259e-06
500 0.0156536013255375
1000 0.0159741187087351
10000 0.0159745038400026
};
\addlegendentry{$d=0.1$}
\addplot [very thick, color3, mark=o, mark size=2, mark options={solid}]
table {%
1 1.52134085596964e-06
5 8.48628556283611e-06
10 2.45391263718042e-06
50 3.53508232407238e-06
100 8.27679008832329e-07
500 0.0152383385679192
1000 0.0155292344485064
10000 0.0156009686724288
};
\addlegendentry{$d=1.0$}
\addplot [very thick, color4, mark=o, mark size=2, mark options={solid}]
table {%
1 2.25901068526364e-06
5 1.76279747704155e-06
10 8.19632320233756e-06
50 4.14237721836772e-06
100 0.0117819798115206
500 0.0149948626906782
1000 0.0150374417321327
10000 0.015060066646183
};
\addlegendentry{$d=10.0$}
\end{axis}

\end{tikzpicture}}
  \caption{Impacts of $\tilde\sigma_Y$ and $d$ on the approximation accuracy of SPD-NNs.}
  \label{fig:compare_d}
\end{figure}

The error plot is shown in \cref{fig:compare_d}, where the error metric is the same as \cref{sect:cnn}. The reported errors are the average errors of 10 simulations with different initial guesses. We can see that the accuracy of the SPD-NN is less sensitive to $d$. Additionally, as long as we choose a small enough $\tilde\sigma_Y$ so that the corresponding plastic form~(\cref{EQ:PLASTICIT_NONLINEAR}) covers the plasticity regime,  the SPD-NN is sufficiently expressive to approximate the constitutive relation. This justifies our choices $d=0.1$ and $\tilde\sigma_Y = 0.1~\textrm{GPa}$.

\subsection{2D Thin Plate with Different Materials}
\label{SEC:PLATE}
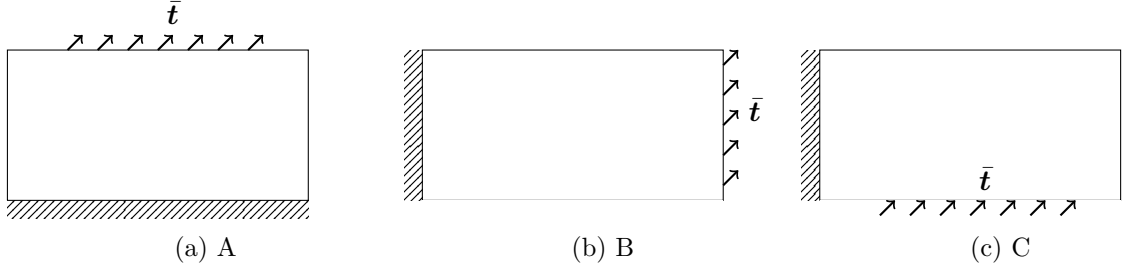
\begin{figure}[htbp]
\centering
\begin{subfigure}[b]{0.32\textwidth}
\begin{tikzpicture}[scale=0.4]
\draw (0,0) -- (10,0) -- (10,5) -- (0,5) -- (0,0);
\fill [pattern = north east lines] (0.0,-0.6) rectangle (10.0, 0.0);
\draw[thick,->] (2, 5) -- (2.5, 5.5);
\draw[thick,->] (3, 5) -- (3.5, 5.5);
\draw[thick,->] (4, 5) -- (4.5, 5.5);
\draw[thick,->] (5, 5) -- (5.5, 5.5) node[anchor=south]{$\bar{\bm{t}}$};
\draw[thick,->] (6, 5) -- (6.5, 5.5);
\draw[thick,->] (7, 5) -- (7.5, 5.5);
\draw[thick,->] (8, 5) -- (8.5, 5.5);
\end{tikzpicture}
\caption{A}
\end{subfigure}%
\begin{subfigure}[b]{0.32\textwidth}
\begin{tikzpicture}[scale=0.4]
\draw (0,0) -- (10,0) -- (10,5) -- (0,5) -- (0,0);
\fill [pattern = north east lines] (-0.6,0.0) rectangle (0.0, 5.0);
\fill [white] (0.0,-0.6) rectangle (10.0, 0.0);

\draw[thick,->] (10, 0.5) -- (10.5, 1);
\draw[thick,->] (10, 1.5) -- (10.5, 2);
\draw[thick,->] (10, 2.5) -- (10.5, 3) node[anchor=west]{$\bar{\bm{t}}$};
\draw[thick,->] (10, 3.5) -- (10.5, 4);
\draw[thick,->] (10, 4.5) -- (10.5, 5);

\end{tikzpicture}
\caption{B}
\end{subfigure}%
\begin{subfigure}[b]{0.32\textwidth}
\begin{tikzpicture}[scale=0.4]
\draw (0,0) -- (10,0) -- (10,5) -- (0,5) -- (0,0);
\fill [pattern = north east lines] (-0.6,0.0) rectangle (0.0, 5.0);
\fill [white] (0.0,-0.6) rectangle (10.0, 0.0);
\draw[thick,->] (2, -0.5) -- (2.5, 0);
\draw[thick,->] (3, -0.5) -- (3.5, 0);
\draw[thick,->] (4, -0.5) -- (4.5, 0);
\draw[thick,->] (5, -0.5) -- (5.5, 0) node[anchor=south]{$\bar{\bm{t}}$};
\draw[thick,->] (6, -0.5) -- (6.5, 0);
\draw[thick,->] (7, -0.5) -- (7.5, 0);
\draw[thick,->] (8, -0.5) -- (8.5, 0);
\end{tikzpicture}
\caption{C}
\end{subfigure}%

\caption{Schematic of the boundary conditions of the thin plate tests.}
\label{fig:plate_bc}
\end{figure}

In this section, we consider thin plate coupons of size $L_x = 10~\textrm{cm}$ by $L_y =5~\textrm{cm}$ with the plane stress assumption~($L_z = 0.1~\textrm{cm}$). These plates are made of different materials, including hyperelastic material~(finite deformation), elasto-plastic material~(infinitesimal deformation), and fiber-reinforced multiscale elasto-plastic material~(infinitesimal deformation).  
For each case, the plate is tested under $13$ loading conditions as depicted in \cref{fig:plate_bc}. The prescribed time-dependent load force $\bar{\bm{t}}\in \mathcal{R}^{2}$ consists of both loading and unloading parts and takes the form
$$\bar{\bm{t}} = p \sin\left(\frac{t\pi}{T}\right)$$
here $p\in \mathcal{R}^{2}$ is the loading parameter vector, as follows,  
\begin{itemize}
\item A) clamp on the bottom edge and impose force load on the top edge. \\
A1:~$(0, p_1)$, A2:~$(0,-p_1)$, A3:~$(p_3, 0)$, A4:~$(-p_3, 0)$, A5:~$(p_3/\sqrt{2},p_1/\sqrt{2})$, and A6:~$(0.75p_3, 0)$.  
\item B) clamp on the left edge and impose force load on the right edge. \\
B1:~$(p_1, 0)$, B2:~$(-p_1, 0)$, B3:~$(0, p_2)$, B4:~$(0, -p_2)$, B5:~$(p_1/\sqrt{2},p_2/\sqrt{2})$, and B6:~$(0, 0.75p_2)$.
\item C) clamp on the left edge and impose force load on the bottom edge. \\
C1:~$\left(0, \frac{p_2L_x}{\sqrt{2\pi}\sigma_X}\exp{\big(\frac{-(x-x_0)^2}{\sigma_X^2}\big)}\right)$, with $x_0 = \frac{5L_x}{6}\, \textrm{and}\, \sigma_X = 0.2L_x$.
\end{itemize}
The total simulation time is $T = 0.2s$.
We use A1-A5 and B1-B5 as training data and A6, B6, and C1 as test data. 
Both training procedures discussed in \cref{SEC:DIRECT_DATA,SEC:INDIRECT_DATA} are applied. 
For the direct input-output data training, the strain-stress sequential data are extracted from all Gaussian points in the  training sets.
For the indirect data training, the full-filed displacement fields on the 21 by 11 grid~($0.5$~cm interval) from the training data are extracted. 
Therefore, this approach is potentially applicable to experimental data. The pre-training is required to obtain good initial guesses.

\subsubsection{Hyperelasticity}
\label{SEC:PLATE_HYPERELASTICITY}
The plate is made of the incompressible Rivlin-Saunders material~\cite{pascon2019large,rivlin1956rheology}  with the density $\rho=800~\textrm{kg/m}^3$ and the energy density function
\begin{equation*}
w = c_1(I_1 - 3) + c_2(I_2 - 3)
\end{equation*}
Here $c_1 = 0.1863~\textrm{MPa}$ and $c_2 = 0.00979~\textrm{MPa}$, and $I_1$, $I_2$, $I_3$ are three scalar invariants of the right Cauchy-Green stretch tensor $\bm{C} = \bm{F}\bm{F}^T = 2\bm\epsilon + 1$, where
\begin{equation*}
I_1 = \mathrm{tr}\bm{C} \quad I_2 = \frac{1}{2}[(\mathrm{tr} \bm{C})^2 - \mathrm{tr} \bm{C}^2] \quad I_3 = J^2 = \det \bm{C}
\end{equation*}
The incompressibility implies that $J = 1$. 
The plate is assumed to undergo finite deformations~(see \cref{SEC:FINITE_DEFORMATION}), and the second Piola-Kirchhoff stress tensor reads
\begin{equation}
\label{EQ:HYPERELASTICITY}
\bm{S} = \frac{\partial w}{\partial \bm{\epsilon}} + \lambda_J\frac{\partial J}{\partial \bm{\epsilon}}
\end{equation}
here $\lambda_J$ is the Lagrangian multiplier, which can be calculated based on the plane stress assumption $\bm{S}_{33} = 0$.
The plate domain is discretized by $20 \times 10$ quadratic quadrilateral elements. The time step size is $\Delta t = 0.001~\mathrm{s}$.
The data sets are generated with load parameters $(p_1, p_2, p_3) = (44800, 4480, 16800)~\textrm{N/m}$.

Both the direct data training approach~\cref{SEC:DIRECT_DATA} and the indirect data training approach~\cref{SEC:INDIRECT_DATA} are applied to train a SPD-NN:
\begin{equation}
\bm S^{n+1} = \mathsf{L}_{\bt}(\bm\epsilon^{n+1}) \mathsf{L}_{\bt}(\bm\epsilon^{n+1})^T(\bm\epsilon^{n+1} - \bm\epsilon^{n}) + \bm S^{n}
\end{equation}
where the neural network consists of 4 hidden layers and 20 neurons in each layer.

\begin{figure}[htpb]
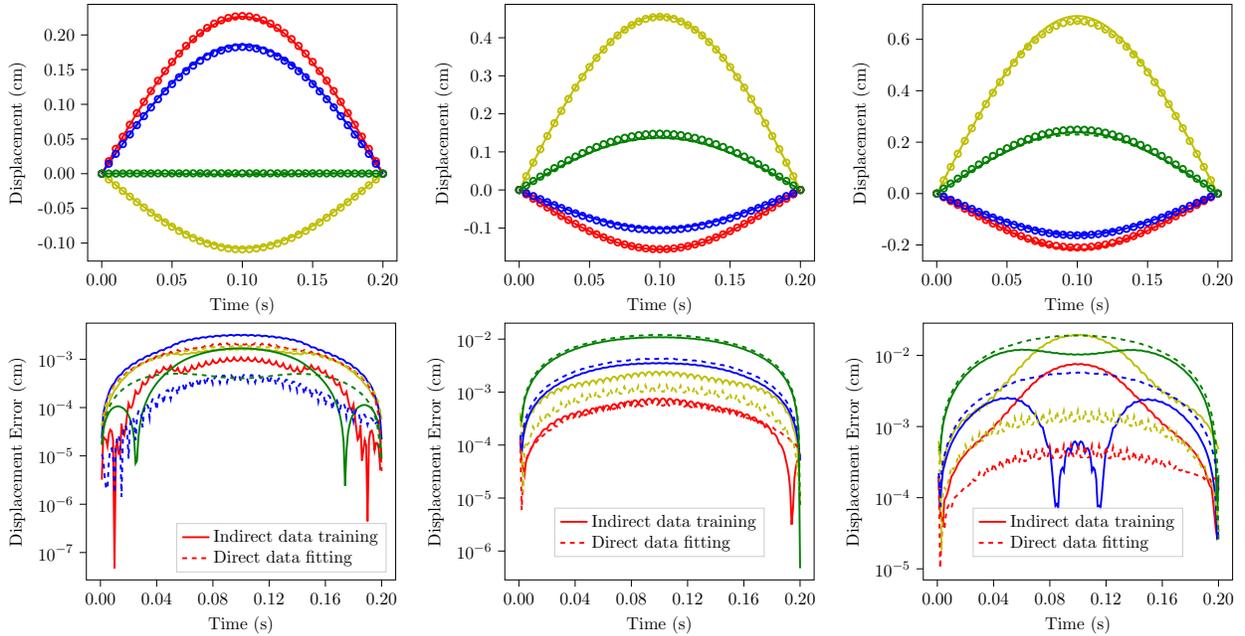

\centering
\scalebox{0.6}{\input{figures/plate_hyperelasticity_disp_nn1_stiffmat_from3_test106.tex}}~
\scalebox{0.6}{\input{figures/plate_hyperelasticity_disp_nn1_stiffmat_from3_test206.tex}}~
\scalebox{0.6}{\input{figures/plate_hyperelasticity_disp_nn1_stiffmat_from3_test300.tex}}
\scalebox{0.6}{\input{figures/diff106.tex}}~
\scalebox{0.6}{\input{figures/diff206.tex}}~
\scalebox{0.6}{\input{figures/diff300.tex}}
  \caption{Top: Trajectories of displacement at top-right~(red: $u_x$, yellow: $u_y$) and top-middle points~(blue: $u_x$, green: $u_y$) of the 2D hyperelastic plate for test A6~(left), B6~(middle) and C1~(right), defined on page~\pageref{SEC:PLATE}. The reference solutions are marked by empty circles, the solutions obtained by the SPD-NN trained using indirect data are marked by solid lines, and the solutions obtained by SPD-NN trained with direct data are marked by dashed lines. Bottom: The absolute errors of displacements for each cases.}
  \label{FIG:HYPERELASTIC_DISP}
\end{figure}

\begin{figure}[htpb]
\centering
 \includegraphics[width=0.33\textwidth]{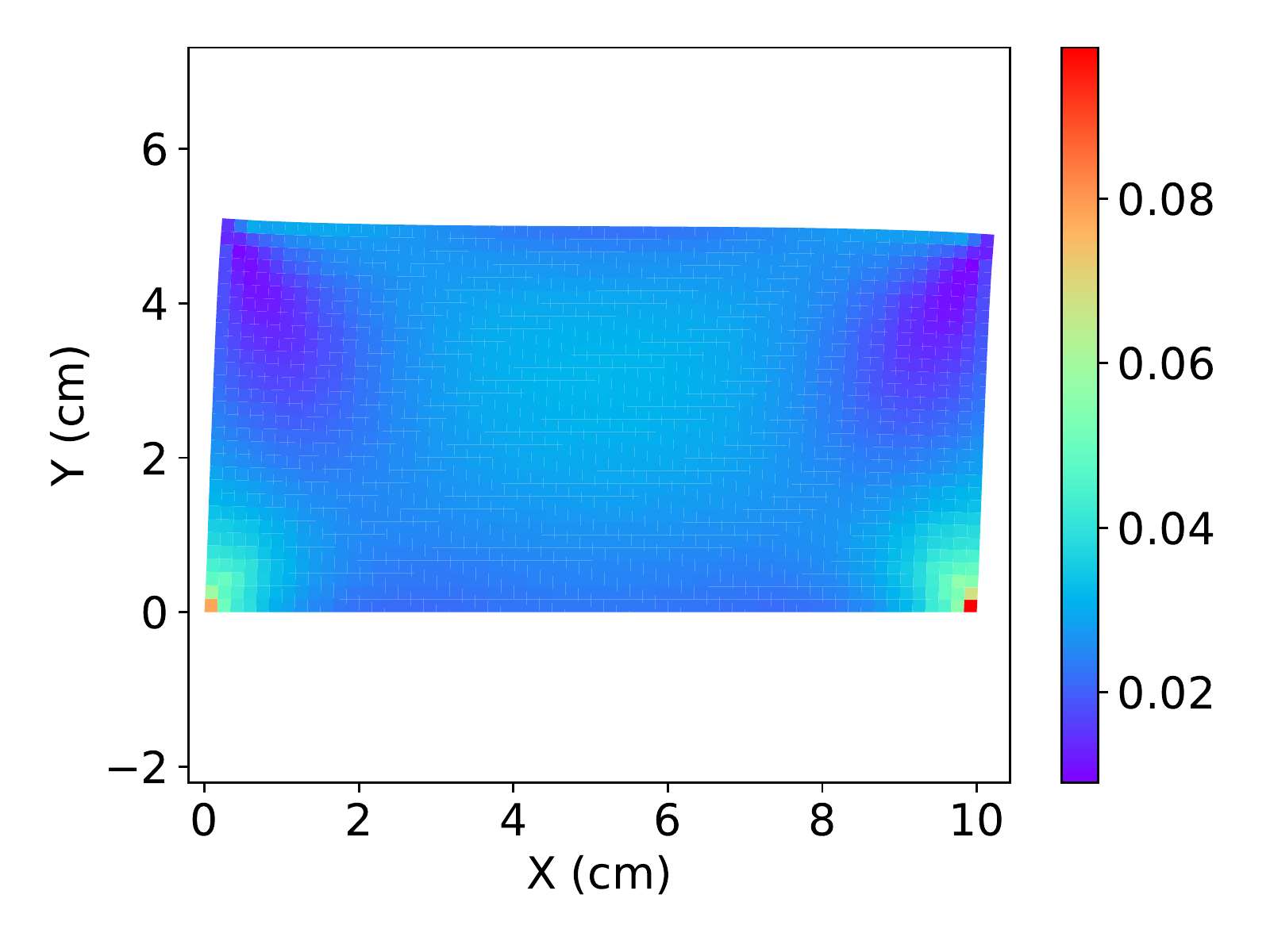}~
 \includegraphics[width=0.33\textwidth]{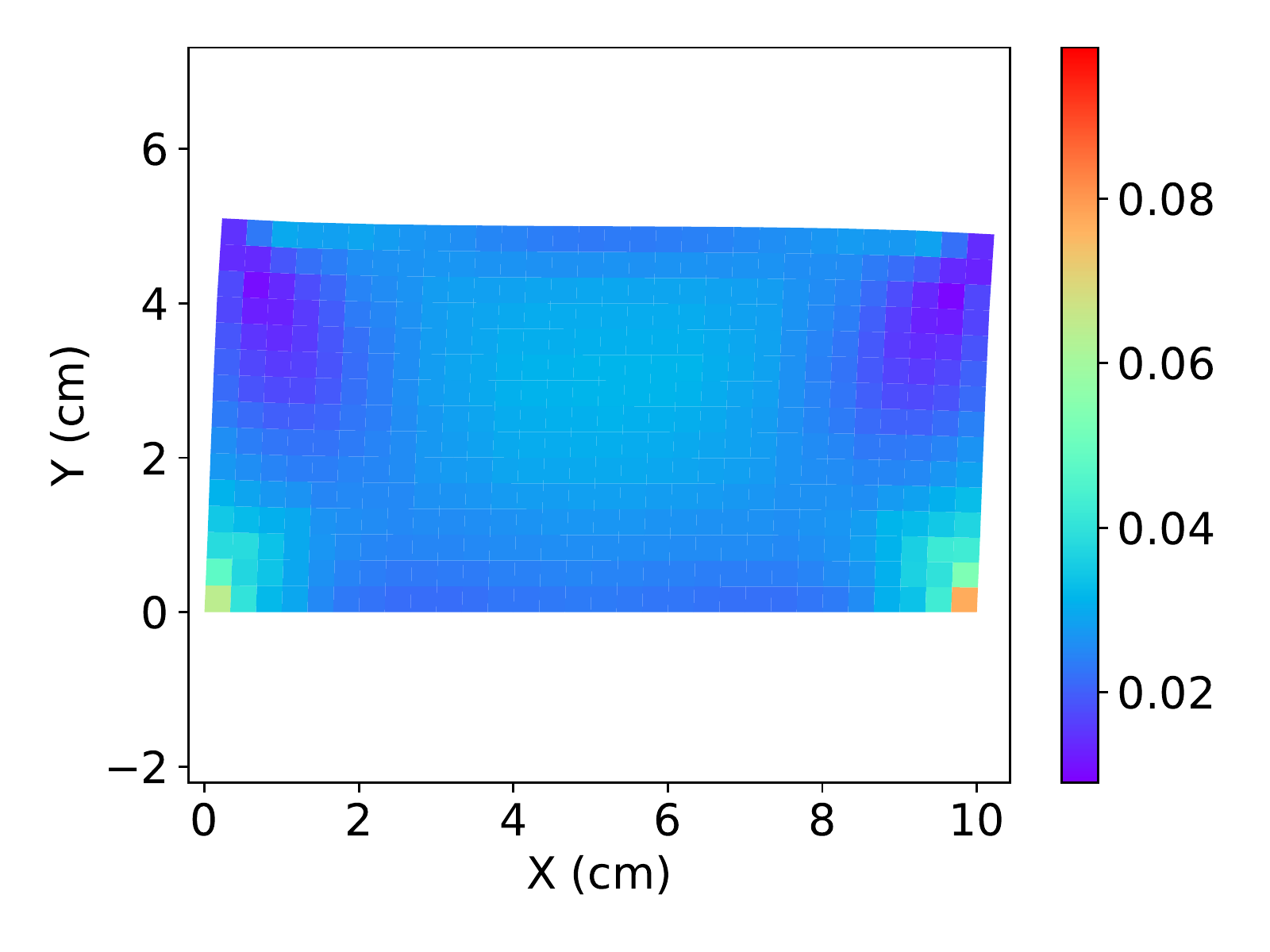}~\includegraphics[width=0.33\textwidth]{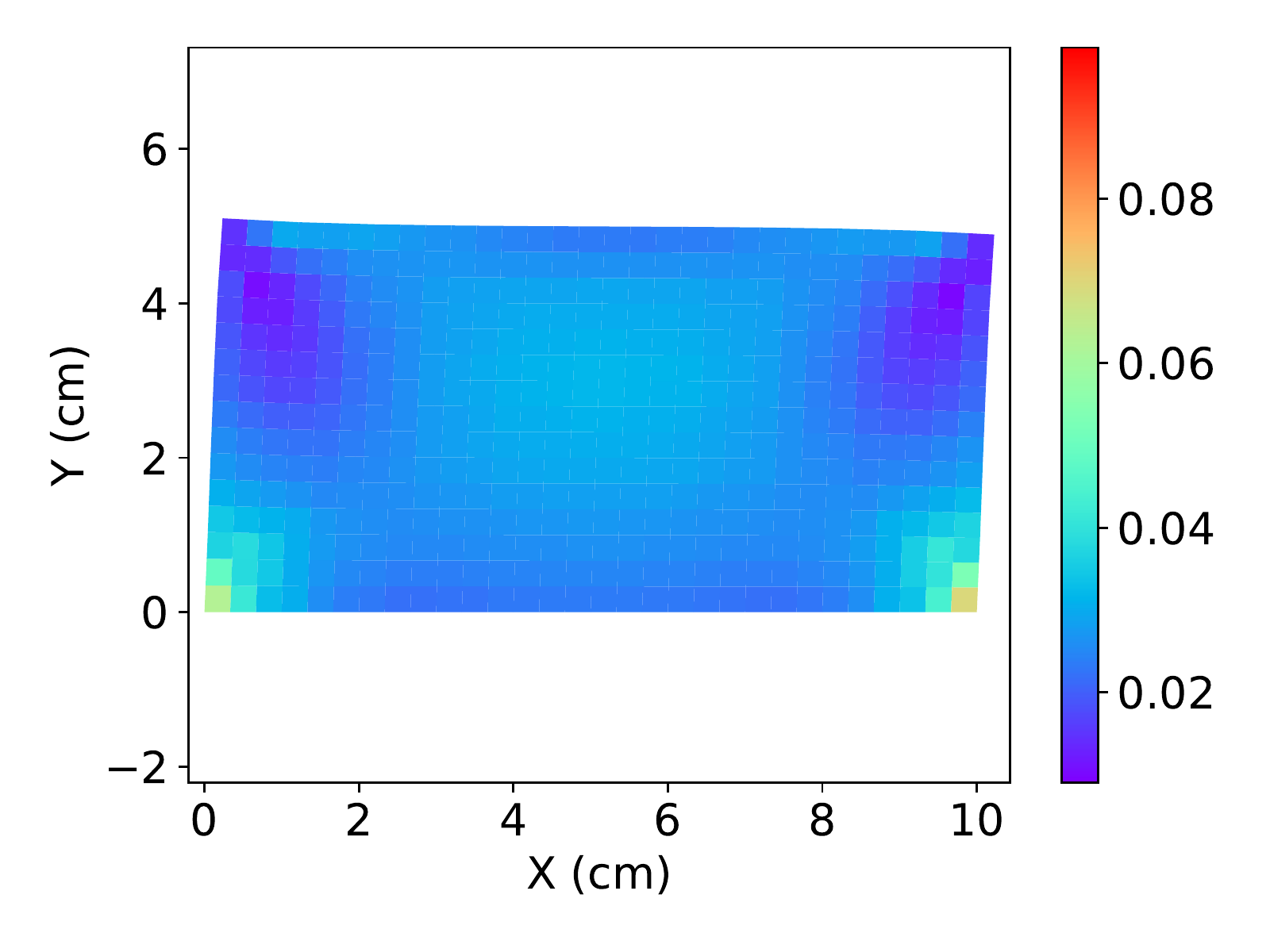}
  
  \includegraphics[width=0.33\textwidth]{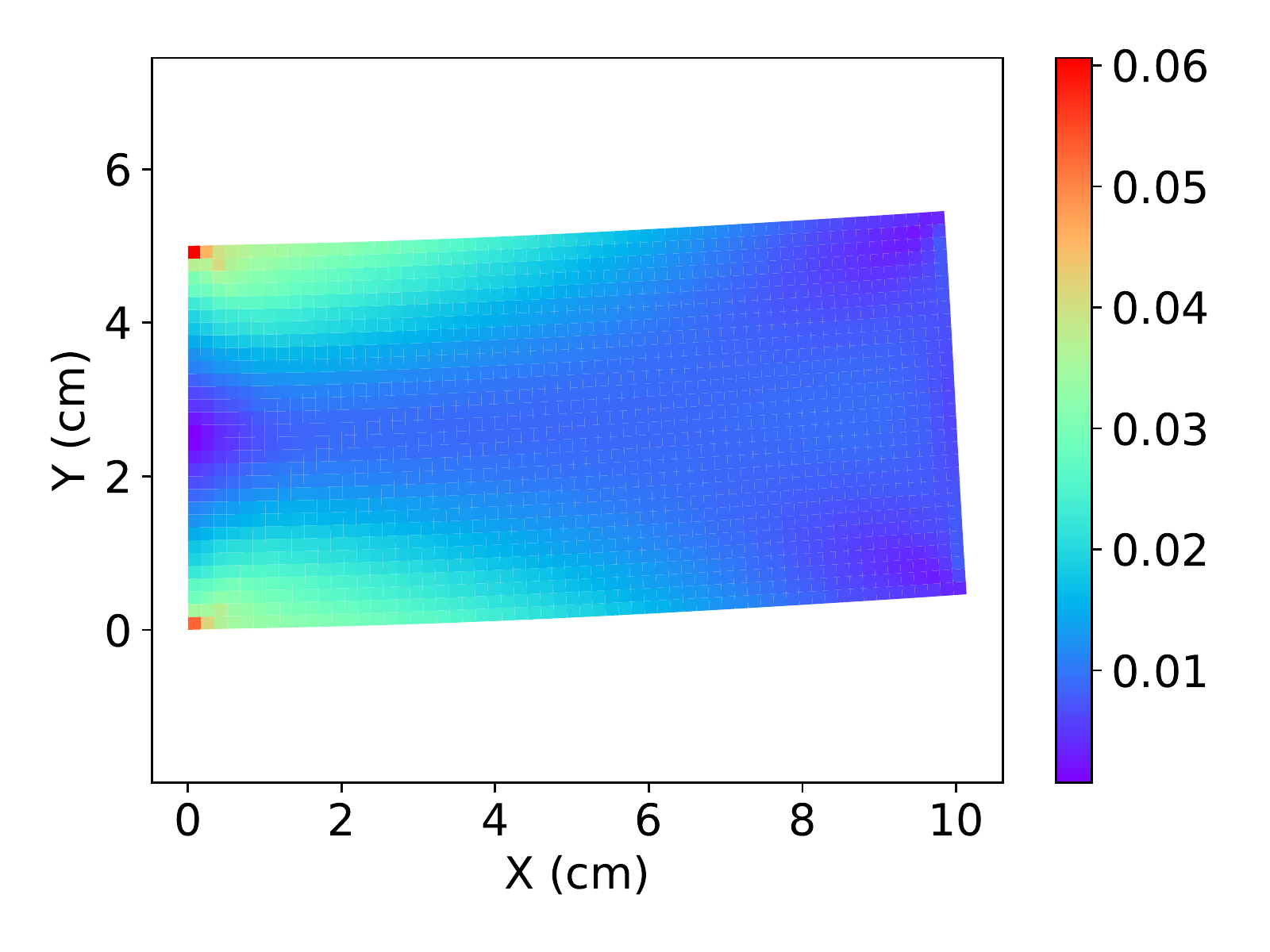}~
  \includegraphics[width=0.33\textwidth]{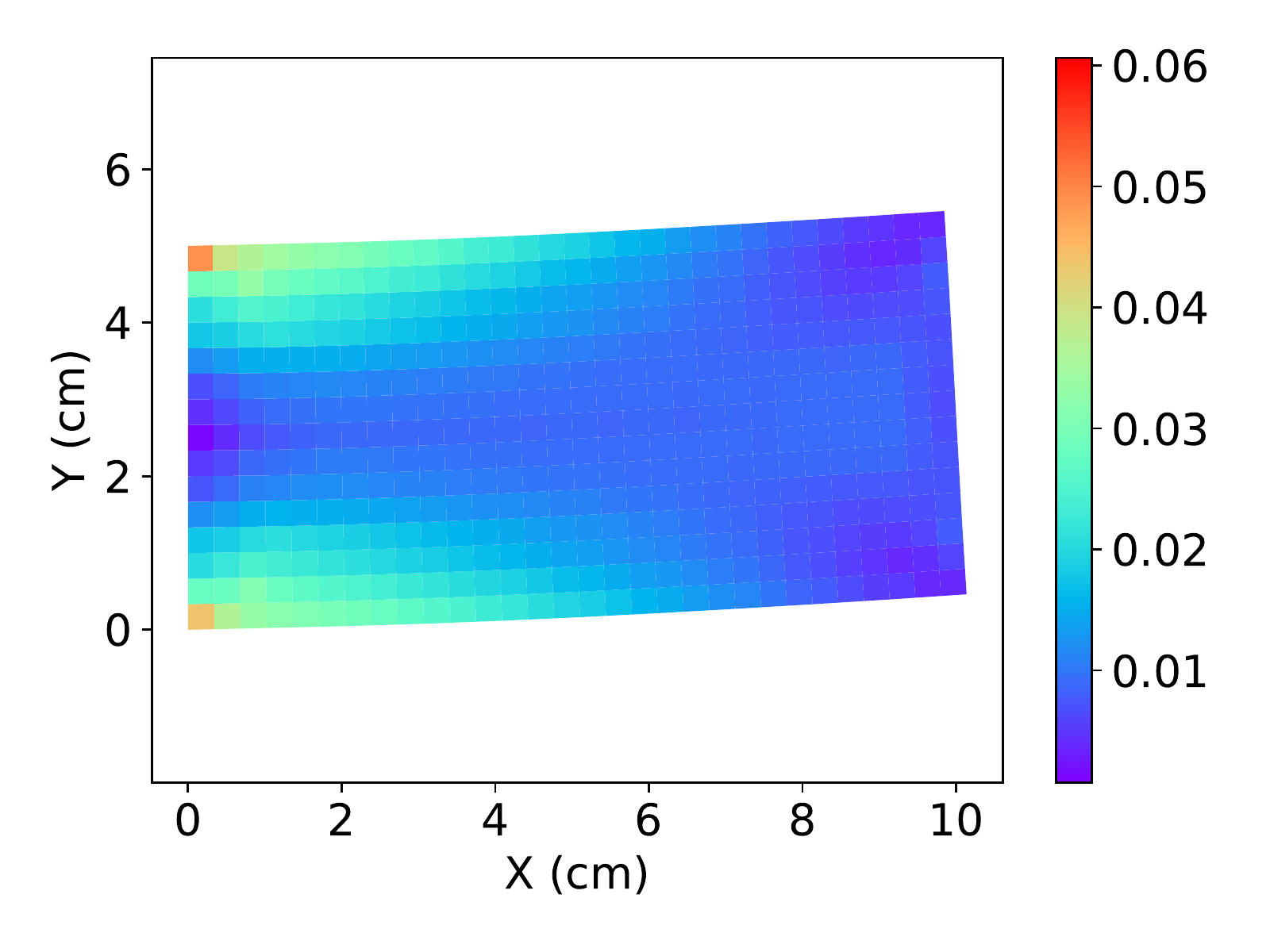}~\includegraphics[width=0.33\textwidth]{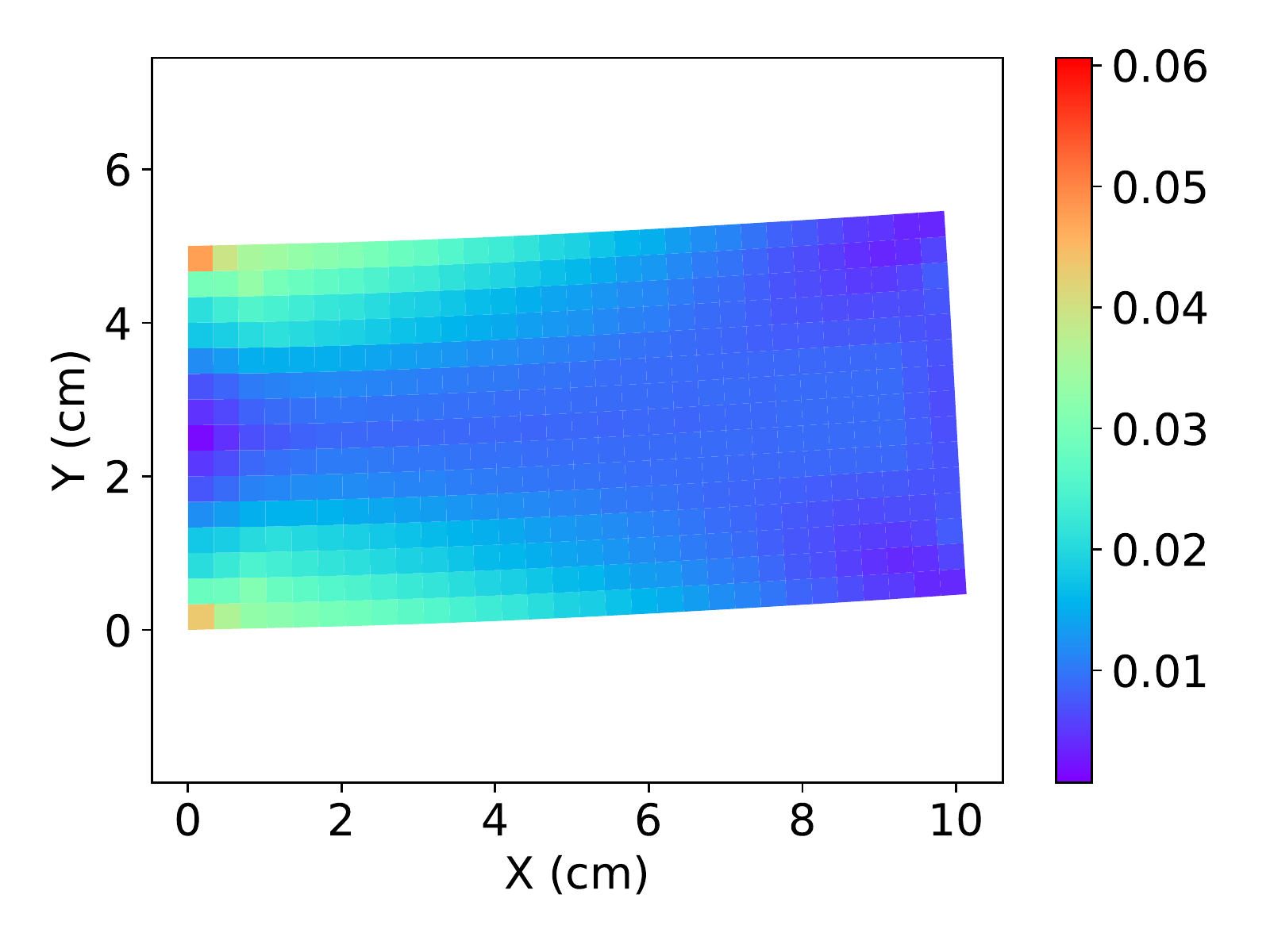}

\includegraphics[width=0.33\textwidth]{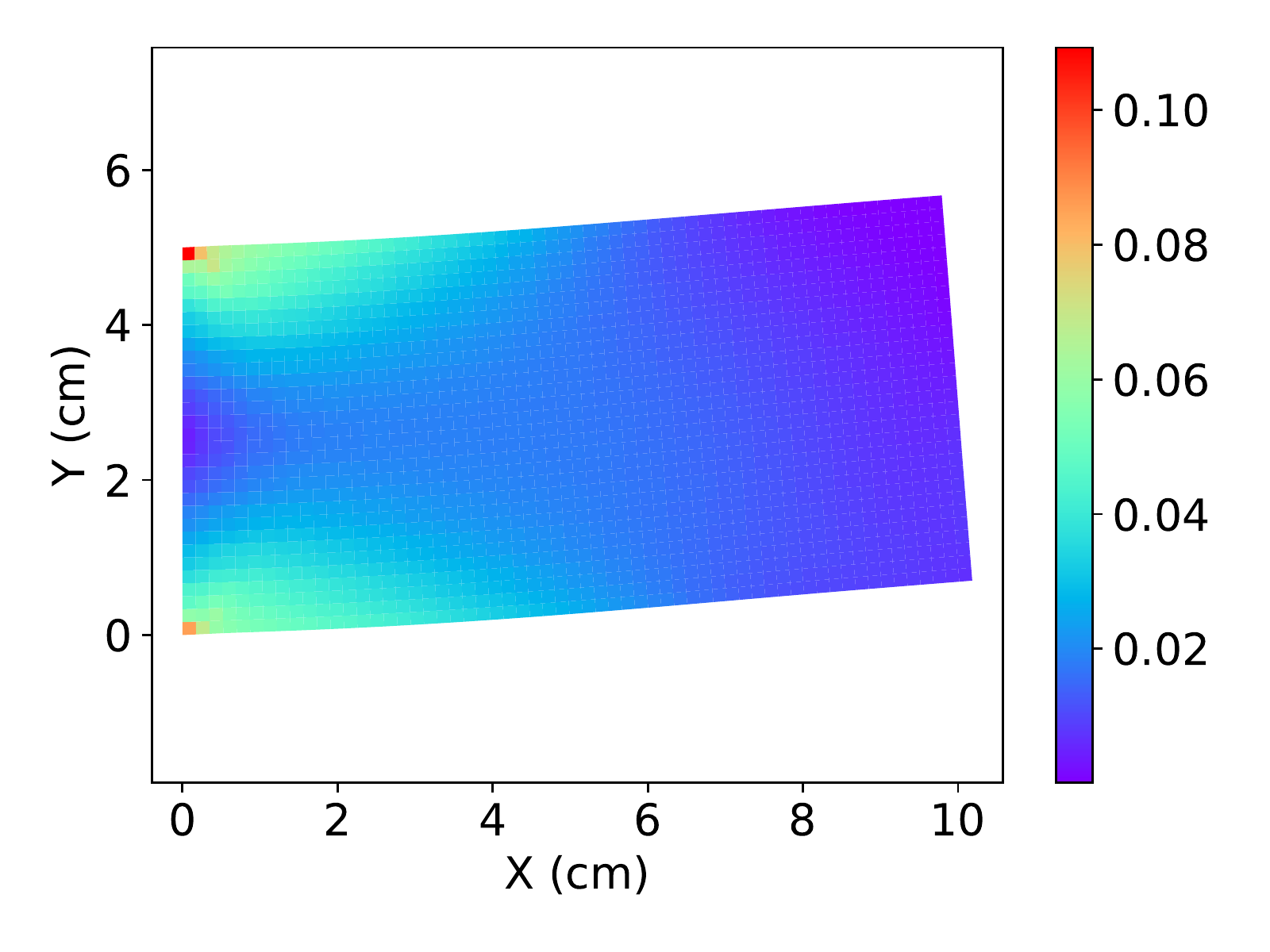}~
  \includegraphics[width=0.33\textwidth]{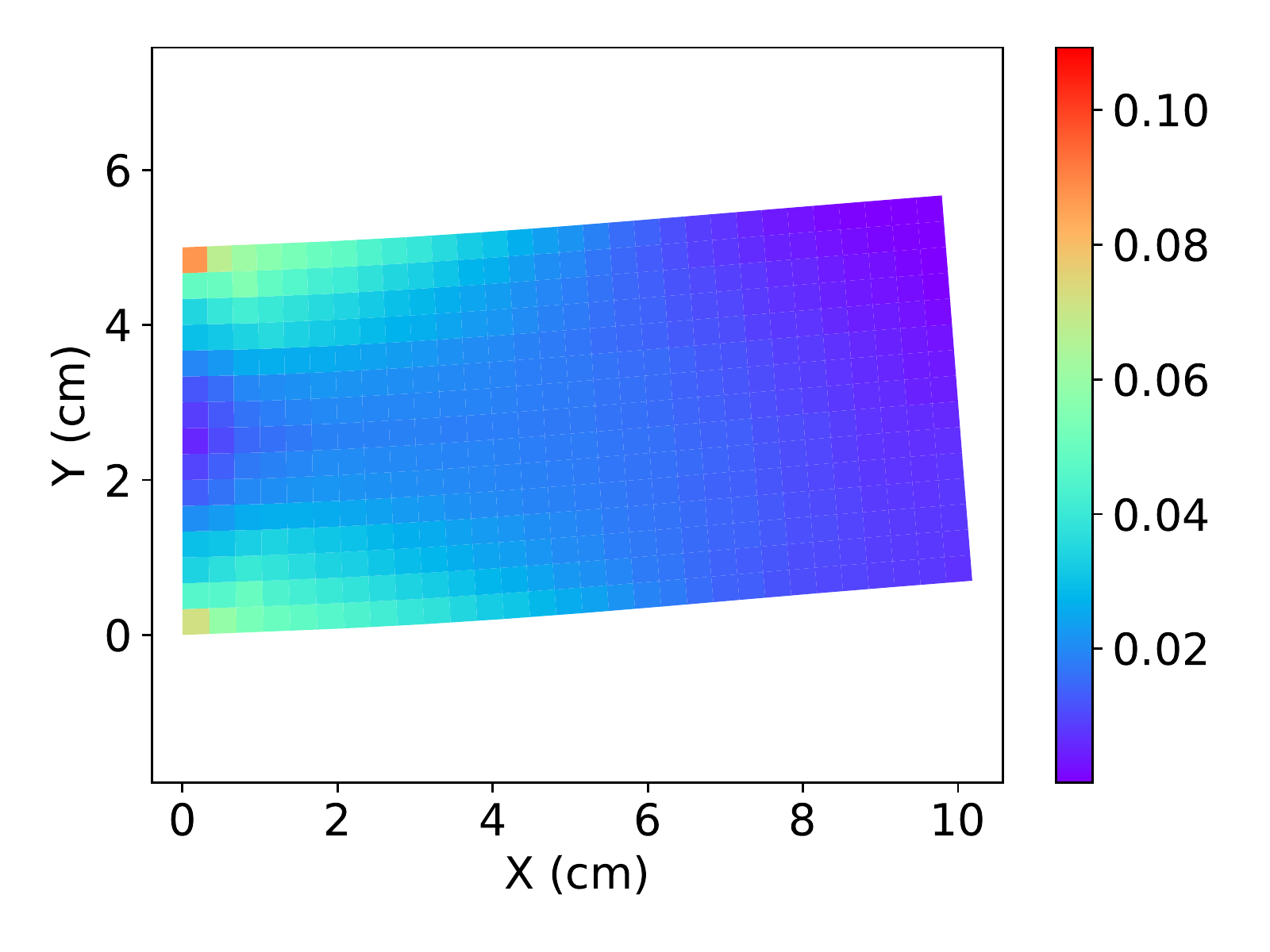}~\includegraphics[width=0.33\textwidth]{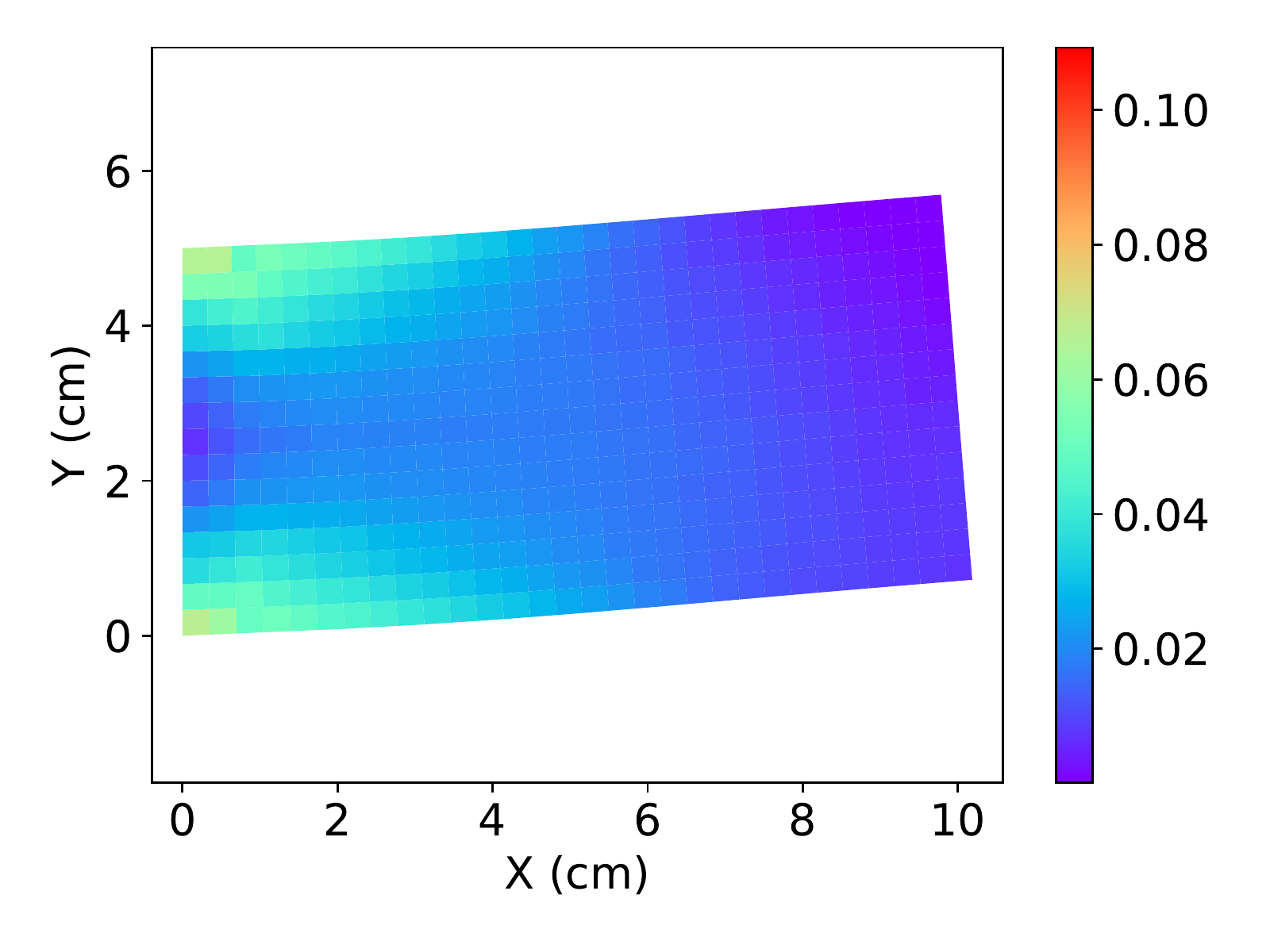}

  \caption{The von Mises stress~(MPa) fields at $t=\frac{T}{2}$ for the 2D hyperelastic plate for test A6~(top), B6~(middle) and C1~(bottom), defined on page~\pageref{SEC:PLATE}.  
  From left to right: reference solutions~(on a fine mesh), solutions obtained by SPD-NN trained with direct data, solutions obtained by SPD-NN trained with indirect data.}
    \label{FIG:HYPERELASTIC_STRESS}
\end{figure}

The predicted trajectories of displacements at the top-right  and top-middle points as a function of time and the references for all test cases are depicted in \cref{FIG:HYPERELASTIC_DISP}.  All predicted results are in good agreement with the references, and the direct input-output data training approach leads to slightly better results for case C1.
The predicted von Mises stress fields at $t=\frac{T}{2}$ for all test cases and the references are depicted in \cref{FIG:HYPERELASTIC_STRESS}.

\subsubsection{Elasto-Plasticity}
\label{SEC:PLASTICITY-PLATE}
The plate is made of titanium, which is assumed to be elasto-plastic material with density $\rho = 4200~\textrm{kg/m}^3$. The constitutive relation is 
$$\bm{\sigma} = \mathsf{C} \bm{\epsilon}$$
here $ \mathsf{C}$ denotes the isotropic plane stress stiffness tensor with Young's modulus $E = 100~\textrm{GPa} $ and Poisson's ratio $\nu = 0.35$. The von Mises yield function with isotropic hardening has the form
\begin{equation}
f = \sqrt{\sigma_{11}^2 - \sigma_{11}\sigma_{22} + \sigma_{22}^2 + 3\sigma_{12}^2}  - \sigma_Y - K\alpha
\end{equation}
The yield strength  $\sigma_Y = 0.97~\textrm{GPa}$ and the plastic modulus $K = 10~\textrm{GPa}$, the internal hardening variable $\alpha$ follows the simplest evolutionary equation
\begin{equation}
\dot{\alpha} = \dot{\lambda}
\end{equation}
This plate domain is discretized by $20 \times 10$ quadratic quadrilateral elements. And the time step size is $\Delta t = 0.001s$.

As for the SPD-NN~\cref{EQ:CHOL_NN_2}, the estimated yield strength is $\tilde{\sigma}_Y = 0.32\textrm{GPa}$, the transition function is 
\begin{equation*}
D(\sigma_{\mathrm{vm}}^{n}, \tilde{\sigma}_Y) = \texttt{sigmoid}\left(\frac{{\sigma_{\mathrm{vm}}^{n}}^2 - \tilde{\sigma}_Y^2}{d\tilde{\sigma}_Y^2} \right)
\end{equation*}
here $\sigma_{\mathrm{vm}}^{n}$ is the computed von Mises stress at the previous time step and $d=0.1$ denotes the nondimensional parameter.
The tangent stiffness matrix $\mathsf{C}_{\bt}$ in the linear region is first estimated as following, and then used as constant in \cref{EQ:CHOL_NN_2}.

\paragraph{Linear region}
The data sets are generated with load parameters 
$$(p_1, p_2, p_3) = (0.16, 0.016, 0.06)~\textrm{GN/m}$$
which are small enough to maintain the deformations in the linear region.

The indirect data training approach in~\cref{SEC:INDIRECT_DATA} is applied to extract the tangent stiffness matrix, and we obtain the following estimation 
\begin{equation}
 \mathsf{C}_{\theta}
 = \begin{bmatrix}
1.04064\times10^6 &   2.09077\times10^5      &   0.0 \\
      2.09077\times10^5 &  1.041146\times10^6  &     0.0\\
      0.0           &       0.0           &        4.19057\times10^5\\
\end{bmatrix}
\end{equation}
Based on Young's modulus and Poisson's ratio of the material, the tangent stiffness matrix is 
\begin{equation}
 \mathsf{C} 
 = \begin{bmatrix}
1.04167\times10^6 & 2.08333\times10^5  & 0.0 \\     
2.08333\times10^5 & 1.04167\times10^6  & 0.0 \\     
0.0             & 0.0     &   4.16667\times10^5 
\end{bmatrix}
\end{equation}
For each components, the relative error is less than one percent. These errors are introduced by the discretization, including the interpolation of the displacement field on the observation grid and the estimation of the acceleration~\cref{EQ:ACCELERATION}. It is worth mentioning, the direct input-output data training delivers exactly the same stiffness matrix as the reference.

\begin{figure}[htpb]
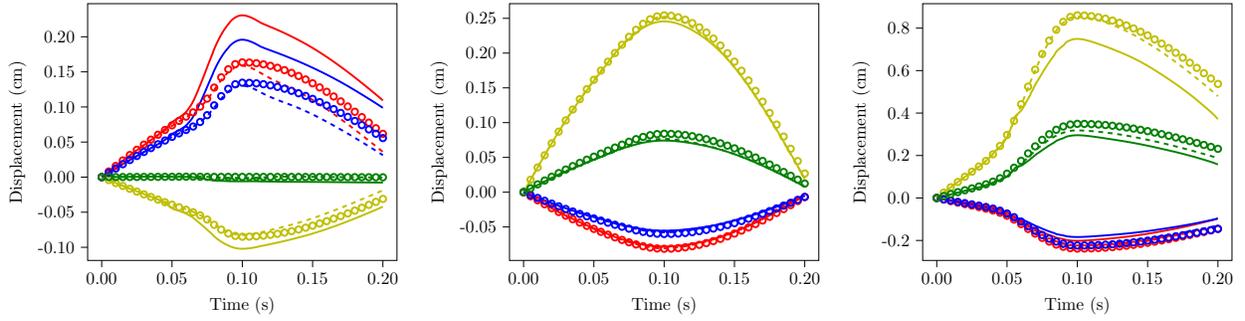

\centering
\scalebox{0.6}{\input{figures/plate_plasticity_disp_nn2_piecewise_from3_test106.tex}}~
\scalebox{0.6}{\input{figures/plate_plasticity_disp_nn2_piecewise_from3_test206.tex}}~
\scalebox{0.6}{\input{figures/plate_plasticity_disp_nn2_piecewise_from3_test300.tex}}
  \caption{Trajectories of displacement at top-right ~(red: $u_x$, yellow: $u_y$) and top-middle points~(blue: $u_x$, green: $u_y$) of the 2D elasto-plastic plate for test A6~(left), B6~(middle) and C1~(right), defined on page~\pageref{SEC:PLATE}.  The reference solutions are marked by empty circles, the solutions obtained by the CholNN trained using indirect data are marked by solid lines, and the solutions obtained by CholNN trained with direct data are marked by dashed lines.}
  \label{FIG:ELASTO_PLASTICITY_DISP}
\end{figure}

\begin{figure}[htpb]
\label{FIG:ELASTO_PLASTICITY_STRESS}
\centering
 \includegraphics[width=0.33\textwidth]{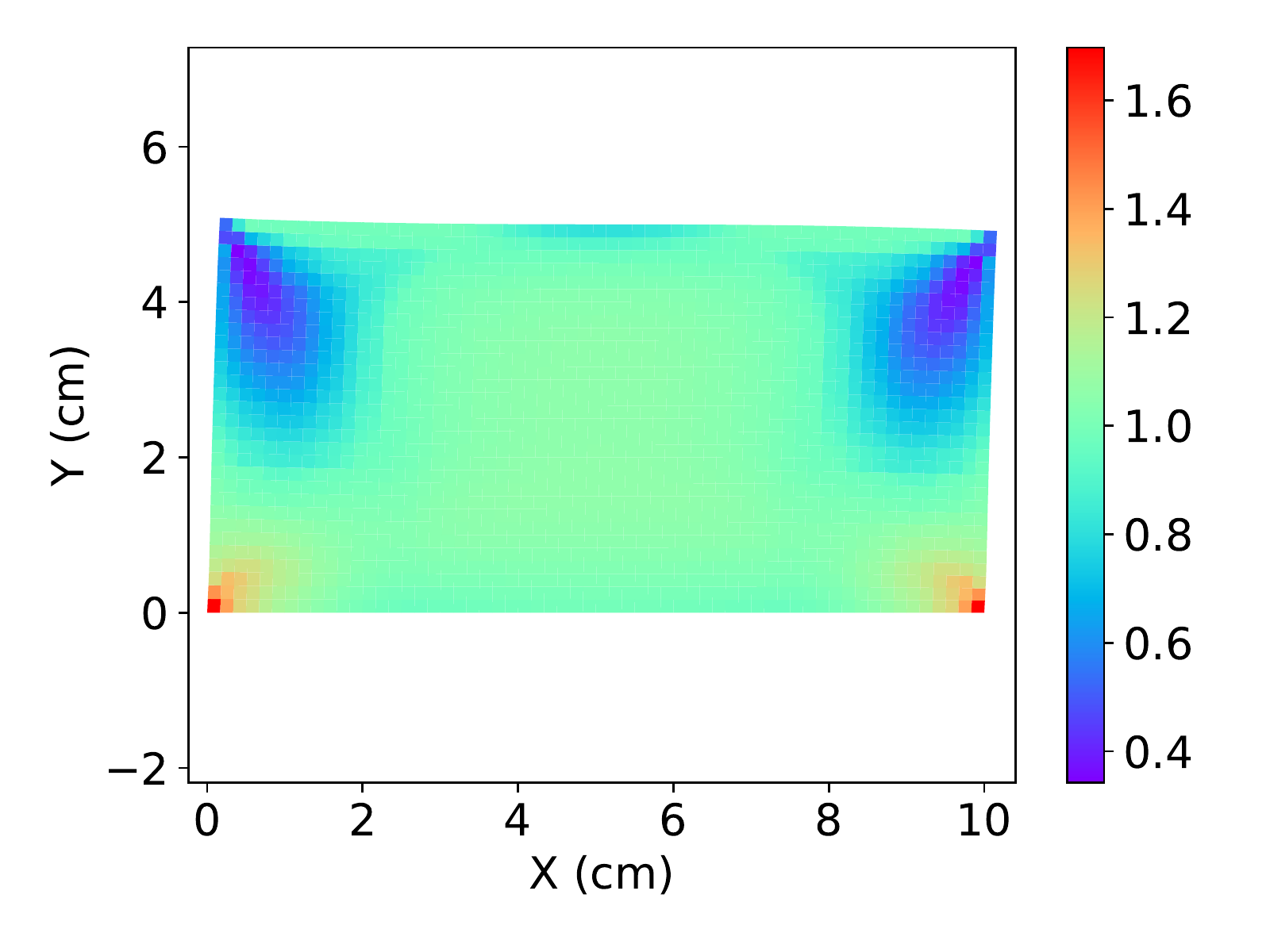}~
  \includegraphics[width=0.33\textwidth]{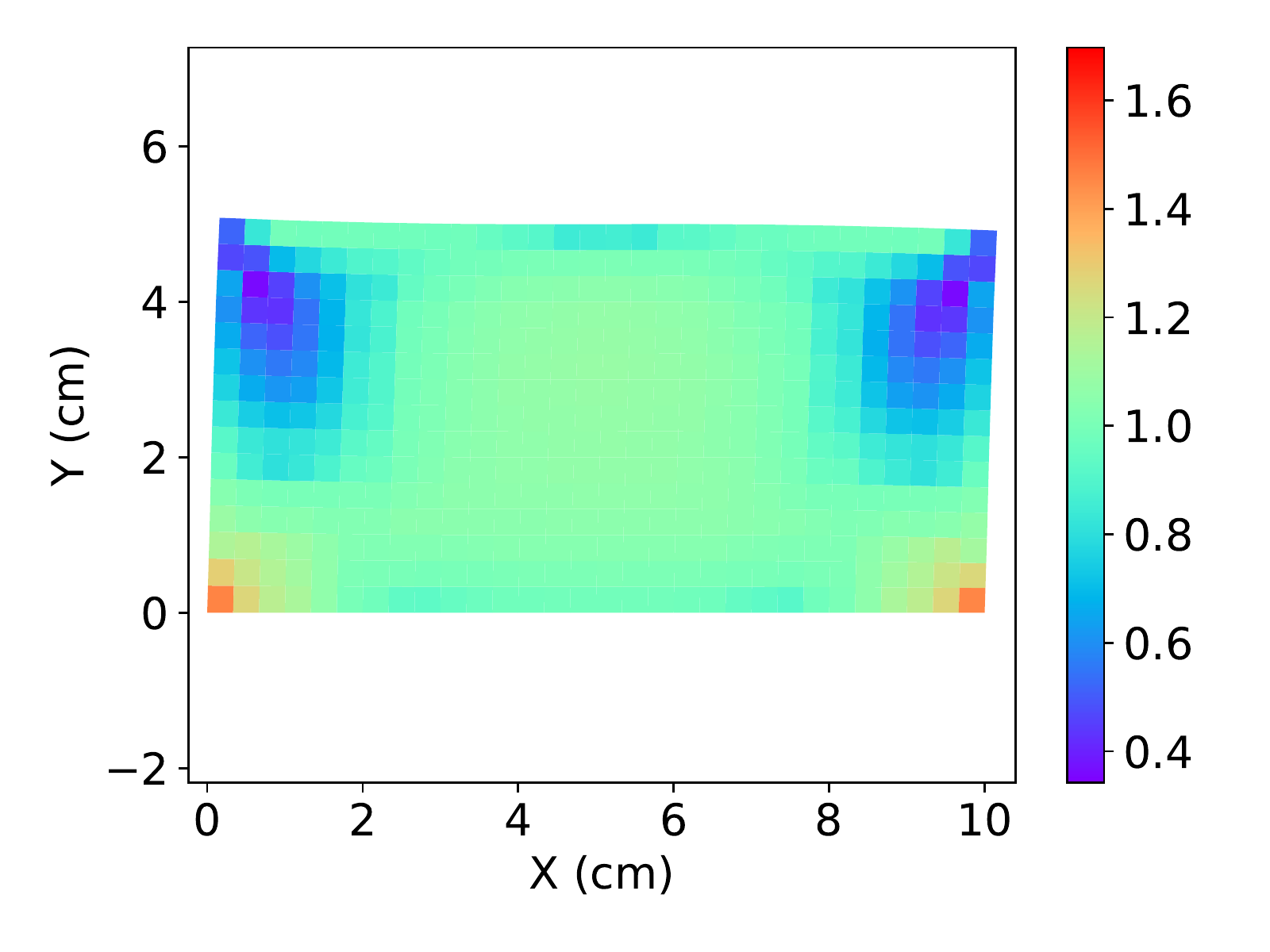}~
  \includegraphics[width=0.33\textwidth]{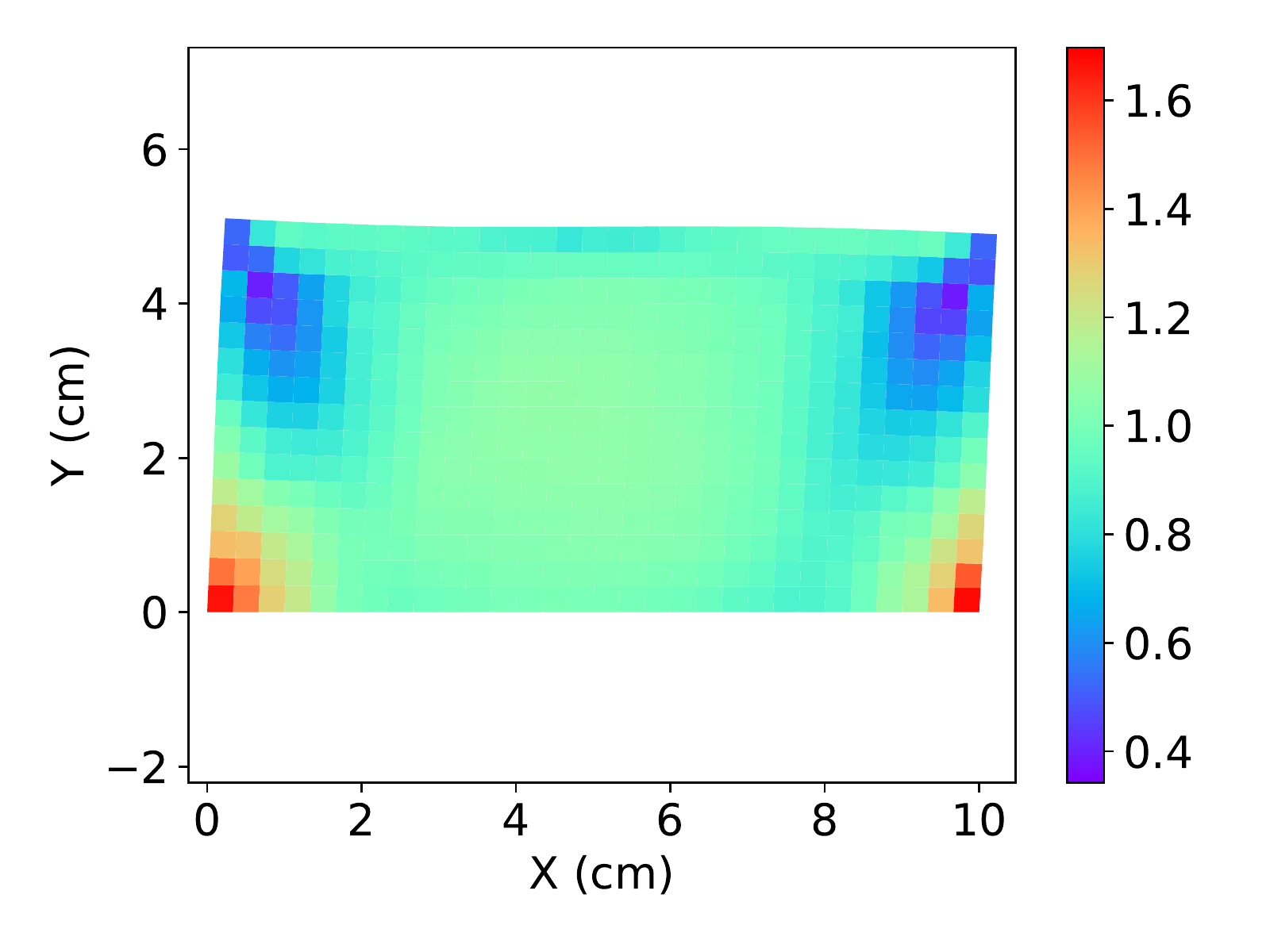}
  \includegraphics[width=0.33\textwidth]{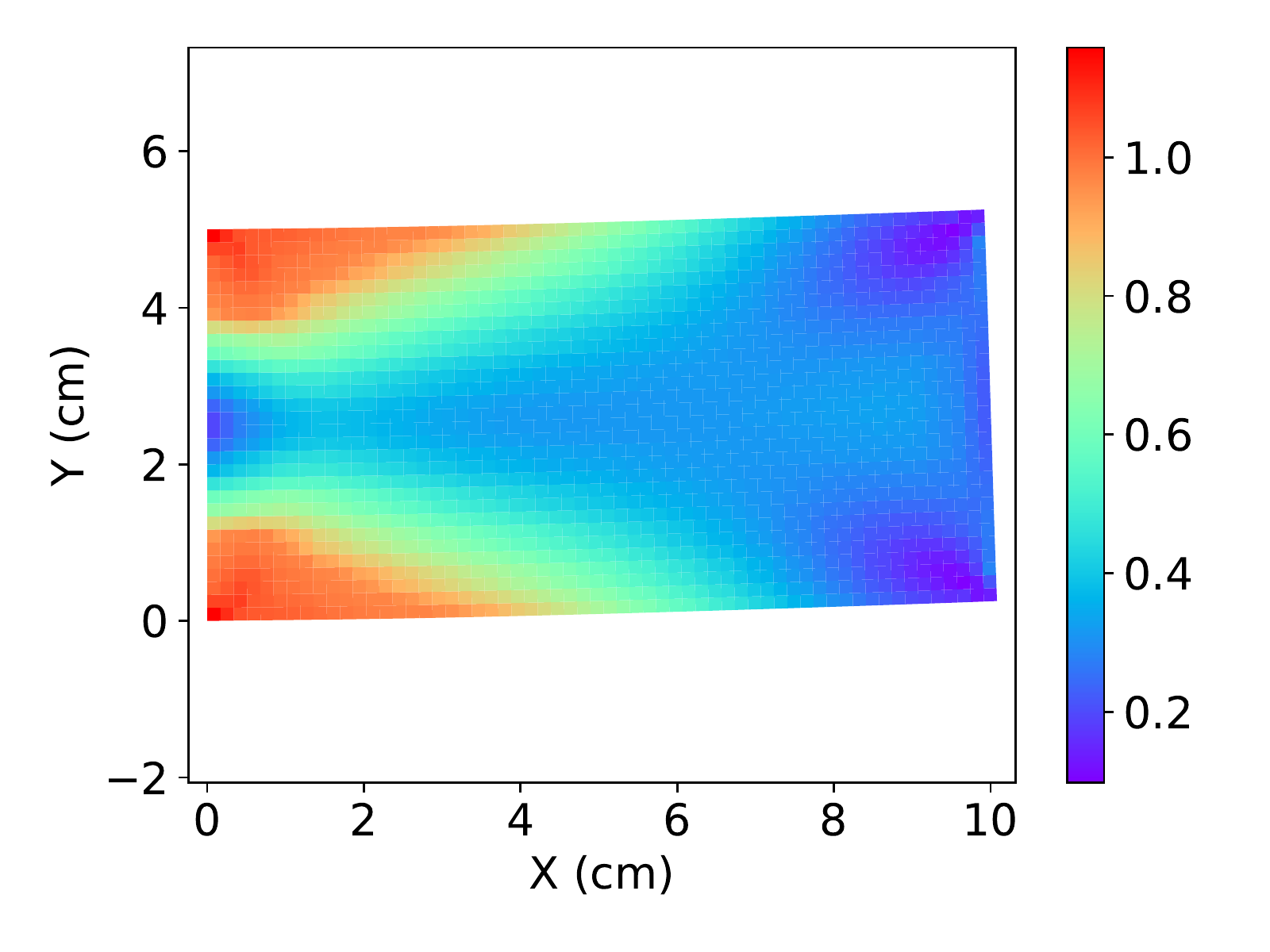}~
  \includegraphics[width=0.33\textwidth]{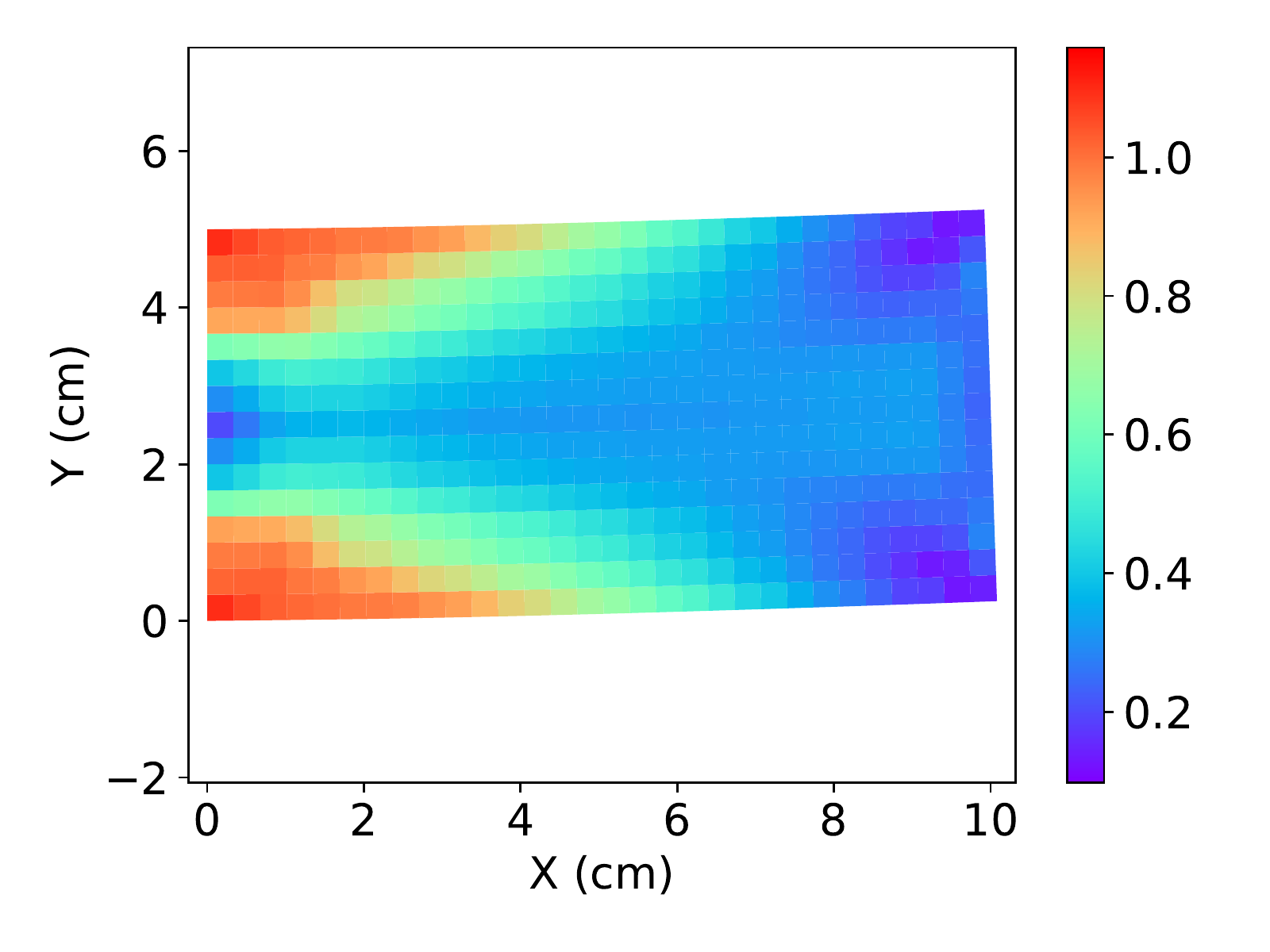}~
  \includegraphics[width=0.33\textwidth]{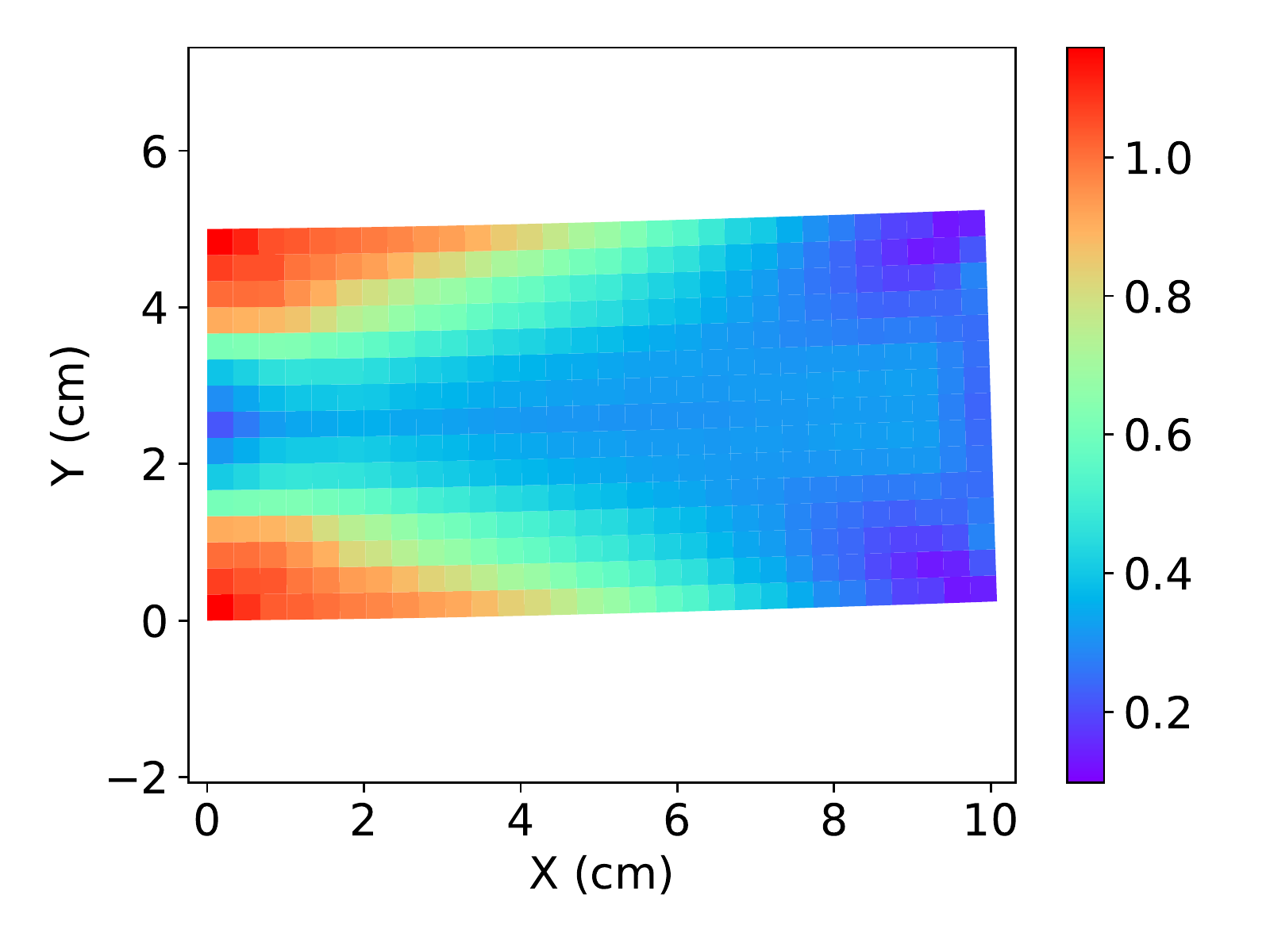}
\includegraphics[width=0.33\textwidth]{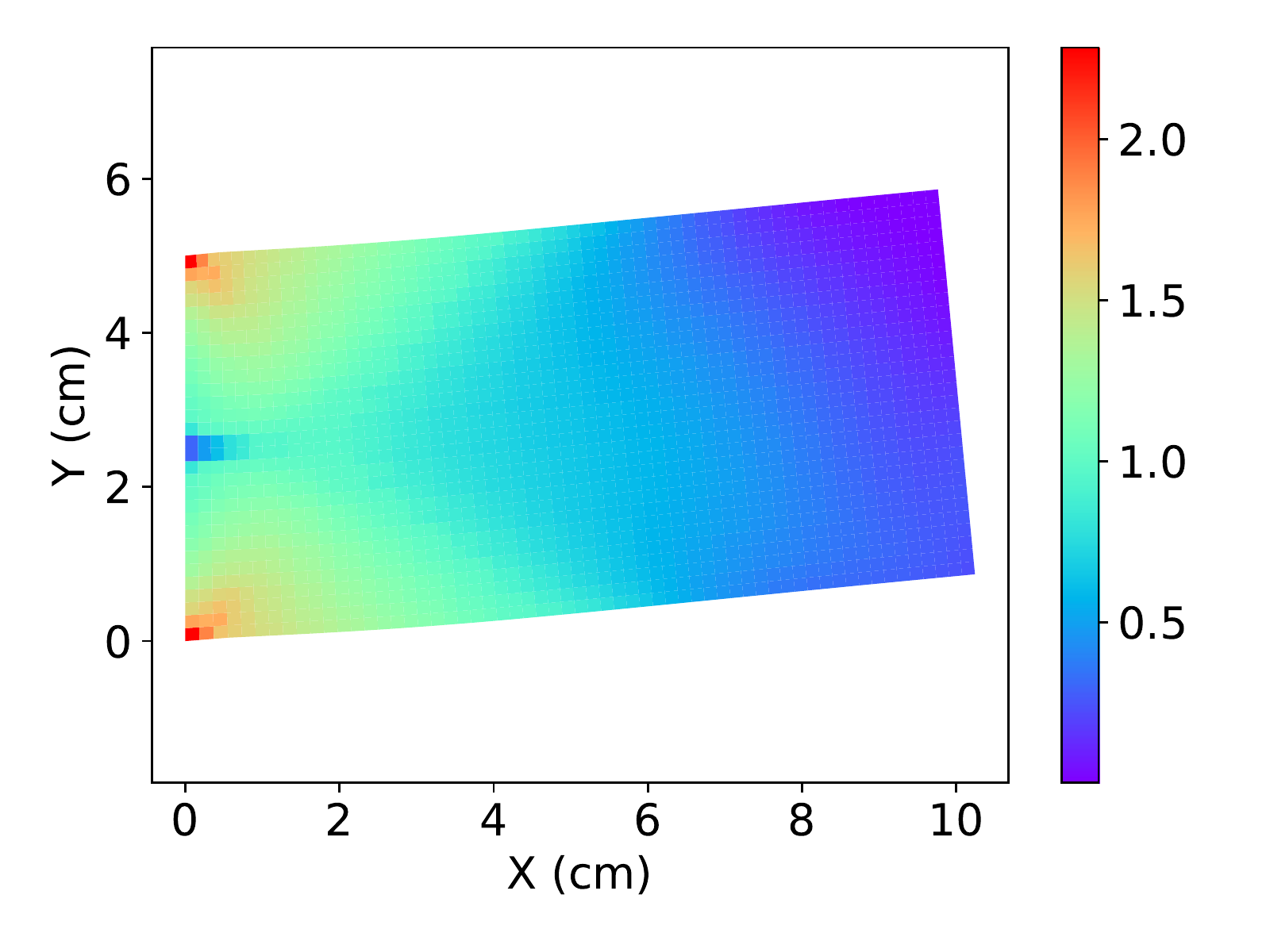}~
  \includegraphics[width=0.33\textwidth]{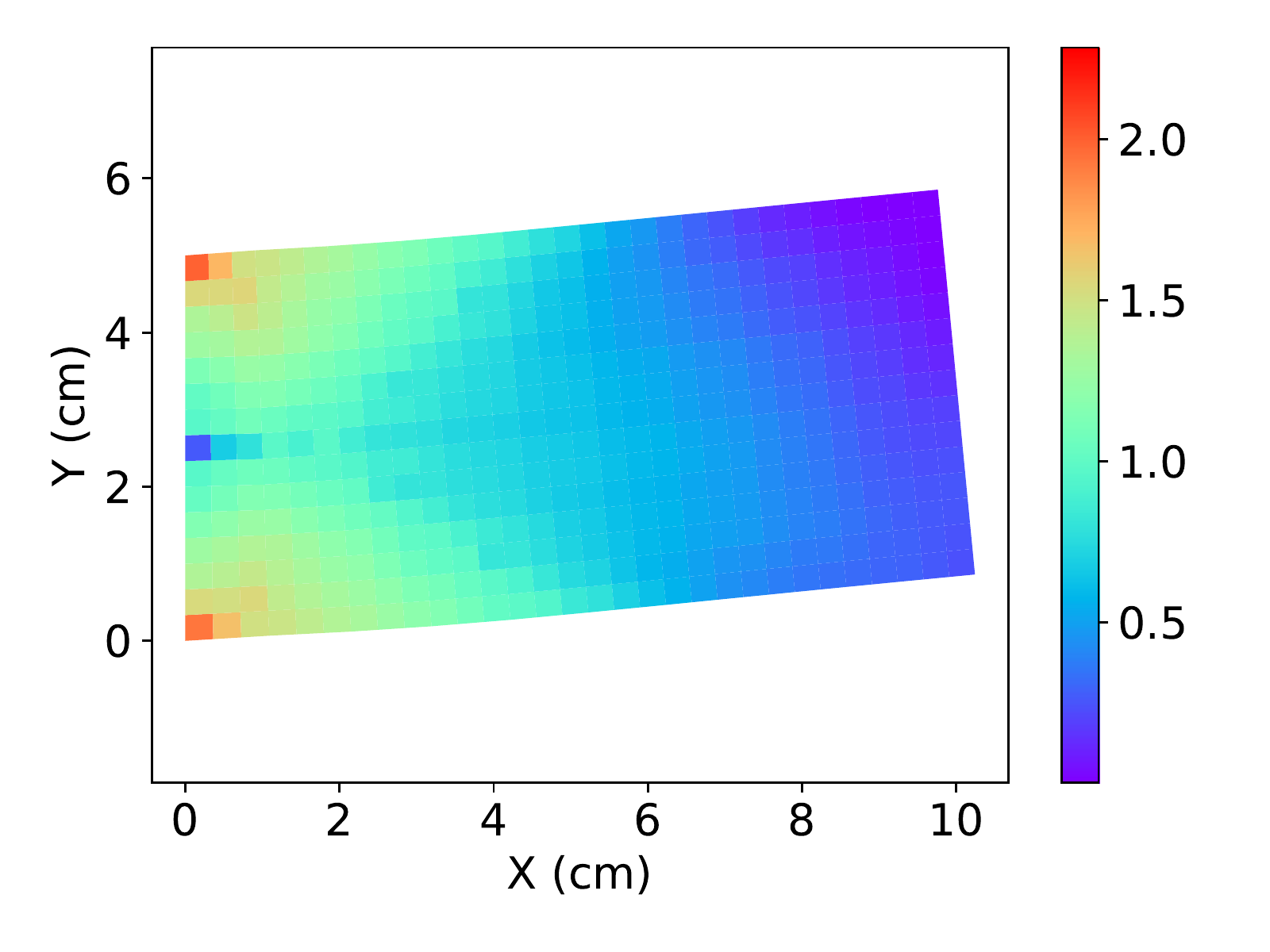}~
  \includegraphics[width=0.33\textwidth]{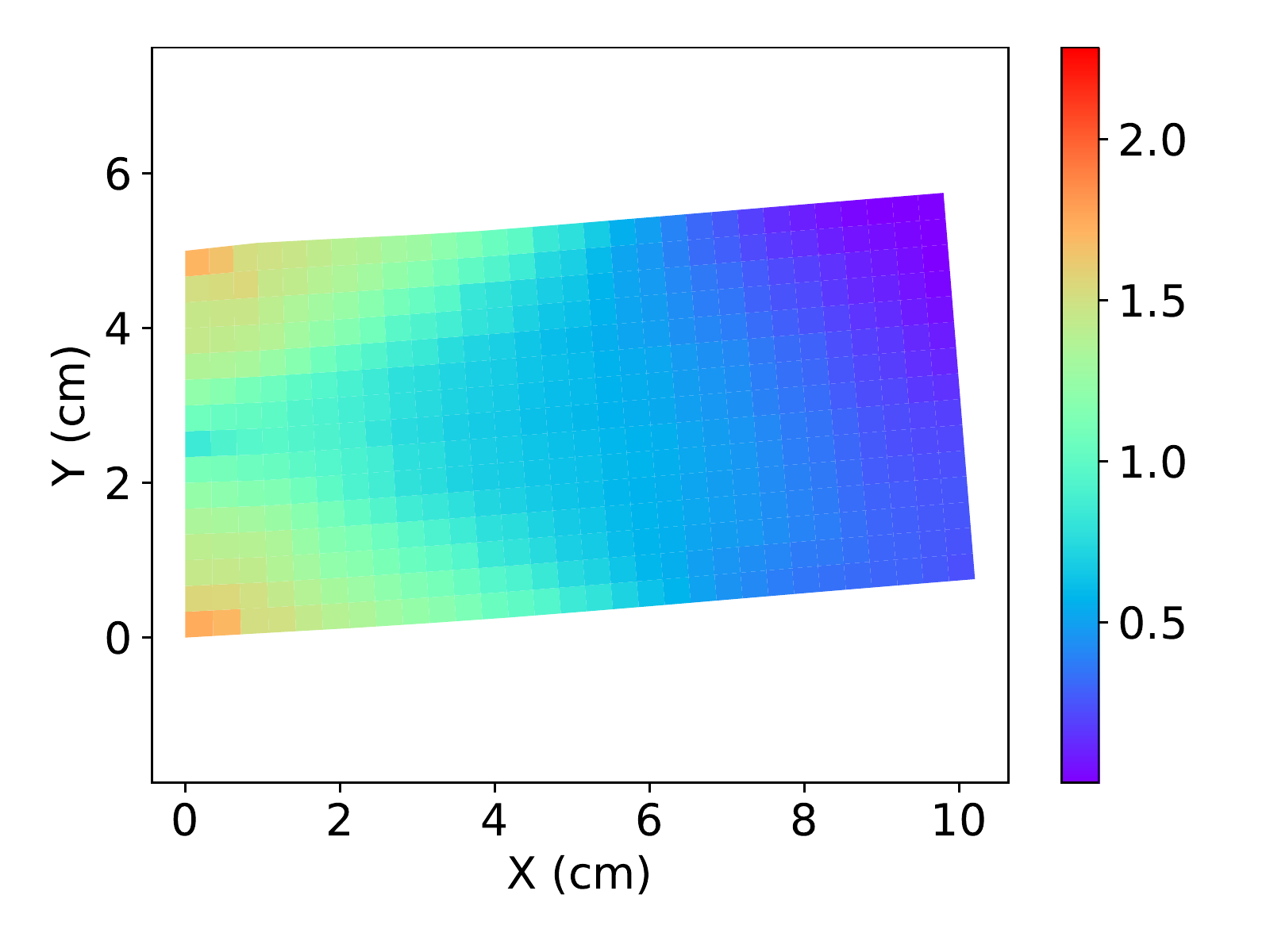}
  \caption{The von Mises stress~(GPa) fields at $t=\frac{T}{2}$ for the 2D elasto-plastic plate for test A6~(top), B6~(middle) and C1~(bottom), defined on page~\pageref{SEC:PLATE}.  From left to right: reference solutions~(on a fine mesh), solutions obtained by SPD-NN trained with direct data, solutions obtained by SPD-NN trained with indirect data.}
\end{figure}

\paragraph{Nonlinear region}
The data sets are generated with load parameters 
$$(p_1, p_2, p_3) = (1.6, 0.16, 0.6)~\textrm{GN/m}$$ 
which are 10 times larger than these in the linear region.
Both direct input-output data training~(\cref{SEC:DIRECT_DATA})  and indirect  data training~(\cref{SEC:INDIRECT_DATA}) are applied to train a SPD-NN  with 5 hidden layers and 20 neurons in each layer.
The predicted trajectories of the displacements at top-right and top-middle points for all test cases are depicted in \cref{FIG:ELASTO_PLASTICITY_DISP}, along with the references.  SPD-NNs trained with both methods are able to predict the initial elastic behavior, the strain-hardening region, and the unloading behavior.  
The SPD-NN obtained by the direct input-output data training performs better especially for the prediction of the yield strength the strain-hardening behavior. 
The predicted von Mises stress fields at $t=\frac{T}{2}$ and the references for all test cases are depicted in \cref{FIG:ELASTO_PLASTICITY_STRESS}. Reasonable agreements are achieved.

\subsubsection{Multiscale Fiber Reinforced Elasto-plasticity}
\label{SEC:MULTISCALE-PLATE}
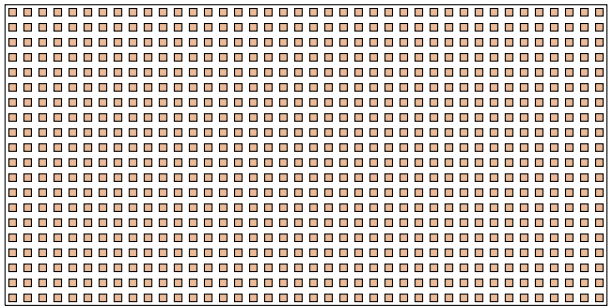
\begin{figure}[htbp]
\centering
\begin{tikzpicture}
  \draw (0,0) -- (8,0) -- (8,4) -- (0,4) -- (0,0);
\foreach \i in {0,...,39}{
    \foreach \j in {0,...,19}
    {
        \pgfmathsetmacro{\cubexl}{\i*0.2  + 0.05}
        \pgfmathsetmacro{\cubeyl}{\j*0.2  + 0.05}
        \pgfmathsetmacro{\dcubex}{0.1}
        \pgfmathsetmacro{\dcubey}{0.1}
        \draw[draw=black,fill=BurntOrange!40!white] (\cubexl,\cubeyl) -- ++(\dcubex,0) -- ++(0,\dcubey) -- ++(-\dcubex,0) -- cycle;
    }
  }
\end{tikzpicture}
    \caption{Schematic of the fiber (orange) reinforced thin plate.}
\label{FIG:PLATE_FIBERS}
\end{figure}

The plate is made of the titanium---the same as \cref{SEC:PLASTICITY-PLATE}---but reinforced by fibers made of SiC, which are assumed to be isotropic and elastic with
$$\rho = 3200~\textrm{kg/m}^3, \; E = 400~\textrm{GPa},\; \textrm{and} \; \nu = 0.35$$
These fibers are square shaped and uniformly distributed in the plate, with a diameter $d= 0.25~\textrm{cm}$ and a fraction $25\%$. There are in total 800 fibers, as shown in \cref{FIG:PLATE_FIBERS}.  This plate domain is discretized by $200 \times 400$ quadratic quadrilateral elements, and 25 elements for each fiber. And the time step size is $\Delta t = 0.001s$.

As for training SPD-NNs, the estimated yield strength and the transition function are the same as these in~\cref{SEC:PLASTICITY-PLATE}. The tangent stiffness matrix $\mathsf{C}_{\bt}$ in the linear region is first calibrated and then fixed as constant in \cref{EQ:CHOL_NN_2} when training the SPD-NN.  

\paragraph{Linear region}

\begin{figure}[htpb]
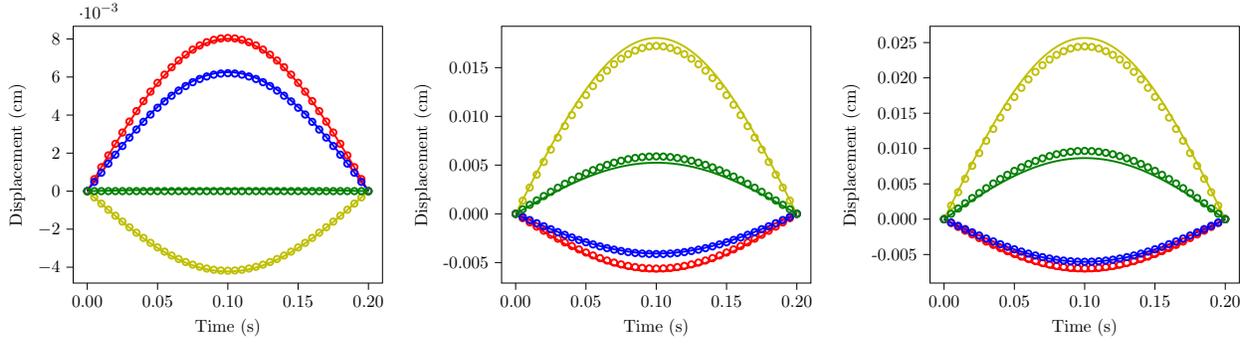

\label{FIG:MULTISCALE_LINEAR_DISP}
\centering
\scalebox{0.6}{\input{figures/plate_multiscale_disp_linear_106.tex}}~
\scalebox{0.6}{\input{figures/plate_multiscale_disp_linear_206.tex}}~
\scalebox{0.6}{\input{figures/plate_multiscale_disp_linear_300.tex}}~
  \caption{Trajectories of displacement at  top-right~(red: $u_x$, yellow: $u_y$) and top-middle points~(blue: $u_x$, green: $u_y$) of the 2D multiscale plate in the linear region for test A6~(left), B6~(middle) and C1~(right), defined on page~\pageref{SEC:PLATE}.  The reference solutions are marked by empty circles and the solutions obtained by the SPD-NN trained with indirect data are marked by solid lines.}
\end{figure}

\begin{figure}[htpb]
\label{FIG:MULTISCALE_LINEAR_STRESS}
\centering
 \includegraphics[width=0.33\textwidth]{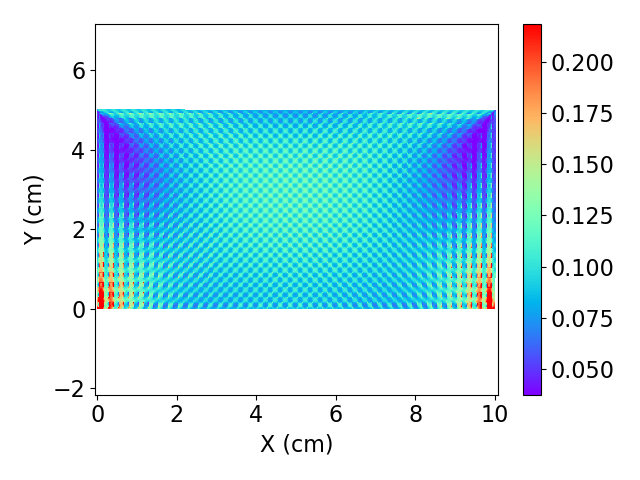}~
  \includegraphics[width=0.33\textwidth]{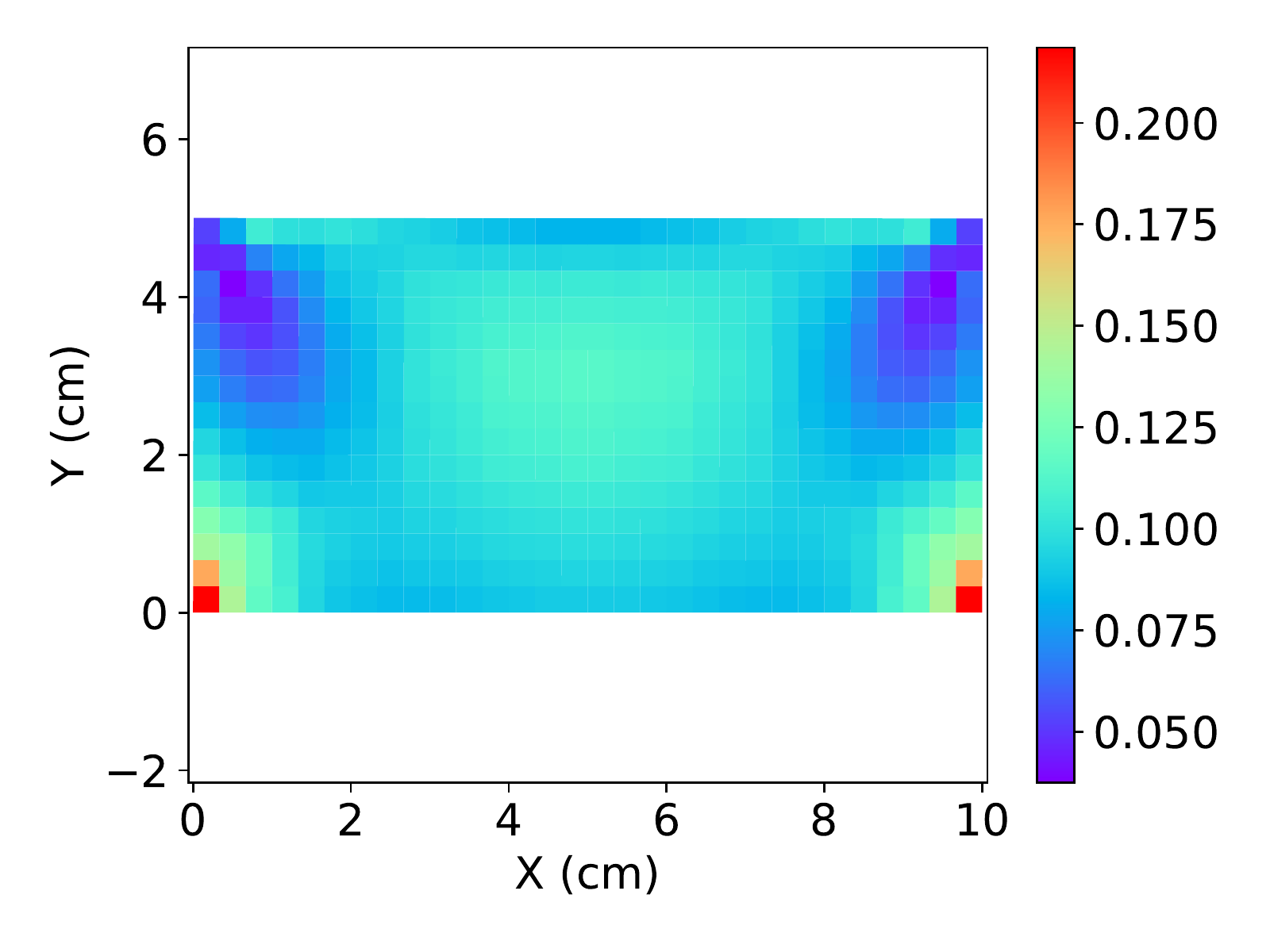}
  
  \includegraphics[width=0.33\textwidth]{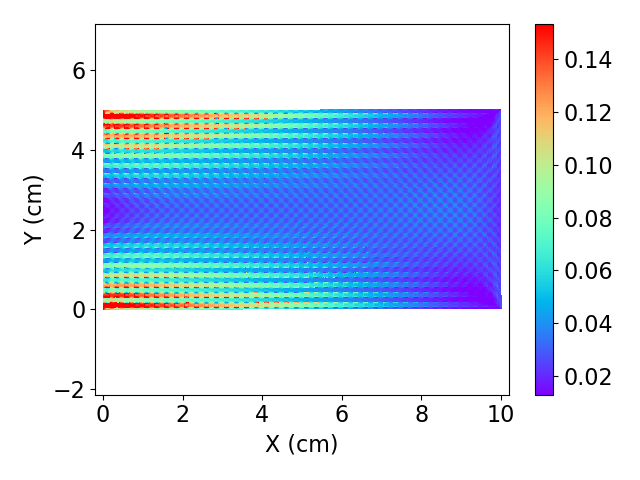}~
  \includegraphics[width=0.33\textwidth]{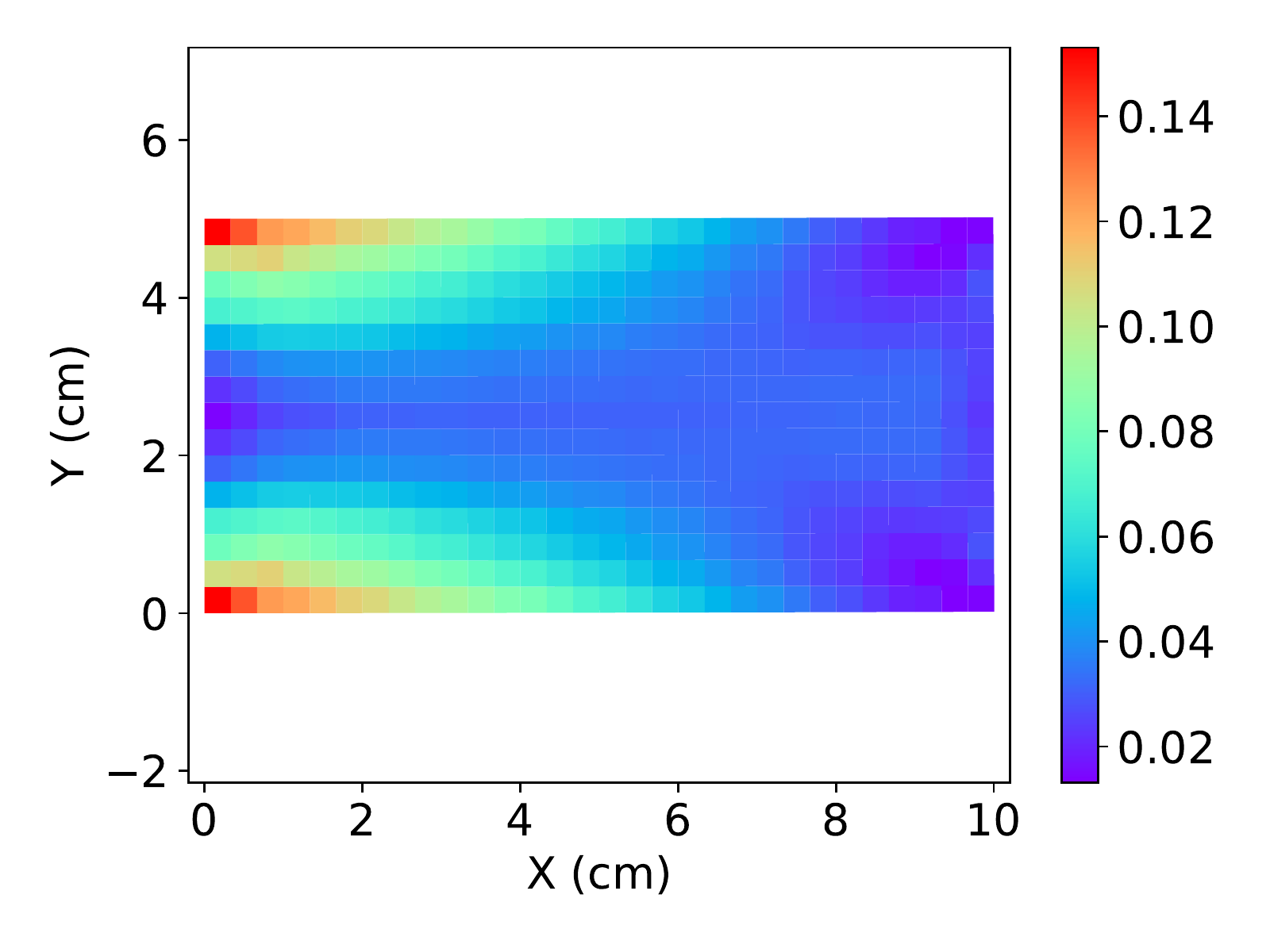}

\includegraphics[width=0.33\textwidth]{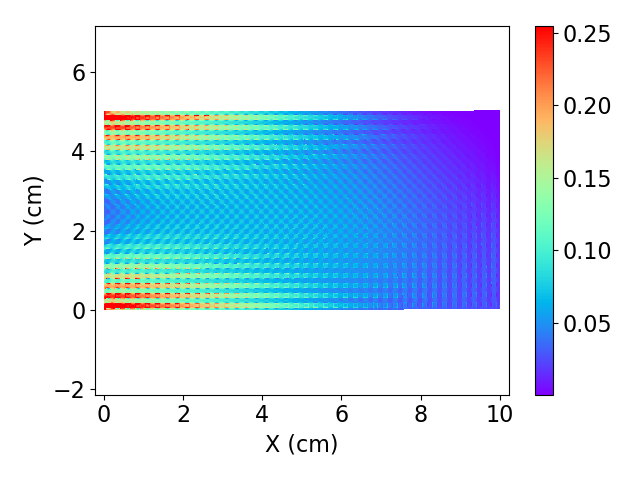}~
  \includegraphics[width=0.33\textwidth]{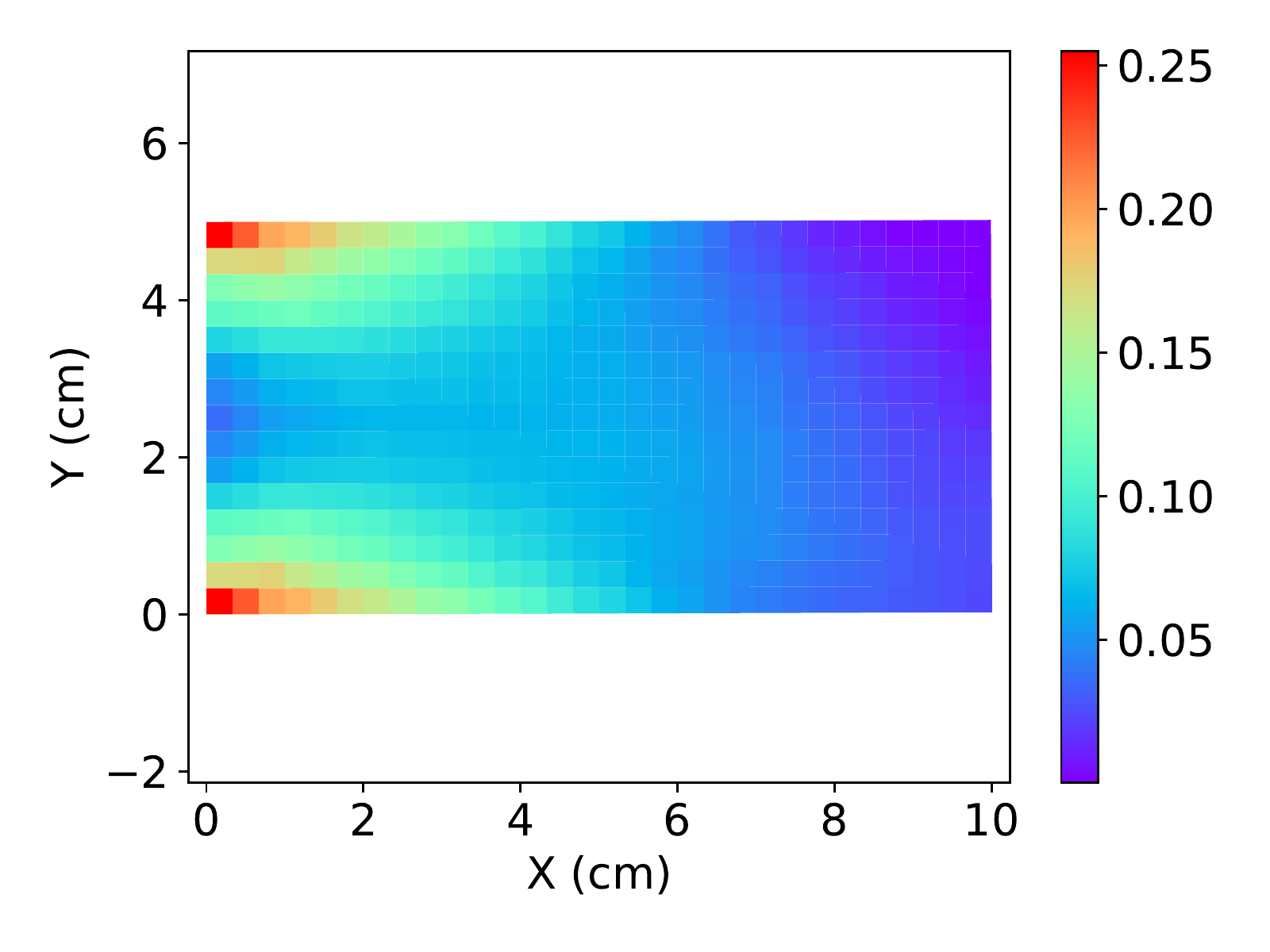}

  \caption{The von Mises stress~(GPa) fields at $t=\frac{T}{2}$ for the 2D multiscale plate in the linear region for test A6~(top), B6~(middle) and C1~(bottom), defined on page~\pageref{SEC:PLATE}.  From left to right: reference solutions~(on a fiber-resolved mesh) and solutions obtained by SPD-NN trained with indirect data.}
\end{figure}

The data sets are generated with load parameters 
$$(p_1, p_2, p_3) = (0.16, 0.016, 0.06)~\textrm{GN/m}$$
which are small enough to maintain the deformations in the linear region. 
The indirect data training approach in~\cref{SEC:INDIRECT_DATA} is applied to extract the following predicted tangent stiffness matrix,
\begin{equation}
 \mathsf{C}_{\theta}
 = \begin{bmatrix}
1.335174\times10^6 & 3.26448\times 10^5 &0.0  \\
3.26448\times10^5 & 1.326879\times 10^6 &0.0 \\
0.0         & 0.0 & 5.26955\times 10^5 \\
\end{bmatrix}
\end{equation}
The predicted linear constitutive relation~\cref{EQ:LINEAR} is verified on the test set. 
The predicted displacements at top-right and top-middle points as a function of time and the references for all test cases are depicted in \cref{FIG:MULTISCALE_LINEAR_DISP}. The corresponding von Mises stress fields at $t=\frac{T}{2}$ are reported in \cref{FIG:MULTISCALE_LINEAR_STRESS}.  The SPD-NN based homogenized model delivers similar results as the high-resolution multiscale model.

\paragraph{Nonlinear region}
The data sets are generated with load parameters 
$$(p_1, p_2, p_3) = (1.6, 0.16, 0.6)~\textrm{GN/m}$$ 
which are 10 times larger than those in the linear region.
Only the indirect data training approach is applied to train a SPD-NN with 5 hidden layers and 20 neurons in each layer. To enable direct input-output data training, homogenization is required to generate strain-stress data from RVE simulations~\cite{le2015computational, bessa2017framework, ling2016machine, wang2018multiscale}, or extract strain-stress data from direct numerical simulations. This may be challenging in the context of experimental data, which necessitates the indirect training approach.  
   
\begin{figure}[htpb]
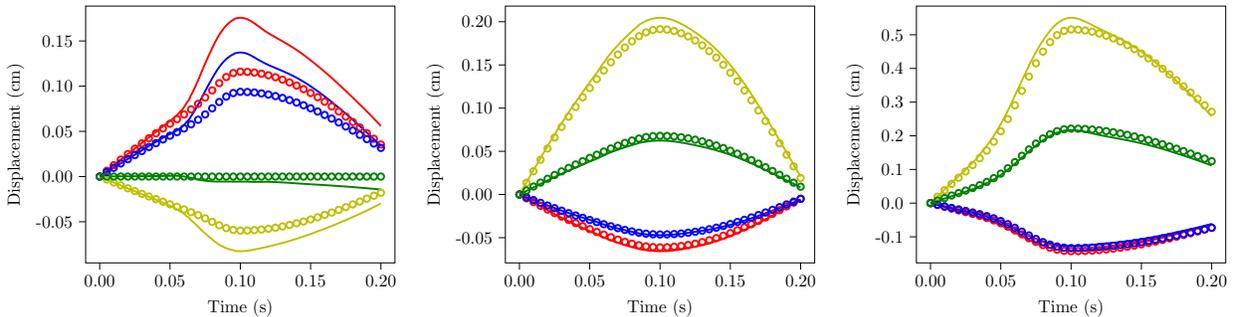

\centering
\scalebox{0.6}{\input{figures/plate_multiscale_disp_nn2_piecewise_from3_test106}}~
\scalebox{0.6}{\input{figures/plate_multiscale_disp_nn2_piecewise_from3_test206}}~
\scalebox{0.6}{\input{figures/plate_multiscale_disp_nn2_piecewise_from3_test300}}~
  \caption{Trajectories of displacement at top-right ~(red: $u_x$, yellow: $u_y$) and top-middle points~(blue: $u_x$, green: $u_y$) of the 2D multiscale plate in nonlinear region for test A6~(left), B6~(middle) and C1~(right), defined on page~\pageref{SEC:PLATE}.  The reference solutions are  marked by empty circles and the solutions obtained by the SPD-NN trained with indirect data are marked by solid lines.}
    \label{FIG:MULTISCALE_NONLINEAR_DISP}
\end{figure}

\begin{figure}[htpb]
\centering
 \includegraphics[width=0.33\textwidth]{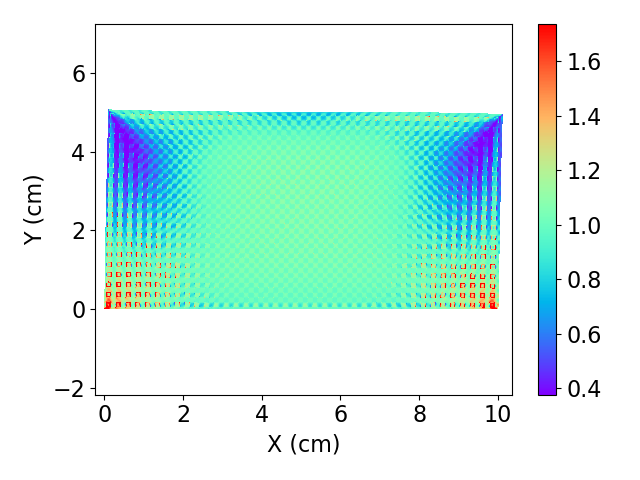}~
  \includegraphics[width=0.33\textwidth]{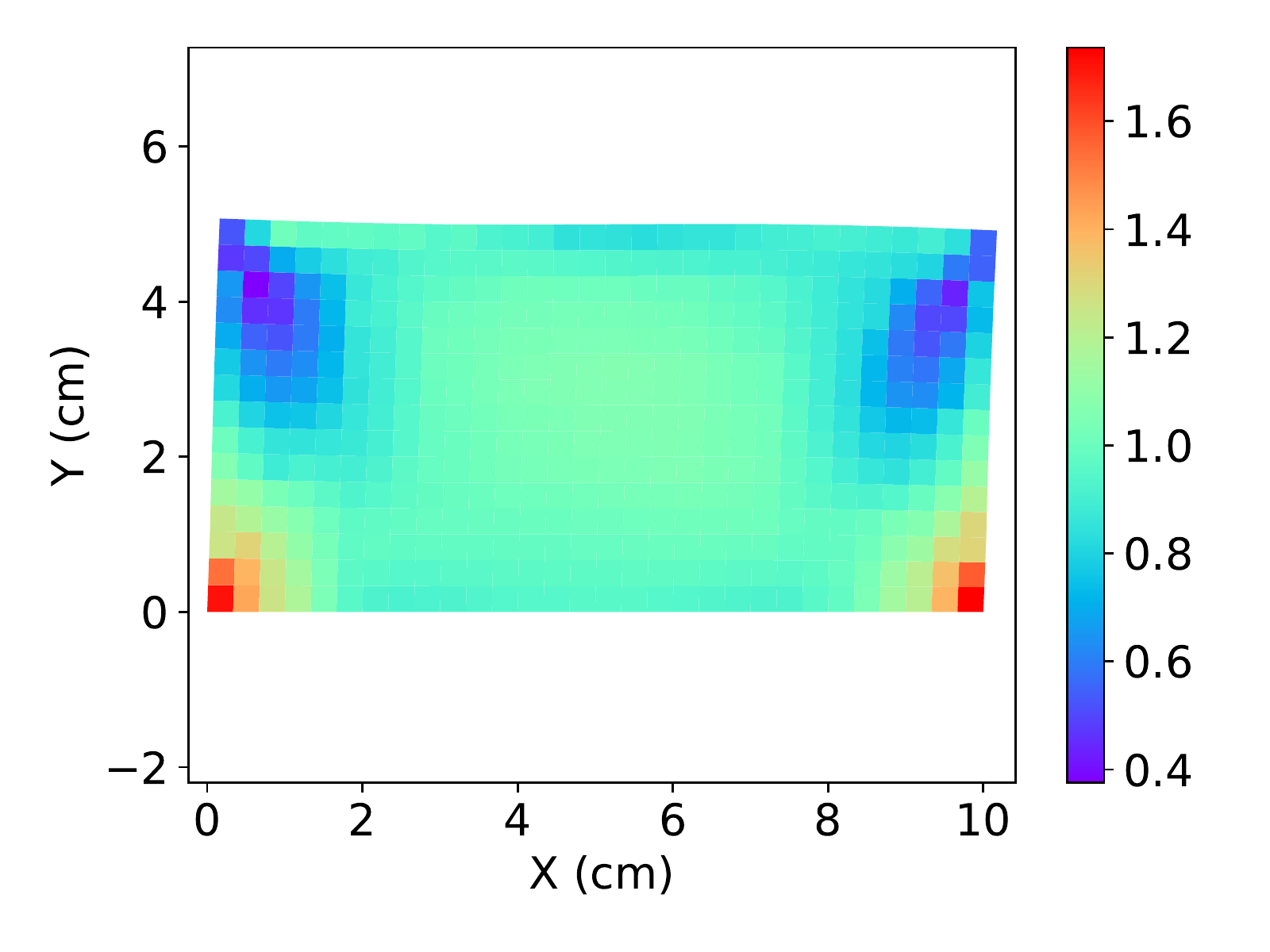}
  
  \includegraphics[width=0.33\textwidth]{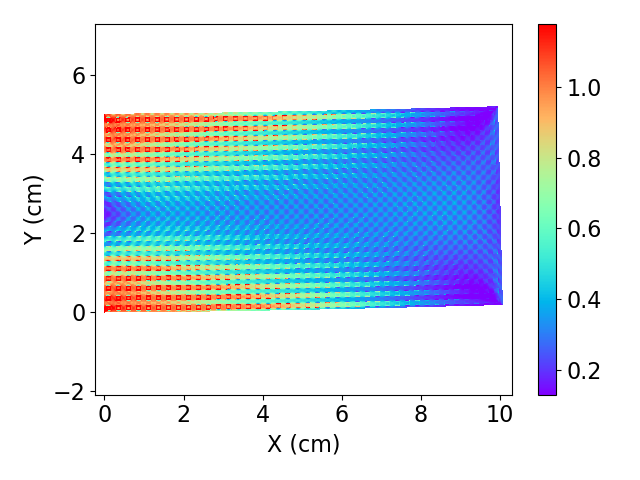}~
  \includegraphics[width=0.33\textwidth]{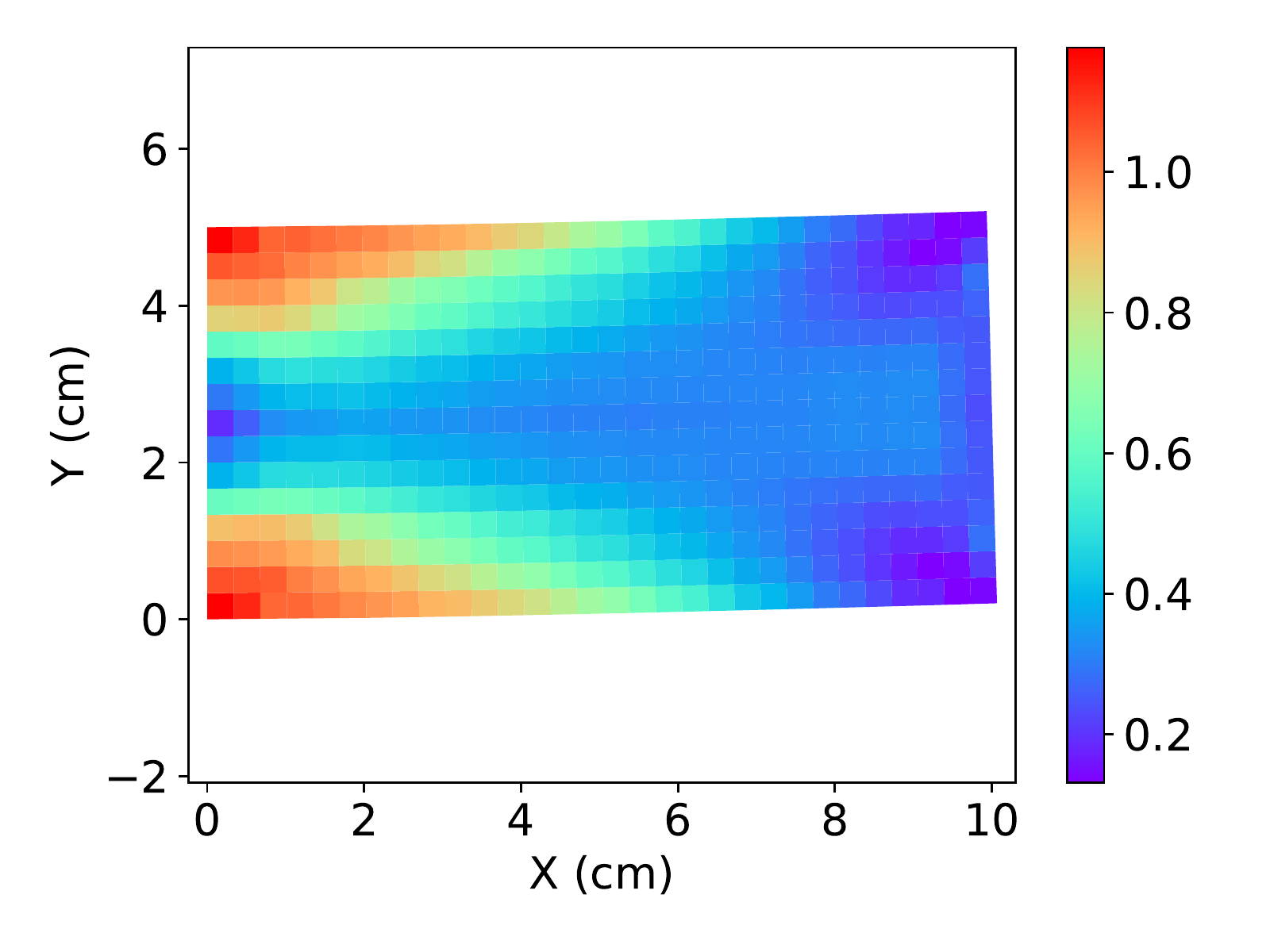}

  \includegraphics[width=0.33\textwidth]{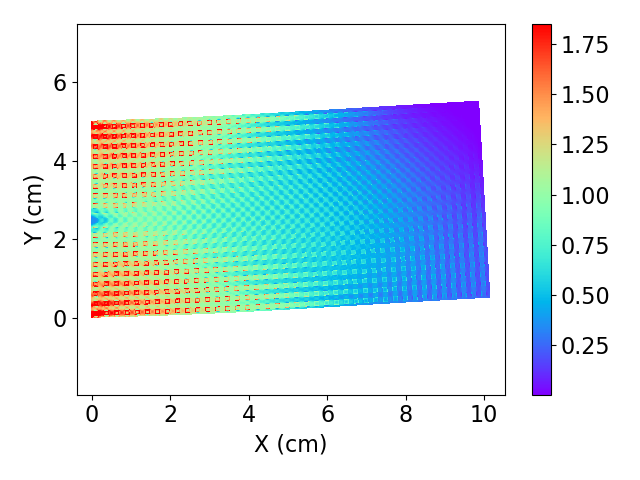}~
  \includegraphics[width=0.33\textwidth]{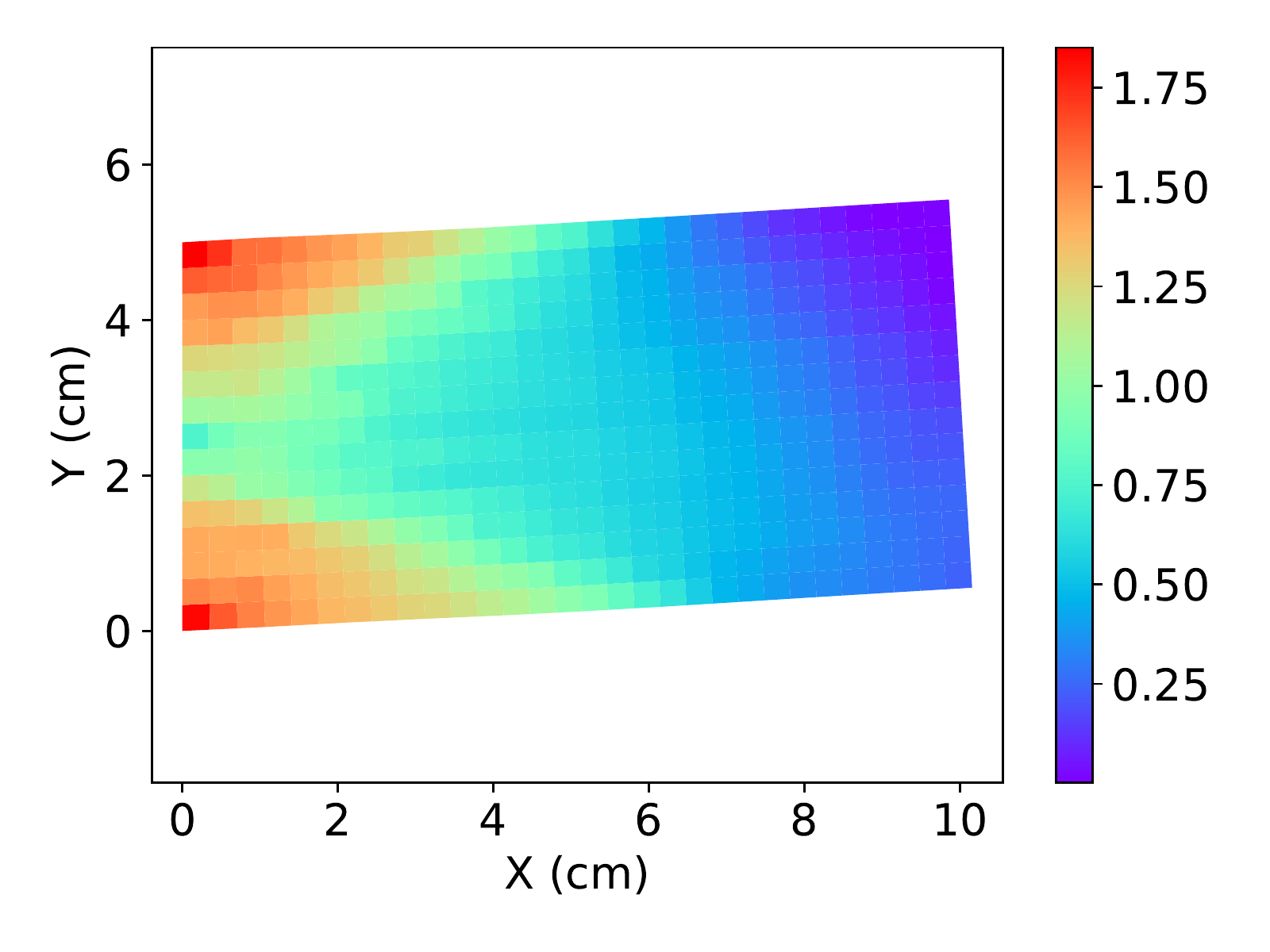}

  \caption{The von Mises stress~(GPa) fields at $t=\frac{T}{2}$ for the 2D multiscale plate in nonlinear region for test A6~(top), B6~(middle) and C1~(bottom), defined on page~\pageref{SEC:PLATE}.  From left to right: reference solutions~(on a fiber-resolved mesh) and solutions obtained by SPD-NN trained with indirect data.}
  \label{FIG:MULTISCALE_NONLINEAR_STRESS}
\end{figure}

The predicted trajectories of displacements at top-right and top-middle points for all test cases are depicted in \cref{FIG:MULTISCALE_NONLINEAR_DISP}, along with the references. 
Comparing with the previous elaso-plasticity case, all displacements are smaller due to the SiC fiber reinforcement.
The proposed SPD-NN gives a satisfactory approximation of the initial elastic behavior, the strain-hardening region, and the following unloading behavior.  
The predicted von Mises stress fields at $t=\frac{T}{2}$ and the references for all test cases are depicted in \cref{FIG:MULTISCALE_NONLINEAR_STRESS}. The predicted and simulated homogenized stress fields are in reasonably good agreement.  
Although the solutions obtained by SPD-NNs do not capture local large stress concentrations near each fiber (at the level of the microstructure), 
the local recovery techniques~\cite{kanoute2009multiscale} can be applied to estimate these local stress concentrations.
For example, a RVE simulation can be conducted to estimate local stresses  using the local strain from the coarse homogenized solution.

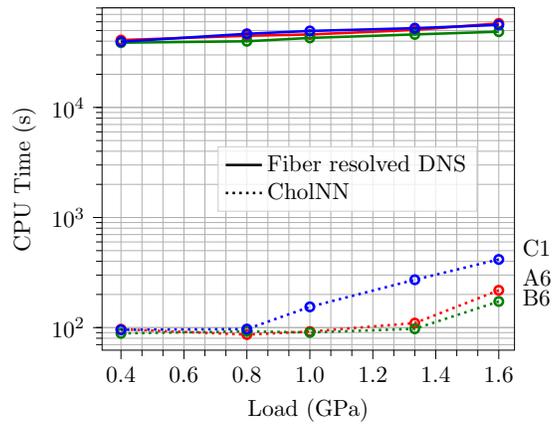
\begin{figure}[htpb]
\centering
	\scalebox{0.8}{
\begin{tikzpicture}

\begin{axis}[
legend cell align={left},
legend style={fill opacity=0.8, draw opacity=1, text opacity=1, at={(0.91,0.5)}, anchor=east, draw=white!80!black},
log basis y={10},
tick align=outside,
tick pos=left,
x grid style={white!69.0196078431373!black},
xlabel={Load (GPa)},
minor tick num=2,
xmajorgrids,
xminorgrids,
xmin=0.34, xmax=1.65,
clip mode=individual,
xtick style={color=black},
xtick={0.2,0.4,0.6,0.8,1,1.2,1.4,1.6,1.8},
xticklabels={0.2,0.4,0.6,0.8,1.0,1.2,1.4,1.6,1.8},
y grid style={white!69.0196078431373!black},
ylabel={CPU Time (s)},
ymajorgrids,
yminorgrids,
ymin=62.3649238054306, ymax=80369.0998834116,
ymode=log,
ytick style={color=black},
]
\addplot [very thick, red, mark=o, mark size=2, mark options={solid}, forget plot]
table {%
1.6 58038.837961276
1.33333333333333 50593.067638013
1 46014.154744781
0.8 44768.697642296
0.4 40926.322249248
};
\addplot [very thick, red, dotted, mark=o, mark size=2, mark options={solid}, forget plot]
table {%
1.6 218.205717485
1.33333333333333 109.871857105
1 91.967749878
0.8 86.359633766
0.4 96.631165397
};
\addplot [very thick, green!50!black, mark=o, mark size=2, mark options={solid}, forget plot]
table {%
1.6 48749.580876096
1.33333333333333 46044.787294992
1 42793.291550977
0.8 39900.831641923
0.4 38758.784705406
};
\addplot [very thick, green!50!black, dotted, mark=o, mark size=2, mark options={solid}, forget plot]
table {%
1.6 172.72724317
1.33333333333333 97.451604466
1 90.329985264
0.8 92.315171475
0.4 88.666778581
};
\addplot [very thick, blue, mark=o, mark size=2, mark options={solid}, forget plot]
table {%
1.6 56211.987075903
1.33333333333333 52497.121672471
1 49507.305923038
0.8 46675.113248691
0.4 39564.816549516
};
\addplot [very thick, blue, dotted, mark=o, mark size=2, mark options={solid}, forget plot]
table {%
1.6 415.933555361
1.33333333333333 272.028080746
1 154.317654317
0.8 96.965845698
0.4 95.472751521
};
\addplot [very thick, black]
table {%
1.6 1
1.33333333333333 nan
1 nan
0.8 nan
0.4 nan
};
\addlegendentry{Fiber resolved DNS}
\addplot [very thick, black, dotted]
table {%
1.6 1
1.33333333333333 nan
1 nan
0.8 nan
0.4 nan
};
\addlegendentry{CholNN}
\draw (axis cs:1.65,242.705838839) node[
  scale=1.0,
  anchor=base west,
  text=black,
  rotate=0.0
]{A6};
\draw (axis cs:1.65,160.337932452) node[
  scale=1.0,
  anchor=base west,
  text=black,
  rotate=0.0
]{B6};
\draw (axis cs:1.65,451.833125732) node[
  scale=1.0,
  anchor=base west,
  text=black,
  rotate=0.0
]{C1};
\end{axis}

\end{tikzpicture}}
  \caption{CPU time for simulations using fiber resolved DNS (direct numerical simulation) and the SPD-NN surrogate model. The same color denotes the same test case (C1, A6, or B6).}
  \label{fig:cpu}
\end{figure}

\paragraph{Accelerating Simulations with Neural Network Surrogates for Constitutive Modeling} It is also worth noting that each SPD-NN-based simulation is several order magnitude faster than the corresponding fiber resolved simulations. The CPU time for both the SPD-NN-based simuations and fiber resolved simulations are shown in \cref{fig:cpu}. Note as we increase the external load, the CPU time increases because more elements undergo plastic deformations, which requires more expensive Raphson-Newton iterations in the numerical simulations. Still, the dramatic acceleration from around 24 hours to just a few minutes is impressive. However, the acceleration should be carefully interpreted in the context of surrogate models. First, the current benchmarks are based on serial execution, where the state-of-the-art FEM simulations usually involve parallelization. Nevertheless, these parallelization techniques for fiber resolved simulations are also directly applicable to SPD-NNs. Second, the speed of the fiber resolved simulations depend on the required resolution (e.g., the resolution on each fiber). If we use very coarse grids, the speed of the fiber resolved simulations may be comparable to or even faster than the SPD-NN-based simulations, though the coarse grids raise concerns on accuracy. Third, we have not compared the present model with other state-of-the-art accelerating techniques for multiscale modeling, for example projection-based model order reduction~\cite{zahr2017multilevel} and self-consistent clustering~\cite{liu2016self}. Finally, as with any NN-based surrogate models, the SPD-NN only works on test data that does not deviate too much from the training data, where fiber resolved simulations are usually considered to be general and applicable in a much wider context.

\section{Discussion on Neural Networks}
\label{SEC:DISCUSS_NN}
The introduction and benchmarking of SPD-NNs raise many questions, some of which are partially answered in this paper while some others are worthwhile to be investigated further. We believe that the following discussion is important and central to the application of NN-based constitutive modeling and should be careful dealt for the development of SPD-NNs and other NN-based approaches. 

\paragraph{Function approximators} In the present work, neural network is used as a basis function to approximate a complex constitutive relation, the strain-stress relation. 
Many other basis functions, including piecewise linear functions, radial basis functions, and radial basis function networks, 
can also be applied for the approximation of the Cholesky factor. The choice of neural networks is well-versed and justified by the following three reasons. 
\begin{itemize}
    \item First, the distribution of the strain data---the input to the neural network---is determined by experiment records, which in general does not evenly spread over the domain. The comparison study presented in \cite{huang2019predictive} illustrates that  
neural network outperforms the other basis functions in terms of regularization and generalization properties when the data distribution is ill-behaved.
\item Second,  the strain-stress curves are non-smooth in the context of plastic deformations and elastic unloading. As a result of the non-smoothness, the approximation efficiency of typical basis functions is usually compromised. Nevertheless, neural networks exhibit the potential to capture the sharp transitions in the strain-stress relations and performs reasonably well for the non-smooth data.
\item Third, the input and the output dimensions are relatively high (e.g., the input is 9D and the output is 4D in the elasto-plasticity and the multi-scale cases), which poses a chanllenge for traditional basis functions. For example, if we were to use the linear basis functions, to discretize the high dimensional input space, even if we only have 10 grid points per dimension, the total degrees of freedom is $10^9$, which is too costly regarding both the computation and the storage. Yet, neural networks are particularly convenient and useful for expressing mappings between high dimensional spaces~\cite{han2018solving} and requires no mesh in the input parameter domain. 
\end{itemize}  All the aforementioned reasons motivate us to use neural networks in this work.

\paragraph{Optimization method} Most of neural networks in literature, especially in computer science communities, are trained with stochastic gradient methods (SGD). In our present work and the previous work~\cite{huang2019predictive},
the optimization of the loss functions~\cref{EQ:LOSS_DIRECT,EQ:LOSS_INDIRECT} is done by the Limited-memory BFGS (L-BFGS-B) method~\cite{byrd1995limited} with the line search routine in~\cite{more1994line}, which attempts to enforce the Wolfe conditions~\cite{byrd1995limited}  by a sequence of polynomial interpolations. Note BFGS is applicable in our case since the data sets are typically small and the neural network is reasonably deep and wide; otherwise, the memory requirement of L-BFGS-B is so high that SGD or other similar first-order methods for training neural networks should be adopted. It is worth mentioning that the choice of BFGS optimizer is well-motivated and has long been adopted for scientific and engineering applications due to its fast convergence and robustness.

\paragraph{Data scaling} We observed that input data scaling significantly helps train SPD-NNs faster, reduces overfitting, and makes better predictions.
In the present work, the inputs and outputs of the neural network are scaled to a similar magnitude. 
Specifically, we introduce a strain reference~$\epsilon_\textrm{ref}$ and a stress reference~$\sigma_\textrm{ref}$ to scale strains and stresses. 
We also scale the tangent stiffness matrix by  $\frac{\sigma_\textrm{ref}}{\epsilon_\textrm{ref}}$. For example, in the elasto-plasticity case, the Cholesky factor has the form 
$$\ChoL_{\theta}\left(\frac{\bm\epsilon^{i}}{\epsilon_\textrm{ref}}, \frac{\bm\epsilon^{i-1}}{\epsilon_\textrm{ref}}, \frac{\bm\sigma^{i-1}}{\sigma_{\textrm{ref}}}\right)$$
The most important thing is that  $\frac{\sigma_\textrm{ref}}{\epsilon_\textrm{ref}} \sim O(E)$, here $E$ is the estimated Young's modulus. This guarantees the 
inputs for SPD-NN, especially for elasto-plasticitiy,  $(\bm{\epsilon}^{i}, \bm{\epsilon}^{i-1}, \bm{\sigma}^{i-1})$ have similar magnitudes.

\paragraph{Local minima} Local minima are observed in~\cref{SEC:TRUSS}, since training neural networks involves highly non-convex optimization problems. 
Although neural networks with minimal losses on the training set from 10 different initial weights perform well, neural networks with median losses are less satisfactory.
This reveals the uncertainty with respect to the initial weights for most NN-based data-driven approaches.
Wider neural networks will be considered in the future, since some theoretical and computational results~\cite{livni2014computational, kawaguchi2019effect} show quality of local minima 
tends to improve toward the global minimum value as depths and widths increase.
As for the indirect data training approach, pre-training approach~(see \cref{SEC:INDIRECT_DATA}) produces acceptable initial weights and thus relieves the concern. 
For all these training processes,  the optimization is terminated when the objective function is called 3000 times and 50000 times for pre-training and training, respectively.
The direct input-output data training is about 4 times faster than the indirect data training.

\section{Conclusion}
\label{SEC:CONCLUSION}

Data-driven approaches continue to gain popularity for constructing constitutive models from high-fidelity simulations and high-resolution experiments. The incorporation of data-driven constitutive models into conservation equations leads to a hybrid model, namely a coupled system with differential equations to describe conservation laws and neural networks to describe the material properties. 

To make these hybrid models numerically more robust, we introduce a novel neural network architecture, SPD-NN, where the neural network outputs the Cholesky factor of the tangent stiffness matrix instead of the stress or the stress increment. 
This neural network architecture weakly imposes convexity on the strain energy function~(i.e., SPD tangent stiffness matrix). The incremental form of SPD-NN also preserves the time consistency.
SPD-NN-based constitutive relations are tested on a 1D elasto-plastic truss problem and several 2D plate problems in which the plate is made of hyperelastic,  elasto-plastic, and multiscale fiber-reinforced materials.
When contrasting the SPD-NN with two other neural network architectures, we show that the SPD-NN  exhibits better numerical stability in the resulting hybrid models. 
The general training approach and the improved numerical stability allow for the potential to extend SPD-NNs to other time-dependent physical systems, such as viscoelastic materials, where the constitutive relations are rate-dependent, and plastic materials with stronger hysteresis, where more history-dependent variables are required. 

However, one limitation of the current approach is that the training process requires full field data, either for strain-stress pairs or displacement measurements. This may not be always possible; for example, the measurements may only be made on the surface of a 3D solid body. NN-based constitutive modeling with incomplete data remains to be investigated in the future.


\bibliographystyle{unsrt}
\bibliography{ref}
\end{document}